%% file: Diplomarbeit_revised.tex
\documentclass[11pt,titlepage,twoside,openright,a4paper]{scrbook}
\usepackage{a4}
\usepackage[latin1]{inputenc}
\usepackage[ngerman]{babel}
\usepackage{graphicx}
\usepackage{amsfonts}
\usepackage{amsmath}
\usepackage{amssymb}
\usepackage{amsthm}
\usepackage{cite}
\usepackage{fancyhdr}
\usepackage{nomencl}
\usepackage{longtable}
\usepackage{listings}
\usepackage{color}

\newtheorem{satz}{Satz}[chapter]
\newtheorem{lemma}[satz]{Lemma}
\newtheorem{satzdef}[satz]{Satz und Definition}
\newtheorem{prop}[satz]{Proposition}
\newtheorem{coro}[satz]{Korollar}
\theoremstyle{definition}
\newtheorem{definition}[satz]{Definition}
\newtheorem{bem}[satz]{Bemerkung}
\theoremstyle{remark}
\newenvironment{bew}{\begin{proof}[Beweis]}{\end{proof}}

\newcommand{\BIGOP}[1]{\mathop{\mathchoice
{\raise-0.22em\hbox{\huge $#1$}}
{\raise-0.05em\hbox{\Large $#1$}}{\hbox{\large $#1$}}{#1}}}
\newcommand{\bigtimes}{\BIGOP{\times}}

\newcommand{\E}{\mathbb{E}}
\newcommand{\Var}{\mbox{Var}}
\newcommand{\Cov}{\mbox{Cov}}
\newcommand{\Prob}{\mathbb{P}}
\newcommand{\R}{\mathbb{R}}
\newcommand{\Q}{\mathbb{Q}}
\newcommand{\Z}{\mathbb{Z}}
\newcommand{\N}{\mathbb{N}}
\newcommand{\argsup}{\mbox{argsup}}

\fancypagestyle{plain}{
  \fancyhf{}
  \cfoot{\thepage}
  
}

\begin{document}

\renewcommand{\baselinestretch}{1.25}

\include{Titelblatt}
\newpage\thispagestyle{empty}\hspace{1cm}\newpage

\pagenumbering{roman}
\renewcommand{\baselinestretch}{1.0}
\setlength{\parindent}{0pt} \setlength{\parskip}{1.5ex}
\include{Vorwort}
\include{Notation}

\tableofcontents
\newpage\thispagestyle{empty}\hspace{1cm}\newpage

\setcounter{secnumdepth}{3}
\setlength{\parindent}{0pt} \setlength{\parskip}{1.5ex}
\setcounter{tocdepth}{5}

\pagenumbering{arabic}

\include{Einleitung}
\include{Grundlagen}
\include{Fehlerabschaetzungen}
\include{Simulation_revised}


\begin{bibliography}{lit}
\addcontentsline{toc}{chapter}{Literaturverzeichnis}
\bibliographystyle{alphadin}
\end{bibliography}

\end{document}

%% file: Titelblatt.tex
\begin{titlepage}
\begin{center}
  \hspace{3cm}\\
  \vspace{1.2cm}
  \LARGE{\textbf{Simulationsverfahren f"ur Brown-Resnick-Prozesse}}\\
  \vspace{3.2cm}
\textsf{
    \begin{Large}\\
    \vspace{0.5cm}
    \textbf{Marco Oesting}\\
    \vspace{4cm}
    \textbf{Institut f"ur Mathematische Stochastik}\\
    \textbf{Georg-August-Universit"at zu G"ottingen}\\
    \vspace{0.5cm}
    \end{Large}
  }
\end{center}
\end{titlepage}

%% file: Vorwort.tex
\chapter*{Vorwort}

Im Rahmen der Vorlesung \glqq{}Einf"uhrung in die Extremwerttheorie\grqq{} von Prof. Dr. Martin Schlather lernte ich im Wintersemester 2008/2009 \glqq{}Brown-Resnick"=Prozesse\grqq{} als eine spezielle Klasse station"arer, max-stabiler stochastischer Prozesse kennen.

Diese haben vielf"altige Anwendungen. So wird z.\,B. in \cite{buishand-2008} ein Brown"=Resnick"=Prozess (mit Standard"=Fr\'echet"=R"andern) zum Variogramm $\gamma((h_1,h_2))= \frac 1 2 ( |h_1| + |h_2| )$ verwendet, um die r"aumliche Verteilung von extremen Regenf"allen zu modellieren. Aufgrund derartiger Anwendungen sind Brown-Resnick"=Prozesse und ihre Simulation auch von starkem wirtschaftlichen Interesse, wie z.\,B. f"ur auf Monte-Carlo-Methoden basierenden Kalkulationen von Versicherern.

Bei der Simulation von Brown-Resnick"=Prozessen nach dem sich aus der Definition kanonisch ergebenden Verfahren w"achst jedoch der Approximationsfehler exponentiell mit der Intervallgr"o"se. Daher verh"alt sich auch die Laufzeit entsprechend schlecht, wenn der Fehler auf gro"sen Intervallen klein gehalten werden soll.

Daraus ergeben sich die Problemstellung und das Ziel dieser Arbeit. Es sollen verschiedene Simulationsverfahren vorgestellt, untersucht und miteinander verglichen werden. Gesucht ist eine Anwort auf die Frage: Welches Verfahren ist -- je nach Anforderungen an die Simulation -- hinsichtlich des Fehlers \glqq{}optimal\grqq{}?
\medskip

Die Arbeit gliedert sich dabei wie folgt: In Kapitel 1 werden die Grundbegriffe f"ur die Definition von Brown-Resnick-Prozessen bereitgestellt. Diese erfolgt in Kapitel 2 ebenso wie die Vorstellung stochastisch "aquivalenter Darstellungen der Prozesse. Das Kapitel basiert im Wesentlichen auf dem Artikel \glqq{}Stationary Max-Stable Fields Associated to Negative Definite Functions\grqq{} von Sachar Kablutschko, Martin Schlather und Laurens de Haan (im Weiteren als \cite{kabluchko-2008} bezeichnet) sowie weiteren Ideen der Autoren.

Anhand der unterschiedlichen Darstellungen ergeben sich verschiedene Simulationsverfahren f"ur Brown-Resnick"=Prozesse. Diese werden in Kapitel 3 n"aher beschrieben und der bei den Verfahren auftretende Approximationsfehler wird abgesch"atzt. Schlie"slich werden die Ergebnisse der Simulationsverfahren in $\texttt{R}$ in Kapitel 4 dargestellt und verglichen. Programmcode, Dokumentation und eine detaillierte Auflistung der Ergebnisse der zugrundeliegenden Simulationsstudie finden sich im Anhang.
\medskip

Sehr herzlich danken m"ochte ich Prof. Martin Schlather f"ur meine Ausbildung in vielf"altigen Bereichen der Stochastik, insbesondere stochastischen Prozessen und Extremwerttheorie, f"ur die Vermittlung wesentlicher Programmierkenntnisse in \texttt{R}, die Bereitstellung des Packages \texttt{RandomFields} sowie f"ur die spannende Fragestellung und die hervorragende Betreuung der Arbeit.

Weiterhin gilt mein Dank Dr. Sachar Kablutschko f"ur seine Hilfe bei Messbarkeitsproblemen sowie seine Bereitschaft zur Zweitkorrektur der Arbeit.
\medskip

G"ottingen, im Juli 2009

%% file: Notation.tex
\chapter*{Notation}

\begin{tabular}{ll}
$\N$ & Menge der nat\"urlichen Zahlen $\N = \{1,\,2,\,\ldots\}$ \vspace{4mm}\\
{$\Z$, $\Q$, $\R$} & {Menge der ganzen/rationalen/reellen Zahlen} \vspace{4mm}\\
{$p\Z$} & {$\{pz: z \in \Z\}$} \vspace{4mm}\\
{$M_{\geq x}$} & {$\{y \in M : y \geq x\}$ mit $M \subset \R$} \vspace{4mm}\\
{$M_{> x}$} & {$\{y \in M : y > x\}$ mit $M \subset \R$} \vspace{4mm}\\
{$M^d$} & {$\{(x_1, \ldots, x_d): \ x_i \in M, \ i=1,\ldots, d\}$ f\"ur eine Menge $M$}\vspace{4mm}\\
{$\mathcal{B}$} & {Borel-$\sigma$-Algebra von $\R$}\vspace{4mm}\\
{$\mathcal{B}^d$} & {Borel-$\sigma$-Algebra des $\R^d$}\vspace{4mm}\\
{$\R^{m \times n}$} & {Menge aller $m \times n$-Matrizen mit reellen Eintr\"agen}\vspace{4mm}\\
{$C(\R^d)$} & {Menge aller stetigen Funktionen von $\R^d$ nach $\R$}\vspace{4mm}\\
{$\R^M$} & {Menge aller reellwertigen, auf einer Menge $M$ definierten Funktionen}\vspace{4mm}\\
{$\sigma(\tau)$} & {$\sigma$-Algebra, die von einem Mengensystem $\tau$ erzeugt wird}\vspace{4mm}\\
{$a\vee b$} & {das Maximum von $a$ und $b$}\vspace{4mm}\\
{$\Prob$} & {Wahrscheinlichkeitsma\ss}\vspace{4mm}\\
{$\E X$} & {Erwartungswert einer Zufallsvariablen $X$}\vspace{4mm}\\
{$\Var X$} & {Varianz einer Zufallsvariablen $X$}\vspace{4mm}\\
{$\Cov (X,Y)$} & {Kovarianz der Zufallsvariablen $X$ und $Y$}\vspace{4mm}\\
{$\delta_x$} & {Dirac-Ma\ss \ im Punkt $x$}\vspace{4mm}\\
{$\mathbf{1}_A$} & {Indikatorfunktion einer Menge $A$}\vspace{4mm}\\
{$\sim$} & {ist verteilt wie}\vspace{4mm}\\
{$u.i.v.$} & {unabh\"angig identisch verteilt}
\end{tabular}
\begin{tabular}{ll}
{$\stackrel{d}{=}$} & {Gleichheit in Verteilung}\vspace{4mm}\\
{$\mathcal{N}(\mu, \Sigma)$} & {Normalverteilung mit Erwartungswert $\mu$ und Kovarianzmatrix $\Sigma$}\vspace{4mm}\\
{$\mbox{Exp}(\lambda)$} & {Exponentialverteilung mit Parameter $\lambda$}\vspace{4mm}\\
{$U(M)$} & {Gleichverteilung auf einer Menge $M$}
\end{tabular}

\newpage

%% file: Einleitung.tex
\chapter{Einleitung}

\section{Zufallsfelder}

In diesem Kapitel sollen zun"achst die Begriffe, die f"ur die Definition von Brown-Resnick-Prozessen ben"otigt werden, erl"autert werden. Dabei handelt es sich u.\,a. um grundlegende m"ogliche Eigenschaften von Zufallsfeldern.

\begin{definition}
Sei $T$ eine Indexmenge. Eine Familie von reellwertigen Zufallsvariablen $\{ Z(t), \ t \in T \}$ hei"st \emph{stochastischer Prozess}. Im Fall $T=\R^d$ nennt man $\{ Z(t), \ t \in T \}$ auch \emph{$d$-dimensionales Zufallsfeld}.
\end{definition}

\begin{definition}
Ein Zufallsfeld hei"st \emph{Gau"ssch} genau dann, wenn alle endlich-dimensionalen Randverteilungen multivariat normal sind.
\end{definition}

\begin{definition}
Ein Zufallsfeld $Z(\cdot)$ auf $\R^d$ hei"st
\begin{itemize}
\item[1. ] \emph{station"ar (im engeren Sinne)} genau dann, wenn alle endlich"=dimensionalen Randverteilungen translationsinvariant sind, d.\,h.
\[ \Prob(Z(x_1+x) \in B_1, \ldots, Z(x_n+x) \in B_n) = \Prob(Z(x_1) \in B_1, \ldots, Z(x_n) \in B_n) \]
f"ur alle $x, x_i \in \R^d$, $B_i \in \mathcal{B}^d$, $i \in \{1,\ldots, n\}$, $n \in \N$.
\item[2. ] \emph{station"ar (im weiteren Sinne)} oder \emph{schwach station"ar} genau dann, wenn Erwartungswert $\E Z(t)$ und Kovarianz $\Cov(Z(x+t), Z(t))$ f"ur alle $x \in \R^d$ nicht von $t \in \R^d$ abh"angen.\\
Die Funktion $C: \R^d \to \R,$  \[ C(h) = \Cov(Z(h),Z(0))\] hei"st \emph{Kovarianzfunktion}.  
\end{itemize}
\end{definition}

\begin{bem}
Da eine multivariate Normalverteilung eindeutig durch Erwartungswert und Kovarianzmatrix bestimmt ist, sind in dem Fall, dass $Z(\cdot)$ ein Gau"ssches Zufallsfeld ist, die Begriffe \glqq{}station"ar im engeren Sinne\grqq{} und \glqq{}station"ar im weiteren Sinne\grqq{} "aquivalent.
\end{bem}

\begin{definition}
Ein Zufallsfeld $Z(\cdot)$ auf $\R^d$ hei"st \emph{intrinsisch station"ar} genau dann, wenn f"ur jedes $h \in \R^d$ das Zufallsfeld $\{Z(x+h) - Z(x), \ x \in \R^d\}$ schwach station"ar ist.\\
Die Funktion $\gamma: \R^d \to [0,\infty),$ \[\gamma(h) = \frac 1 2 \E(Z(h)-Z(0))^2\] hei"st \emph{(Semi-)Variogramm}.
\end{definition}

\begin{bem} \label{variostetig}
Aus der Stetigkeit des Variogramms im Ursprung folgt die Stetigkeit des Variogramms auf $\R^d$, denn es gilt f"ur beliebiges $x \in \R^d$:
\begin{align*}
& \lim_{h \to 0} \E( Z(x+h) - Z(0) )^2 - \E( Z(x) - Z(0) )^2\\
 ={} & \lim_{h \to 0} \E( Z(x+h) - Z(x) )^2 + 2\E(Z(x+h) - Z(x))(Z(x)-Z(0))\\
 & \hspace{1cm} +  \E( Z(x) - Z(0) )^2 -  \E( Z(x) - Z(0) )^2\\
 ={} & \lim_{h \to \infty} \E( Z(h) - Z(0) )^2 + 2\E(Z(x+h) - Z(x))(Z(x)-Z(0))\\
 \leq{} & \lim_{h \to \infty} \E( Z(h) - Z(0) )^2 + 2 \sqrt{\E(Z(x+h) - Z(x))^2 \cdot \E(Z(x) - Z(0))^2}\\
 & \hspace{2.5cm} \text{nach der Cauchy-Schwarz-Ungleichung}\\
 ={} & 0, \quad \text{falls } \lim_{h \to 0} \E(Z(h)-Z(0))^2 = 0.
\end{align*}
Somit ist ein Variogramm $\gamma$ genau dann stetig auf $\R^d$, wenn es im Ursprung stetig ist.
\end{bem}

\begin{definition} \label{max-stabil}
Sei $\{ Z(t), \ t \in \R^d \}$ ein Zufallsfeld auf $\R^d$ mit
Standard-Gumbel-R"andern, d.\,h. es ist $\Prob(Z(t) \leq x) = \exp(-e^{-x})$ f"ur alle $t \in \R^d, \ x \in \R,$ und seien $Z_i(\cdot) \sim_{u.i.v.} Z(\cdot), \ i \in \N$.\\
Dann hei"st $Z(\cdot)$ \emph{max-stabil} genau dann, wenn f"ur alle $n \in \N$ die Beziehung
\[ \max_{i=1,\ldots,n} Z_i(\cdot) - \log{n} \stackrel{d}{=} Z(\cdot) \]
gilt.
\end{definition}
\bigskip

\section{Poisson-Punktprozesse}

Weiterhin treten im Zusammenhang mit Brown-Resnick-Prozessen sogenannte \glqq{}Poisson"=Punktprozesse\grqq{} auf, die auf \glqq{}Polnischen R"aumen\grqq{} definiert sind.

\begin{definition}
Sei $(\Omega, m)$ ein metrischer Raum und $\mathcal{T}$ die von $m$ induzierte Topologie. Dann hei"st $\Omega$ \emph{Polnischer Raum} genau dann, wenn $\Omega$ bez"uglich $m$ vollst"andig ist und $\mathcal{T}$ eine abz"ahlbare Basis besitzt.

Die von der Topologie $\mathcal{T}$ erzeugte $\sigma$-Algebra, $\sigma(\mathcal{T})$, hei"st \emph{Borel-$\sigma$-Algebra}.
\end{definition}

\begin{definition}
Sei $S$ ein lokal kompakter Polnischer Raum mit Borel-$\sigma$-Algebra $\mathcal{S}$, versehen mit einem lokal endlichen Ma"s $\mu$. Ein zuf"alliges Z"ahlma"s $\Pi$ auf $(S, \mathcal{S})$ hei"st \emph{Poisson"=Punktprozess mit Intensit"atsma"s $\mu$} genau dann, wenn
\begin{itemize}
\item[1. ] $\Pi(A_1), \ldots, \Pi(A_n)$ sind stochastisch unabh"angig f"ur alle $A_1, \ldots, A_n \in \mathcal{S}$ mit  $A_i \cap A_j = \emptyset$ f"ur $i \neq j$, $i,j \in \{1, \ldots, n\}$ und $n \in \N$,
\item[2. ] $\Pi(A)$ ist Poisson-verteilt mit Parameter $\mu(A)$ f"ur alle $A \in \mathcal{S}$.
\end{itemize}
\end{definition}
\medskip

Das nachfolgende Lemma behandelt die \glqq{}Erweiterung\grqq{} von bestimmten Poisson"=Punktprozessen um Zufallsvariablen. Eine solche wird im weiteren Verlauf der Arbeit noch mehrfach  auftreten.

\begin{lemma} \label{2dim-ppp}
Sei $\sum_{i \in \N} \delta_{X_i}$ ein Poisson"=Punktprozess auf \ $\R$ mit Intensit"atsma"s $\lambda e^{-x} dx$, $\lambda > 0$. Seien weiterhin $Y_1, Y_2, \ldots$ unabh"angig identisch verteilte Zufallsgr"o"sen (unabh"angig von $\{X_i\}_{i \in \N}$) mit Werten in einem Polnischen Raum $S$. Sei $\Prob$ das zu $Y_i$ geh"orige Wahrscheinlichkeitsma"s.\\
Dann ist $\widehat{\Pi} := \sum_{i \in \N} \delta_{(X_i, Y_i)}$ ein Poisson"=Punktprozess auf $\ \R \times S$ mit Intensit"atsma"s $\lambda e^{-x} dx \Prob(dy)$.
\end{lemma}
\begin{bew}
Zun"achst bemerken wir, dass das Ma"s $\lambda e^{-x} dx \Prob(dy)$ lokal endlich ist, denn f"ur jede beschr"ankte Borel-messbare Menge $A \times B \subset \R \times S$ ist $$\int_{A \times B} \lambda e^{-x} dx \Prob(dy) \leq \int_{A \times S} \lambda e^{-x} dx \Prob(dy) = \int_{A} \lambda e^{-x} dx < \infty,$$ da $A$ beschr"ankt ist.

Seien nun $[a,b] \subset \R$, $A \subset S$ Borelmengen. Dann gilt f"ur $k \in \N$:
\begin{align*}
& \Prob(\widehat{\Pi}([a,b] \times A)=k)\\
 = {}& \sum_{n=k}^\infty \Prob(\Pi([a,b]=n) \Prob(\#\{i \in \{1,\ldots,n\}: Y_i \in A \}=k) \\
 = {} & \sum_{n=k}^\infty \exp(-\lambda(e^{-a}-e^{-b})) \frac{1}{n!} (\lambda(e^{-a}-e^{-b}))^n \binom n k \Prob(A)^k (1-\Prob(A))^{n-k} \\
 = {} & \frac{1}{k!} \exp(-\lambda(e^{-a}-e^{-b})) \left(\Prob(A)\lambda(e^{-a}-e^{-b})\right)^k \\
 & \hspace{6cm} \cdot \sum_{n=0}^\infty \frac{(\lambda(e^{-a}-e^{-b}))^n}{n!} \left(1-\Prob(A)\right)^n \\
 = {} & \frac{1}{k!} \exp\left(-\Prob(A)\lambda(e^{-a}-e^{-b})\right) \left(\Prob(A)\lambda(e^{-a}-e^{-b})\right)^k \\
 = {} & \frac{1}{k!} \exp\left(-\int_{A} \lambda \int_{[a,b]} e^{-x} dx d\Prob(y)\right) \left(\int_{A} \lambda \int_{[a,b]} e^{-x} dx d\Prob(y)\right)^k
\end{align*}
Die Anzahl von Punkten in einer Borel-Menge ist also Poisson-verteilt mit der behaupteten Intensit"at.
\medskip

Die Unabh"angigkeit der Anzahlen von Punkten in  Mengen, die in der ersten Komponente disjunkt sind, folgt bereits aus der Tatsache, dass $\sum_{i \in \N} \delta_{X_i}$ ein Poisson"=Punktprozess ist und dass $Y_1, Y_2, \ldots$ unabh"angig sind.\\
Es gen"ugt daher, dies f"ur Mengen der Form $[a,b] \times A_1, [a,b] \times A_2, \ldots, [a,b] \times A_n$ zu zeigen mit disjunkten Borel-Mengen $A_1, A_2, \ldots, A_n \subset S$ und $[a,b] \subset \R$. Seien dazu $k_1, \ldots, k_n \in \N, \ k_1 + \ldots + k_n = k$ und $\Lambda_{a,b} = \lambda(e^{-a}- e^{-b})$.
Dann gilt
\begin{align*}
&\Prob(\widehat{\Pi}(A_1 \times [a,b]) = k_1, \ldots, \widehat{\Pi}(A_n \times [a,b]) = k_n)\\
= {} & \sum_{j=k}^\infty \frac 1 {j!} \Lambda_{a,b}^j \exp(-\Lambda_{a,b}) \frac {j!} {(j-k)!} \left(1 - \sum_{i=1}^n \Prob(A_i)\right)^{j-k} \prod_{i=1}^n \frac 1 {k_i!} \left(\Prob(A_i)\right)^{k_i} \\
= {} & \prod_{i=1}^n \frac 1 {k_i!} \left( \Prob(A_i)\Lambda_{a,b} \right)^{k_i} \exp\left(-\Lambda_{a,b} \Prob(A_i)\right)\\
= {} & \prod_{i=1}^n \Prob(\widehat{\Pi}(A_i \times [a,b]) = k_i),
\end{align*}
also die zu zeigende stochastische Unabh"angigkeit.

Somit ist $\widehat{\Pi}$ ein Poisson-Punktprozess mit dem angegebenen Intensit"atsma"s.
\end{bew}

Ebenfalls werden wir den folgenden Transformationssatz f"ur Poisson"=Punktprozesse ben"otigen.

\begin{satz} \label{ppp-trafo}
Seien $S$, $T$ Polnische R"aume mit Borel-$\sigma$-Algebren $\mathcal{S}$ bzw. $\mathcal{T}$ und $\Pi = \sum_{i \in \N} \delta_{X_i}$ ein Poisson"=Punktprozess auf $(S, \mathcal{S})$ mit lokal endlichem Intensit"atsma"s $\mu$ und $f: S \to T$ messbar so, dass auch das Bildma"s $\mu^* = \mu \circ f^{-1}$ lokal endlich ist. Dann ist durch $\sum_{i \in \N} \delta_{f(X_i)}$ ein Poisson"=Punktprozess auf $(T, \mathcal{T})$ mit Intensit"atsma"s $\mu^*$ gegeben. 
\end{satz}
\begin{bew}
siehe \cite{klenke-2006}, Satz 24.16.
\end{bew}

%% file: Grundlagen.tex
\chapter{Brown-Resnick-Prozesse}

Bevor Brown-Resnick-Prozesse definiert werden k"onnen, sind zun"achst noch einige Vorbemerkungen notwendig: In diesem Kapitel werden wir u.\,a. Poisson"=Punktprozesse auf dem Raum $C(\R^d)$, also dem Raum der reellwertigen, stetigen Funktionen auf $\R^d$, betrachten. Dazu versehen wir $C(\R^d)$  mit der $\sigma$-Algebra $\mathcal{C}$, die von den Mengen der Form
\[ C_{t_1,\ldots,t_m}(B) := \{ f \in C(\R^d): (f(t_1),\ldots,f(t_m)) \in B \}, \]
mit $t_1, \ldots, t_m \in \R^d$, $m \in \N$ und $B \in \mathcal{B}^m$ erzeugt wird.
Diese $\sigma$-Algebra erm"oglicht es uns, Aussagen "uber Funktionswerte an endlich vielen Stellen zu treffen, was sich als sehr n"utzlich herausstellen wird.
\medskip

Nun ist $\Omega_k := C([-k,k]^d)$, versehen mit der Metrik $m_k$ gegeben durch $$m_k(f,g):= \sup_{x \in [-k,k]^d} |f(x)-g(x)|,$$ ein Polnischer Raum. Dabei ist durch $$\mathcal{F}_k := \{ B_{\varepsilon}^{(k)}(f): \ \varepsilon \in \Q_{>0}, f \in \Pi_{\Q} \},$$ wobei $\Pi_\Q$ die Menge aller Polynome auf $\R^d$ mit rationalen Koeffizienten und $B_{\varepsilon}^{(k)}(f) := \{g \in \Omega_k: \ m_k(f,g) < \varepsilon \}$ sei, nach dem Satz von Stone-Weierstra"s (reelle Fassung) eine abz"ahlbare topologische Basis gegeben.

Damit ist auch der Produktraum $\Omega:=\bigtimes_{k=1}^{\infty} \Omega_k$ mit der Metrik $$m := \sum_{k=1}^\infty 2^{-k} \frac {m_k}{1+m_k}$$ ein Polnischer Raum.
\medskip

\begin{lemma} \label{prodtopo}
Die von $m$ auf $\Omega$ induzierte Topologie ist die Produkttopologie aus den von den Metriken $m_k$ auf $\Omega_k$, $k \in \N$, induzierten Topologien.
\end{lemma}
\begin{bew}
Es gen"ugt zu zeigen, dass die Basis der einen Topologie jeweils in der anderen Topologie enthalten ist.

Eine Basis der Produkttopologie bilden die Mengen der Form $\bigtimes_{k \in \N} B_{\varepsilon_k}^{(k)}(f_k)$ mit $f_k \in \Omega_k$, wobei nur endlich viele $\varepsilon_k$ endlich seien und wir $B_{\infty}^{(k)}(f) := \Omega_k$ definieren.
F"ur $g \in \bigtimes_{k \in \N} B_{\varepsilon_k}^{(k)}(f_k)$ mit $g = (g_1, g_2, \ldots)$ gibt es nun $\delta = (\delta_1, \delta_2, \ldots)$ mit $0 < \delta_k < \infty$, falls $\varepsilon_k < \infty$, und $\delta_k = \infty$, falls $\varepsilon_k = \infty$, so, dass $$\bigtimes_{k \in \N} B_{\delta_k}^{(k)}(g_k) \subset \bigtimes_{k \in \N} B_{\varepsilon_k}^{(k)}(f_k).$$
Sei $\varepsilon := \min_{\{k \in \N: \ \delta_k < \infty\}} \left(2^{-k-1}\cdot \min(\delta_k,1) \right)$, dann folgt -- wegen $\frac{1}{2^k} \frac{m_k(\cdot,\cdot)}{1+m_k(\cdot,\cdot)} \leq m(\cdot,\cdot)$ -- f"ur $h=(h_1, h_2, \ldots) \in B_{\varepsilon}(g) := \{ h \in \Omega: \ m(h,g) < \varepsilon \}$, dass
$$m_k(g_k,h_k) \leq \frac{m(g,h)2^k}{1 - m(g,h)2^k} \leq \frac{m(g,h)2^k}{\frac 1 2} \leq \delta_k,$$
 falls $\delta_k < \infty$.
Somit ist $$B_{\varepsilon}(g) \subset \bigtimes_{k \in \N} B_{\delta_k}^{(k)}(g_k) \subset \bigtimes_{k \in \N} B_{\varepsilon_k}^{(k)}(f_k),$$ d.\,h. $\bigtimes_{k \in \N} B_{\varepsilon_k}^{(k)}(f_k)$ ist auch offen bez"uglich der von $m$ erzeugten Topologie.

Sei andererseits nun ein Basiselement $B_\varepsilon(f)$ der von $m$ induzierten Topologie gegeben. Dann existiert f"ur jedes $g \in B_\varepsilon(f)$ ein $\delta > 0$ mit $B_{\delta}(g) \subset B_\varepsilon(f)$. W"ahle nun $N \in \N$ so gro"s, dass $2^{-N} < \frac{\delta} 2$.

Dann gilt f"ur jedes $h \in B_{\delta/2}^{(1)}(g_1) \times B_{\delta/2}^{(2)}(g_2) \times \ldots \times B_{\delta/2}^{(N)}(g_N) \times \bigtimes_{k=N+1}^{\infty} \Omega_{k}$, dass
$$m(g,h) \leq \sum_{k=1}^{N} 2^{-k} m_k(g_k,h_k) + \sum_{j=N+1}^{\infty} 2^{-j} \leq \sum_{k=1}^{N} 2^{-k} \frac {\delta} 2 + 2^{-N} < \frac{\delta} 2 + \frac{\delta} 2 = \delta.$$
Also ist $B_{\delta/2}^{(1)}(g_1) \times B_{\delta/2}^{(2)}(g_2) \times \ldots \times B_{\delta/2}^{(N)}(g_N) \times \bigtimes_{k=N+1}^{\infty} \Omega_{k} \subset B_{\delta}(g) \subset B_\varepsilon(f)$, d.\,h. jedes Element $g \in B_\varepsilon(f)$ hat eine Umgebung aus der Basis der Produkttopologie, die in $B_\varepsilon(f)$ enthalten ist. Somit l"asst sich $B_\varepsilon(f)$ als Vereinigung von Basiselementen der Produkttopologie schreiben und ist daher in dieser enthalten.
\end{bew}

Da wir $C(\R^d)$ mit der abgeschlossenen Teilmenge $$F := \{ (f_1, f_2, \ldots) \in \Omega: \ f_{n+1} \big|_{[-n, n]^d} = f_n \ \forall n \in \N\} \subset \Omega$$ identifizieren k"onnen, ist auch $C(\R^d)$ mit der Metrik $m$ ein Polnischer Raum.
Nun ist der Zusammenhang zwischen der Borel-$\sigma$-Algebra von $C(\R^d)$ und der $\sigma$-Algebra $\mathcal{C}$ zu kl"aren.

\begin{prop}
Sei $\tau$ die von $m$ auf $C(\R^d)$ induzierte Topolgie. Dann gilt $\sigma(\tau) = \mathcal{C}$.
\end{prop}
\begin{bew}
Wir stellen fest, dass $\mathcal{C}$ bereits von den Mengen der Form $C_{t_1, \ldots, t_m}(B)$, wobei $t_1, \ldots, t_m \in \R^d$, $m \in \N$ und $B \subset \R^m$ offen seien, erzeugt wird. Da $B$ offen ist, existiert zu jedem $f \in C_{t_1, \ldots, t_m}(B)$ ein $\varepsilon > 0$ mit $$C_{t_1, \ldots, t_m}\left((f(t_1)-\varepsilon, f(t_1) + \varepsilon) \times \ldots \times (f(t_m)-\varepsilon, f(t_m) + \varepsilon)\right) \subset C_{t_1, \ldots, t_m}(B),$$ denn diese Teilmenge ist ein offener Ball um $f$ in der Supremumsnorm auf $\R^m$.

W"ahlen wir nun $N \in \N$ derart, dass $t_1, \ldots, t_m \in [-N, N]^d$ ist, so gilt nach Konstruktion $$\left(\left(\bigtimes_{k=1}^N B_{\varepsilon}^{k}(g) \times \bigtimes_{k=N+1}^{\infty} \Omega_k \right) \cap F\right) \subset C_{t_1, \ldots, t_m}(B).$$
Also gibt es zu jedem Element aus $C_{t_1, \ldots, t_m}(B)$ eine in $C_{t_1, \ldots, t_m}(B)$ enthaltene offene Umgebung aus der Basis der Spur der Produkttopologie auf $F$. Diese Topologie ist nach Lemma \ref{prodtopo} gerade $\tau$. Also ist jede erzeugende Menge $C_{t_1, \ldots, t_m}(B) \in \tau$ und damit gilt auch $\mathcal{C} \subset \sigma(\tau)$.
\medskip

Da 
\begin{align*}
\mathcal{F} :={} & \Bigg\{ \left(\bigtimes_{k \in \N} B_{\varepsilon_k}^{(k)}(f_k)\right) \cap F : \ \varepsilon_k \in \Q_{>0} \cup \{\infty\}, \ f_k \in \Pi_{\Q},\\ & \hspace{2cm} \text{wobei nur endlich viele $\varepsilon_k$ endlich seien}  \Bigg\}
\end{align*}
 nach obigen Ausf"uhrungen und Lemma \ref{prodtopo} eine abz"ahlbare Basis von $\tau$ ist, gilt $\sigma(\mathcal{F}) = \sigma(\tau)$ und es gen"ugt $G \in \mathcal{C}$ f"ur alle $G \in \mathcal{F}$ zu zeigen.\\
Sei nun $q_1^{(k)}, q_2^{(k)}, \ldots$ eine Abz"ahlung von $\Q^d \cap [-k,k]^d$. Dann ist
\begin{align*}
&\left(\bigtimes_{k \in \N} B_{\varepsilon_k}^{(k)}(f_k)\right) \cap F\\
 ={}& \left\{ g \in C(\R^d) \ \big| \ \forall k \in \N \ \exists n_k \in \N \ \forall j \in \N :  \ | g(q_j^{(k)}) - f_k(q_j^{(k)}) | < \varepsilon_k - \frac 1 {n_k} \right\}\\
 ={} & \bigcap_{k \in \N} \bigcup_{n_k \in \N} \bigcap_{j \in \N} C_{q_j^{(k)}}\left( \left(f_k(q_j^{(k)}) - \varepsilon_k + \frac 1 {n_k}, f_k(q_j^{(k)}) + \varepsilon_k-\frac 1 {n_k}\right) \right) \in \mathcal{C}
\end{align*}
 -- wie behauptet.
\end{bew}

Somit ist $C(\R^d)$ ein Polnischer Raum mit Borel-$\sigma$-Algebra $\mathcal{C}$. Da $C(\R^d)$ zudem lokal kompakt ist, k"onnen wir im Folgenden auf $(C(\R^d), \mathcal{C})$ Poisson"=Punktprozesse definieren.
\newpage

\section{Brown-Resnick-Prozesse als max-stabile Prozesse}

F"ur den Beweis des Satzes \ref{brown-resnick}, der gleichzeitig den Begriff \emph{Brown-Resnick-Prozess} definiert, ben"otigen wir das folgende Lemma. Dieses l"asst sich leicht nachrechnen, weshalb hier auf seinen Beweis verzichtet werden soll.

\begin{lemma} \label{nvbedingt}
Seien $$(X_1, X_2) \sim {\cal N} \left( \left(\begin{array}{*{1}{c}} \mu_1 \\ \mu_2 \end{array}\right),
\left( \begin{array}{*{2}{c}} \Sigma_{11} & \Sigma_{12} \\ \Sigma_{21} & \Sigma_{22} \end{array} \right) \right)$$ mit
$\mu_1 \in \R^k, \mu_2 \in \R^l, \Sigma_{11} \in \R^{k\times k}, \Sigma_{12} = \Sigma_{21}^T \in R^{k \times l}, \Sigma_{22} \in \R^{l \times l}$.\\
Dann gilt:
\[ X_2 | (X_1 = x_1) \sim {\cal N} ( \mu_2 + \Sigma_{21} \Sigma_{11}^{-1} (x_1 - \mu_1), \Sigma_{22}- \Sigma_{21} \Sigma_{11}^{-1} \Sigma_{12}). \]
\end{lemma}
\bigskip

Der nachfolgende Satz geht auf Brown und Resnick zur"uck, die ihn f"ur die Brownsche Bewegung als zugrundeliegendem Prozess (d.\,h. f"ur den Spezialfall $\gamma(h) = \frac 1 2 |h|$) bewiesen (vgl. \cite{brown-1977}).

\begin{satzdef} \label{brown-resnick}
Sei $W(\cdot)$ ein pfadstetiges, intrinsisch station"ares Gau"ssches Zufallsfeld auf \ $\R^d$ mit
Variogramm $\gamma$, Erwartungswert $\E W(t) = 0$  und Varianz \ $\sigma^2(t) := \emph{\Var}(W(t))$, $t \in \R^d$. \\
Weiterhin seien $W_i(\cdot) \sim W(\cdot),\, i \in \N$ unabh"angig identisch verteilt, sowie\\
$\Pi = \sum \delta_{X_i}$ ein von diesen stochastisch unabh"angiger Poisson-Punktprozess auf \ $\R$ mit Intensit"atsma"s $e^{-x}dx$. Dann ist durch
\[ Z(t) := \max_{i \in \N} \left(X_i + W_i(t) - \frac{\sigma^2(t)}{2}\right), \quad t \in \R^d \]
ein station"arer max-stabiler Prozess auf \ $\R^d$ mit Standard-Gumbel-R"andern gegeben. Die Verteilung von $Z$ h"angt nur von $\gamma$ ab.\\
$Z(\cdot)$ hei"st \emph{Brown-Resnick-Prozess} zum Variogramm $\gamma$.
\end{satzdef}
\begin{bew}
Seien $t_1,\,\ldots,\,t_m \in \R^d, \ y_1,\, \ldots y_m \in \R, \ m \in \N$ beliebig und $\Prob_{\widetilde{t_1},\ldots,\widetilde{t_k}}$ sei das zum Zufallsvektor $(W(\widetilde{t_1}),\ldots,W(\widetilde{t_k}))$ geh"orende Wahrscheinlichkeitsma"s.\\
Wir definieren nun f"ur $i=1,\ldots,m$ die Gr"o"sen $\sigma^2_i := \sigma^2(t_i)$ sowie die Zufallsvariablen $\xi_i := W(t_i) - \frac{\sigma^2_i}{2}$.
Dann ist $\xi_i \sim {\cal N} ( -\frac{\sigma_i^2}{2}, \sigma_i^2)$.
\medskip

Sei 
\begin{align}
\Gamma_{ij} :={}& \Cov(\xi_i - \xi_1, \xi_j - \xi_1) = \E( (W(t_i)-W(t_1))\cdot (W(t_j)-W(t_1)) ) \nonumber \\
={}&-\frac{1}{2} \E \{((W(t_i)-W(t_1))-(W(t_j)-W(t_1)))^2 \nonumber \\
&  \qquad - (W(t_i)-W(t_1))^2 - (W(t_j)-W(t_1))^2 \} \nonumber \\
={}&-\gamma(t_i-t_j) + \gamma(t_i-t_1) + \gamma(t_j-t_1)	\label{eq:Cov}
\end{align}
 und $\Gamma := \left( \Gamma_{ij} \right)_{i=2,\ldots,m, \ j=2,\ldots,m}$.
\medskip

Weiterhin seien $\mu_{2m} := (\E(\xi_j-\xi_1))_{j=2,\ldots,m} = \left(\frac{\sigma_1^2 - \sigma_j^2} 2 \right)_{j=2,\ldots,m}$\\
 und $g:=(\Cov(\xi_j-\xi_1,\xi_1))_{j=1,\ldots,m}$.
 
Nun ist 
\[\left(\begin{array}{*{1}{c}} \xi_2-\xi_1\\ \vdots\\ \xi_m -\xi_1\\ \xi_1 \end{array} \right)= \left( \begin{array}{*{5}{c}}
 -1  & 1 & 0 & \cdots & 0 \\
 \vdots & & \ddots & \ddots \\
 -1 & 0 & \cdots &0 &1 \\
 1 & 0 & \cdots & 0 &0 \end{array} \right)
 \cdot 
\left(\begin{array}{*{1}{c}} \xi_1 \\ \xi_2\\ \vdots \\ \xi_m \end{array} \right) \]
und $(\xi_1, \ldots, \xi_m)$ ist normalverteilt.
\medskip

Daher ist auch $(\xi_2 - \xi_1, \ldots, \xi_m - \xi_1, \xi_1)$ normalverteilt und nach obigen "Uberlegungen gilt
\begin{equation*}
\left(\begin{array}{*{1}{c}} \xi_2-\xi_1\\ \vdots\\ \xi_m -\xi_1\\ \xi_1 \end{array} \right) \sim {\cal N} \left(
\left( \begin{array}{*{1}{c}} \mu_{2m} \\ -\frac{\sigma_1^2} 2 \end{array} \right), \left( \begin{array}{*{2}{c}} \Gamma & g \\
 g^T & \sigma_1^2 \end{array} \right) \right).
\end{equation*}
Dann ist nach Lemma \ref{nvbedingt}
\begin{equation}
\xi_1 \ | \ (\xi_j - \xi_1)_{j=2,\ldots, m} \sim {\cal N} \left( -\frac{\sigma_1^2}{2} + g^T \Gamma^{-1} ((\xi_j - \xi_1)_{j=2,\ldots, m} - \mu_{2m}),
\sigma_1^2 - g^T \Gamma^{-1} g \right) \label{eq:bedingt}.
\end{equation}
Weiterhin ist $\Phi := \sum_{i=1}^n \delta_{(X_i, W_i(\cdot))}$ nach Lemma \ref{2dim-ppp} ein Poisson"=Punktprozess auf $(\R \times C(\R^d), \mathcal{B} \times \mathcal{C})$ mit Intensit"atsma"s $e^{-x}dx \times \Prob$ und unter Verwendung der "Aquivalenz
\begin{align*}
& Z(t_1) \leq y_1, \ldots, Z(t_m) \leq y_m \\
\iff {} &  X_i \leq \min_{k=1,\ldots, m} \left(y_k - W_i(t_k) + \frac{\sigma^2(t_k)}{2}\right) \quad \forall i \in \N \\
\iff {} &  \Phi \left( \left\{ (x,w): x \geq \min_{k=1,\ldots, m} \left(y_k - w(t_k) + \sigma^2(t_k)/2 \right)\right\} \right) = 0
\end{align*}
folgt daraus
\begin{align*}
&-\log(\Prob(Z(t_1) \leq y_1, \ldots, Z(t_m) \leq y_m)) \\
={}& \int_{\R^d} \int_{\min_{k=1,\ldots, m} \left(y_k - w_k + \sigma^2(t_k)/2\right)}^\infty e^{-x}dx \ 
\Prob_{t_1,\ldots,t_m}(d(w_1,\ldots,w_m)) \\
={}& \int_{\R^d} \max_{k=1,\ldots, m} e^{-(y_k - w_k + \sigma^2(t_k)/2)} \Prob_{t_1,\ldots,t_m}(d(w_1,\ldots,w_m))
\end{align*}
\begin{align*}
={}& \E \left(\max_{k=1,\ldots, m} e^{-(y_k - \xi_k)} \right) = \E \left( e^{\xi_1} \max_{k=1,\ldots, m} e^{(\xi_k-\xi_1)-y_k} \right)\\
={}&\E \left( \E\left( e^{\xi_1} \max_{k=1,\ldots, m} e^{(\xi_k-\xi_1)-y_k} \ \big| \ (\xi_j - \xi_1)_{j=2,\ldots,m} \right) \right) \\
={}& \E \left(\max_{k=1,\ldots, m} e^{(\xi_k-\xi_1)-y_k} \E\left( e^{\xi_1}\ \big| \ (\xi_j - \xi_1)_{j} \right)\right) \\
\stackrel{\eqref{eq:bedingt}}{=}&  \E \left(\max_{k=1,\ldots, m} e^{(\xi_k-\xi_1)-y_k} \cdot \exp \left(g^T \Gamma^{-1} ((\xi_j - \xi_1)_j-\mu_{2m}) - \frac{1}{2} g^T \Gamma^{-1} g \right) \right),\\
& \hspace{5.5cm} \text{da $U \sim {\cal N}(\mu, \sigma^2) \Rightarrow \E(e^{U}) = e^{\mu + \frac{\sigma^2}{2}}$}\\
={}& \int_{\R^{m-1}} \max_{k=2,\ldots, m} \{ e^{-y_1}, e^{w_k-y_k}\} \cdot \exp \left(g^T \Gamma^{-1} (w -\mu_{2m}) - \frac{1}{2} g^T \Gamma^{-1} g \right) \\ 
& \quad \cdot \frac{1}{(2\pi)^{\frac{m-1}2} \det(\Gamma)} \exp \left(-\frac{1}{2} (w-\mu_{2m})^T\Gamma^{-1} (w-\mu_{2m}) \right) \ dw_2 \cdots dw_{m}\\
={}& \int_{\R^{m-1}} \max_{k=2,\ldots, m} \{ e^{-y_1}, e^{w_k-y_k}\} \cdot \frac{1}{(2\pi)^{\frac{m-1}2} \det(\Gamma)} \\
& \quad \cdot \exp \left(-\frac{1}{2} (w-\mu_{2m}-g)^T\Gamma^{-1} (w-\mu_{2m}-g)\right) \ dw_2 \cdots dw_{m}.
\end{align*}

Mit 
\begin{align*}
 &(\mu_{2m} + g)_j = \frac{1}{2}(\sigma_1^2-\sigma_j^2) + \E((\xi_j - \E\xi_j - \xi_1 + \E\xi_1)(\xi_1-\E\xi_1)) \\
 ={}&-\frac{1}{2}\E \left( -(\xi_1 - \E\xi_1)^2 + (\xi_j - \E\xi_j)^2 - 2(\xi_j - \E\xi_j)(\xi_1-\E\xi_1) + 2(\xi_1-\E\xi_1)^2\right)\\
 ={}&-\frac{1}{2}\E(\xi_j - \E\xi_j - \xi_1 + \E\xi_1)^2 = -\frac{1}{2}\E(W(t_j) - W(t_1))^2 = -\gamma(t_j-t_1)
\end{align*}
erh"alt man also
\begin{align*}
&-\log(\Prob(Z(t_1) \leq y_1, \ldots, Z(t_m) \leq y_m)) \\
={}&\int_{\R^{m-1}} \max_{k=2,\ldots, m} \{ e^{-y_1}, e^{w_k-y_k}\} \cdot \frac{1}{(2\pi)^{(m-1)/2} \det(\Gamma)} \\
& \quad \cdot \exp \left(-\frac{1}{2} (w+\gamma(t_\cdot-t_1))^T\Gamma^{-1} (w-\gamma(t_\cdot-t_1))\right) \ dw_2 \cdots dw_{m}.
\end{align*}

In dieser Gleichung treten nach \eqref{eq:Cov} nur Terme der Form $\gamma(t_j-t_i)_{i,j=1,\ldots,m}$ auf, d.\,h. $Z(\cdot)$ ist station"ar und die Verteilung h"angt nur von $\gamma$ ab.

Im Fall $m=1$ ist das Integral degeneriert und wir erhalten 
\[-\log(\Prob(Z(t)\leq y) = e^{-y},\]
d.\,h. $Z(\cdot)$ hat Standard-Gumbel-R"ander.

Seien nun $Z^{(l)}(\cdot) \sim_{u.i.v.} Z(\cdot)$ f"ur $l \in \{1,\ldots, n\}$. Dann gilt
\begin{align*}
&-\log(\Prob(\max_{l=1,\ldots,n} Z^{(l)}(t_1) - \log n \leq y_1, \ldots,\max_{l=1,\ldots,n}  Z^{(l)}(t_m) - \log n \leq y_m)) \\
={}&-\log(\Prob(Z(t_1) \leq y_1 + \log n, \ldots, Z(t_m) \leq y_m+ \log n)) \cdot n \\
={}&\int_{\R^{m-1}} \max_{k=2,\ldots, m} \{ e^{-y_1-\log n}, e^{w_k-y_k- \log n}\} \cdot \frac{1}{(2\pi)^{(m-1)/2} \det(\Gamma)} \\
& \quad \cdot \exp \left(-\frac{1}{2} (w+\gamma(t_\cdot-t_1))^T\Gamma^{-1} (w-\gamma(t_\cdot-t_1))\right) \ dw_2 \cdots dw_{m} \cdot n\\
={}&\frac{1}{n} \int_{\R^{m-1}} \max_{k=2,\ldots, m} \{ e^{-y_1}, e^{w_k-y_k}\} \cdot \frac{1}{(2\pi)^{(m-1)/2} \det(\Gamma)} \\
& \quad \cdot \exp \left(-\frac{1}{2} (w+\gamma(t_\cdot-t_1))^T\Gamma^{-1} (w-\gamma(t_\cdot-t_1))\right) \ dw_2 \cdots dw_{m} \cdot n\\
={}&-\log(\Prob(Z(t_1) \leq y_1, \ldots, Z(t_m) \leq y_m)),
\end{align*}
also $\max_{l=1,\ldots,n} Z^{(l)}(\cdot) - \log n \stackrel{d}{=} Z(\cdot)$ -- wie behauptet.
\end{bew}
\bigskip

\begin{bem}\hspace{8cm}
\begin{itemize}
\item[1. ] Im Folgenden wird mit $Z(\cdot)$ stets der oben definierte Brown-Resnick-Prozess bezeichnet.
\item[2. ] Sei $W(\cdot)$ ein intrinsisch-station"ares Gau"ssches Zufallsfeld auf $\R^d$ mit Erwartungswert 0 und Variogramm $\gamma$. Existieren $ \alpha, \beta, C_k > 0$ derart, dass \begin{equation} \E ||W(h)-W(0)||^\alpha \leq C_k ||h||^{d+\beta} \text{ f"ur alle } h \in [-k,k]^d, \ k \in \N,\label{eq:kolmogorov}\end{equation} so gibt es nach dem Satz von Kolmogorov-Chentsov (Verallgemeinerung auf $\R^d$, vgl. Abschnitt \MakeUppercase{\romannumeral 1}.\MakeUppercase{\romannumeral 1}.1.3 in \cite{borodin-2002}) eine pfadstetige Modifikation von $W(\cdot)$. Ob diese Voraussetzung erf"ullt ist, h"angt nur vom Variogramm $\gamma$ ab, denn es ist $W(h) - W(0) \sim \mathcal{N} (0 , 2\gamma(h))$. Wir k"onnen daher f"ur die entsprechenden Variogrammtypen die Pfadstetigkeit ohne Einschr"ankung voraussetzen (vgl. hierzu Abschnitt 2.3.1 in \cite{chiles-1999}).
\item[3. ] Definition und Satz \ref{brown-resnick} lassen sich auch auf nicht-pfadstetige Zufallsfelder erweitern. Aus technischen Gr"unden wird in dieser Arbeit jedoch darauf verzichtet.
\item[4. ] Wir k"onnen auch Brown-Resnick-Prozesse mit Standard"=Fr\'echet"=R"andern (d.\,h. $\Prob(Z(t) \leq x) = \exp(- \frac 1 x) \ \forall x \in \R, \ t \in \R^d$) erhalten durch $$Z(t) = \max_{i \in \N} \ \exp\left(X_i + W_i(t) - \frac{\sigma^2(t)}{2}\right), \ t \in \R^d.$$ Ebenso erreichen wir Standard"=Fr\'echet"=R"ander in den weiteren S"atzen, indem wir auf die dortigen Prozesse vor der Maximumsbildung die Exponentialfunktion anwenden.
\end{itemize}
\end{bem}
\bigskip

\section{Stochastisch "aquivalente Darstellungen mit zuf"alligen Translationen}

\begin{satz} \label{z-trans}
Seien $W$ wie in Satz \ref{brown-resnick}, $W_i \sim_{u.i.v.} W, \ i \in \N$. Sei $Q$ ein Wahrscheinlichkeitsma"s auf \ $\R^d$ und  $\Pi = \sum \delta_{(X_i, H_i)}$ ein von $\{W_i(\cdot), \ i \in \N \}$ unabh"angiger Poisson-Punktprozess auf \ $\R \times \R^d$ mit Intensit"atsma"s $e^{-x}dx \
Q(dh)$. Dann ist
\[ \tilde{Z}(t) :=\max_{i \in \N} \left(X_i + W_i(t-H_i) - \frac{\sigma^2(t-H_i)}{2}\right), \quad t \in \R^d \]
ein Brown-Resnick-Prozess zum Variogramm $\gamma$, d.\,h. $\tilde{Z} \stackrel{d}{=} Z$.
\end{satz}
\begin{bew}
Seien $t_1,\,\ldots,\,t_m \in \R^d, \ y_1,\, \ldots y_m \in \R, \ m \in \N$ beliebig und $\Prob_{\widetilde{t_1},\ldots,\widetilde{t_k}}$ sei das zum Zufallsvektor $(W(\widetilde{t_1}),\ldots,W(\widetilde{t_k}))$ geh"orende Wahrscheinlichkeitsma"s.
Nun ist $\Phi := \sum_{i \in \N} \delta_{(X_i, H_i, W_i)}$ nach Lemma \ref{2dim-ppp} ein Poisson"=Punktprozess auf $(\R \times \R^d \times C(\R^d),\ \mathcal{B} \times \mathcal{B}^m \times \mathcal{C})$ mit Intensit"atsma"s $e^{-x} dx \ Q(dh) \Prob(dW)$ und es gelten die "Aquivalenzen
\begin{align*}
& \tilde{Z}(t_1) \leq y_1, \ldots, \tilde{Z}(t_m) \leq y_m \\
\iff {} & X_i \leq \min_{k=1,\ldots, m} \left(y_k - W_i(t_k-H_i) + \frac{\sigma^2(t_k-H_i)}{2}\right) \quad \forall i \in \N\\
\iff {} & \Phi \Big( \Big\{ (x,h, w) \in \R \times \R^d \times C(\R^d): \\
& \hspace{3cm} x > \min_{k=1,\ldots, m} \left(y_k - w(t_k-h) + \frac{\sigma^2(t_k-h)} 2 \right) \Big\} \Big) = 0.
\end{align*}
Damit ergibt sich
\begin{align*}
&-\log(\Prob(\tilde{Z}(t_1) \leq y_1, \ldots, \tilde{Z}(t_m) \leq y_m)) \\
={}&\int_{\R^d} \int_{\R^m} \int_{\min_{k=1,\ldots, m} \left(y_k - w_k + \frac{\sigma^2(t_k-h)} 2\right)}^\infty e^{-x}dx  \\
& \hspace{5cm} \Prob_{t_1-h,\ldots,t_m-h}(d(w_1,\ldots,w_m)) \ Q(dh)\\
={}&-\log\left(\int_{\R^d} \Prob(Z(t_1-h) \leq y_1, \ldots Z(t_m-h) \leq y_m) \ Q(dh)\right) \\
={}&-\log\left(\int_{\R^d} \Prob(Z(t_1) \leq y_1, \ldots Z(t_m) \leq y_m) \ Q(dh)\right), \\
& \hspace{6.5cm} \text{ da $Z(\cdot)$ station"ar ist nach Satz \ref{brown-resnick}}\\
={}&-\log(\Prob(Z(t_1) \leq y_1, \ldots Z(t_m) \leq y_m)), \\
& \hspace{6.2cm} \text{ da $Q$ ein Wahrscheinlichkeitsma"s ist.}
\end{align*}
Also gilt $\tilde{Z} \stackrel{d}{=} Z$ -- wie behauptet.
\end{bew}
\medskip

Dieser Satz erm"oglicht uns stochastisch "aquivalente Darstellungen des Brown-Resnick-Prozesses. Zwei Spezialf"alle werden in den folgenden Korollaren betrachtet.
\medskip

\begin{coro} \label{z1}
Seien $W$ wie in Satz \ref{brown-resnick}, $W_i^{(j)} \sim_{u.i.v.} W, \ i \in \N, \ j=1,\ldots,n$.
Weiterhin seien \ $\Pi^{(j)} = \sum \delta_{X_i^{(j)}}, \ j=1,\ldots,n$ stochastisch unabh"angige Poisson"=Punktprozesse auf \ $\R$ (unabh"angig von $\{W_i^{(j)}(\cdot),\ i \in \N,\ j \in \{1,\ldots,n\} \})$ mit Intensit"atsma"s \ $\frac{1}{n}\cdot e^{-x}dx$ sowie $h_1,\, \ldots, \,h_n \in \R^d$ fest.
Dann ist 
\[ Z_1(t) := \max_{j=1, \ldots, n} \max_{i \in \N} \left(X_i^{(j)} + W_i^{(j)}(t-h_j) - \frac{\sigma^2(t-h_j)}{2}\right), \quad t \in \R^d \]
ein Brown-Resnick-Prozess zum Variogramm $\gamma$, d.\,h. $Z_1 \stackrel{d}{=} Z$.
\end{coro}
\begin{bew}
Nach Konstruktion ist $\sum_{i \in \N} \delta_{(X_i^{(j)},h_j)}$ f"ur jedes $j=1,\ldots,n$ ein Poisson"=Punktprozess auf $\R \times \R^d$ mit Intensit"atsma"s $\frac{1}{n} \ e^{-x}dx \ \delta_{h_j}(dh)$. Diese Prozesse sind f"ur $j \in \{1,\ldots,n\}$ stochastisch unabh"angig.
\medskip

Daher ist auch $\sum_{j=1}^n \sum_{i \in \N} \delta_{(X_i^{(j)},h_j)}$ ein Poisson"=Punktprozess und das zugeh"orige Intensit"atsma"s ist $\frac{1}{n} e^{-x}dx (\sum_{j=1}^n \delta_{h_j})(dh)$. Mit $Q=\frac 1 n \sum_{j=1}^n \delta_{h_j}$ folgt daher die Aussage des Korollars aus Satz \ref{z-trans}, denn $Q$ ist ein Wahrscheinlichkeitsma"s auf $\R^d$. 
\end{bew}
\medskip

\begin{coro} \label{z2}
Seien $W,\, W_i,\, i \in \N$ wie in Satz \ref{brown-resnick}, $I \subset \R^d$ ein endliches Intervall. Weiterhin sei $\Pi = \sum \delta_{(X_i, H_i)}$ ein davon unabh"angiger Poisson"=Punktprozess auf \ $\R \times I$ mit Intensit"atsma"s \ $|I|^{-1} e^{-x}dx \ dh$.
Dann ist
\[ Z_2(t) :=\max_{i \in \N} \left(X_i + W_i(t-H_i) - \frac{\sigma^2(t-H_i)}{2}\right), \quad t \in \R^d \]
ein Brown-Resnick-Prozess zum Variogramm $\gamma$, d.\,h. $Z_2 \stackrel{d}{=} Z$.
\end{coro}
\begin{bew}
Es ist durch $Q(dh):=\frac 1 {|I|} dh$ das zur Gleichverteilung auf $I$ geh"orende Wahrscheinlichkeitsma"s gegeben. Damit folgt die Aussage aus Satz \ref{z-trans}.
\end{bew}
\bigskip

\section{Brown-Resnick-Stationarit"at}

Der folgende Abschnitt "uber den Begriff der \glqq{}Brown-Resnick-Stationarit"at\grqq{} orientiert sich an Kapitel 2 aus \cite{kabluchko-2008}.
\medskip

Sei $\{\xi(t), \ t \in \R^d \}$ ein stochastischer Prozess mit der Eigenschaft 
\begin{equation} \E(e^{\xi(t)}) < \infty \quad \forall t \in \R^d \label{eq:laplace} \end{equation}
und seien $\xi_i(\cdot) \sim_{u.i.v.} \xi(\cdot), \ i \in \N$.
Sei $\sum_{i \in \N} \delta_{X_i}$ ein davon unabh"angiger Poisson"=Punktprozess auf $\R$ mit Intensit"atsma"s $e^{-x}dx$.
\medskip

\begin{definition} \label{br_stat}
$\xi(\cdot)$ hei"st \emph{Brown-Resnick-station"ar} genau dann, wenn der Prozess $\{\eta(t), \ t \in \R^d\}$, gegeben durch
\[\eta(t) := \max_{i \in \N} (X_i + \xi_i(t)),\]
station"ar ist.
\end{definition}
\medskip

\begin{bem} \label{z_br_stat}
Sei $W(\cdot)$ ein pfadstetiges, intrinsisch station"ares Gau"ssches Zufallsfeld auf \ $\R^d$ mit Variogramm $\gamma$, Erwartungswert 0 und Varianz $\sigma^2(t)$.\\ W"ahlen wir nun $\xi(t) = W(t) - \frac{\sigma^2(t)}{2}$ , dann ist 
\[ \E e^{\xi(t)} = \int_\R e^{w-\frac{\sigma^2(t)}{2}} \frac{1}{\sqrt{2\pi}\sigma(t)} e^{-\frac{1}{2} \left(\frac{w}{\sigma(t)}\right)^2} dw 
 = \int_\R \frac{1}{\sqrt{2\pi}\sigma(t)} e^{-\frac{1}{2} \left(\frac{w-\sigma^2(t)}{\sigma(t)}\right)^2} dw= 1 \]
f"ur alle $t \in \R^d$.
Damit ist $\xi(\cdot)$ Brown-Resnick-station"ar, denn es ist $\eta(\cdot) = Z(\cdot)$ station"ar nach Satz \ref{brown-resnick}.
\end{bem}

Wir wollen nun eine "aquivalente Definition der Brown-Resnick-Stationarit"at geben.
Dazu betrachten wir den Raum $C(\R^d)$ versehen mit der Borel-$\sigma$-Algebra ${\cal C}$.
Nun k"onnen wir $\sum \delta_{X_i + \xi_i(\cdot)}$ nach Satz \ref{ppp-trafo} als Poisson"=Punktprozess auf $C(\R^d)$ auffassen und das zugeh"orige Intensit"atsma"s $\Lambda$ ist das Bildma"s von $e^{-x} dx \times \Prob$ (wobei $\Prob$ das zur Verteilung von $\xi(\cdot)$ geh"orende Wahrscheinlichkeitsma"s sei) unter der messbaren Abbildung $$\pi: \R \times C(\R^d) \rightarrow C(\R^d), \ (x,\xi(\cdot)) \mapsto x + \xi(\cdot),$$ also
\begin{align*}
\Lambda(C_{t_1,\ldots, t_m} (B)) ={}&\int_\R e^{-x} \int_{B-(x,\ldots,x)} \Prob_{t_1,\ldots,t_m}(dw_1,\ldots dw_m) \ dx \\
={}& \int_{\R^d} \int_{\R} \mathbf{1}_{B-(x,\ldots,x)}(w_1,\ldots, w_m) \ e^{-x} dx  \Prob_{t_1,\ldots,t_n}(dw_1,\ldots,dw_m).
\end{align*}

F"ur $t \in \R^d, z \in \Z$ betrachten wir die Menge $A_{t,z}:=\{f \in C(\R^d): f(t) > z\}$. Dann ist
\begin{align}
\Lambda(A_{t,z}) = {} & \int_\R e^{-x} \Prob(\xi(t) > z - x) dx \stackrel{y:=z-x}{=} \int_\R e^{y-z} \Prob(\xi(t) > y) dy  \nonumber\\
\stackrel{u:= e^y}{=}{}& e^{-z} \int_{\R_{> 0}} \Prob(e^{\xi(t)} > u) du = e^{-z} \E e^{\xi(t)} < \infty \label{eq:atzendl}
\end{align}
nach \eqref{eq:laplace}. Weiterhin ist jede beschr"ankte Menge der Form $C_{t_1,\ldots,t_m}(B) \in {\cal C}$ in einem solchen $A_{t,z}$ enthalten. Daher ist $\Lambda$ lokal endlich.
\medskip

\begin{satz} \label{ppp_br_stat}
$\{ \xi(t), \ t \in \R^d \}$ ist Brown-Resnick-station"ar genau dann, wenn der Poisson"=Punktprozess $\sum_{i \in \N} \delta_{X_i + \xi_i(\cdot)}$ auf \ $C(\R^d)$ translationsinvariant ist (d.\,h., dass das zugeh"orige Intensit"atsma"s translationsinvariant ist).
\end{satz}
\begin{bew}
Mit den Notationen von oben gilt
\begin{align}
& \Prob(\eta(t_1) \leq y_1, \ldots, \eta(t_m) \leq y_m) \nonumber \\
={}&\Prob \left( \sum_{i \in \N} \delta_{X_i + \xi_i(\cdot)} \left(\bigcup_{j=1,\ldots, m} \{ f \in C(\R^d): f(t_j) > y_j \}\right)=0\right) \nonumber \\
={}&\exp(-\Lambda(C_{t_1,\ldots,t_m}(B))) \label{eq:etaC}
\end{align}
mit $B= \R^m \backslash \bigtimes_{j=1}^m (-\infty, y_j]$.
\medskip

Ist nun der Punktprozess $\sum \delta_{X_i + \xi_i(\cdot)}$ translationsinvariant, so ist $\Lambda$ translationsinvariant, und f"ur jedes $h \in \R^d$ folgt nach \eqref{eq:etaC}
\begin{align*}
 \Prob(\eta(t_1) \leq y_1, \ldots, \eta(t_m) \leq y_m)
={}&\exp(-\Lambda(C_{t_1,\ldots,t_m}(B))) \\
 = \quad \exp(-\Lambda(C_{t_1+h,\ldots,t_m+h}(B)))
={}&\Prob(\eta(t_1+h) \leq y_1, \ldots, \eta(t_m+h) \leq y_m),
\end{align*}
also ist $\eta(\cdot)$ station"ar.
\medskip

Ist andererseits $\eta(\cdot)$ station"ar, so ergibt sich f"ur alle $h \in \R^d$ wiederum nach \eqref{eq:etaC} die Gleichung
\begin{equation} \exp(-\Lambda(C_{t_1,\ldots,t_m}(B))) = \exp(-\Lambda(C_{t_1+h,\ldots,t_m+h}(B))) \label{eq:lambdastat} \end{equation}
f"ur alle Mengen $B= \R^m \setminus \bigtimes_{j=1}^m (-\infty, y_j]$.

F"ur feste $t_1, \ldots, t_m \in \R^d$ definieren wir nun $\mathcal{C}_{t_1,\ldots,t_m} := \{ C_{t_1,\ldots,t_m}(A), \ A \in \mathcal{A}) \}$ mit $\mathcal{A} := \{ \R^m \setminus \bigtimes_{j=1}^m (-\infty, y_j]: \ y_1, \ldots, y_m \in \R \}$. Da nun $\mathcal{A}$ stabil unter Vereinigungen ist und au"serdem $C_{t_1,\ldots,t_m}(B_1) \cup C_{t_1,\ldots,t_m}(B_2) = C_{t_1,\ldots,t_m}(B_1 \cup B_2)$ f"ur alle $B_1, B_2 \in \mathcal{B}^m$ gilt, ergibt sich die Stabilit"at von $\mathcal{C}_{t_1,\ldots,t_m}$ unter Vereinigung. Mit der Siebformel sowie der Gleichung $C_{t_1,\ldots,t_m}(B_1) \cap C_{t_1,\ldots,t_m}(B_2) = C_{t_1,\ldots,t_m}(B_1 \cap B_2)$ folgt nun aus \eqref{eq:lambdastat}, dass auch der Schnitt zweier Mengen aus $\mathcal{C}_{t_1,\ldots,t_m}$
die Beziehung
\begin{equation}
\Lambda(C_{t_1,\ldots,t_m}(B)) = \Lambda(C_{t_1+h,\ldots,t_m+h}(B)) \ \forall h \in \R^d \label{eq:lambdastat2}
\end{equation} 
erf"ullt.

Induktiv so fortfahrend und dabei die Siebformel benutzend erh"alt man, dass die Gleichung \eqref{eq:lambdastat2} auch f"ur jedes Element $C_{t_1,\ldots,t_m}(B)$ der (nach Konstruktion schnittstabilen) Menge
$$ \mathcal{M}_{t_1, \ldots, t_m} :={} \{ C_1 \cap \ldots \cap C_n, \ n \in \N, C_i \in \mathcal{C}_{t_1,\ldots,t_m} \}$$
gilt.

Damit gilt die Translationsinvarianz auch f"ur die Elemente von $$\mathcal{M} := \bigcup_{\substack{t_1, \ldots, t_m \in \R^d,\\ m \in \N}} \mathcal{M}_{t_1, \ldots, t_m}$$ und $\mathcal{M}$ ist schnittstabil.

Weiterhin ist $\mathcal{A}$ ein Erzeuger von $\mathcal{B}^m$ und somit erzeugt $$\bigcup_{\substack{t_1, \ldots, t_m \in \R^d,\\ m \in \N}} \mathcal{C}_{t_1, \ldots, t_m} \subset \mathcal{M}$$ die $\sigma$-Algebra $\mathcal{C}$. Daher ist $\mathcal{M}$ ein schnittstabiler Erzeuger von $\mathcal{C}$.

F"ur jedes $t \in \R^d$ und $B_n := \R \setminus (-\infty, -n] \in \mathcal{A}$ ist die Folge $\left(C_t(B_n)\right)_{n \in \N} \subset \mathcal{M}$ monoton wachsend gegen $C(\R^d)$ und wegen \eqref{eq:atzendl} gilt $$\Lambda(C_t(B_n)) = \Lambda(A_{t,-n}) < \infty$$ f"ur alle $n \in \N$. Damit folgt nach dem Eindeutigkeitssatz "uber Ma"se die Translationsinvarianz von $\Lambda$ auf ganz $\mathcal{C}$.
\end{bew}
\medskip

\begin{satz} \label{z3}
Seien $W(\cdot)$ ein pfadstetiges, intrinsisch station"ares Gau"ssches Zufallsfeld auf \ $\R^d$ mit Variogramm $\gamma$, Erwartungswert 0 und Varianz $\sigma^2(\cdot)$,\\ $W_i^{(j)} \sim_{u.i.v.} W,\ i\in \N, \ j \in \Z^d$ und $T_i^{(j)} := \inf \left( \emph{\argsup}_{t \in \R^d} \left(W_i^{(j)}(t) - \frac{\sigma^2(t)}{2}\right)\right)$. Es gelte
\begin{equation} \lim_{|t| \to \infty} \left(W(t) - \frac{\sigma^2(t)}{2}\right) = -\infty \quad \Prob-f.s. \label{eq:lim} \end{equation}
Weiterhin seien $\Pi^{(j)} = \sum \delta_{X_i^{(j)}}$ stochastisch unabh"angige Poisson"=Punktprozesse auf \ $\R$ (unabh"angig von $\{W_i^{(j)}(\cdot),\ i \in \N,\ j \in \Z^d \})$  mit Intensit"at $\frac 1{m^d} e^{-x}dx$, $m \in \N$ sowie $p \in \R_{> 0}^d$.
Dann ist
\[ Z_3(t) := \max_{j \in \Z^d} \max_{\substack{i \in \N \\ T_i^{(j)} \in \left(-\frac{m p}{2}, \frac{m p}{2} \right] }} \left(X_i^{(j)} + W_i^{(j)}(t - p \cdot j) - \frac{\sigma^2(t - p \cdot j)}{2}\right), \quad t \in \R^d \]
ein Brown-Resnick-Prozess zum Variogramm $\gamma$, d.\,h. $Z_3 \stackrel{d}{=} Z$.

Dabei ist die Multiplikation $p \cdot j$ als komponentenweise Multiplikation zu verstehen, d.\,h. $(p_1,\ldots, p_d) \cdot (j_1,\ldots,j_d) := (p_1 j_1, \ldots, p_d j_d)$.
\end{satz}
\begin{bew}
Sei $\xi_i^{(j)}(t) := W_i^{(j)}(t) - \frac{\sigma^2(t)}{2}$.\\
Zun"achst bemerken wir, dass jedes $T_i^{(j)}$ aufgrund von Bedingung \eqref{eq:lim} $\Prob$-f.s. endlich ist und auch
$M_i^{(j)} = \sup_{t \in \R^d} (X_i^{(j)} + \xi_i^{(j)}(t))$ wohldefiniert ist.
\medskip

Wir zeigen nun, dass die Transformation 
\[ \Phi: C(\R^d) \to \R^d \times C(\R^d), \ \ X_i^{(j)} + \xi_i^{(j)}(\cdot) \mapsto (T_i^{(j)}, X_i^{(j)} + \xi_i^{(j)}(\cdot))\]
messbar ist.\\
Sei dazu $t_1, t_2, \ldots$ eine Abz"ahlung von $\Q^d$. Wir definieren nun  f"ur $k \in \N$ die messbare Abbildung $M_k: C(\R^d) \to \R, \ \xi \mapsto \xi(t_k)$. Dann ist auch die Abbildung $M: C(\R^d) \to \R, \ \xi \mapsto \sup_{k \in \N} M_k(\xi)$ als punktweises Supremum abz"ahlbar vieler messbarer Funktionen wieder messbar und es ist $M(\xi_i^{(j)}) = M_i^{(j)}$. Damit ist weiterhin die Abbildung $\widetilde{\Phi}: C(\R^d) \to C(\R^d), \ \ X_i^{(j)} + \xi_i^{(j)}(\cdot) \mapsto X_i^{(j)} + \xi_i^{(j)}(\cdot) - M_i^{(j)}$ messbar. Nun ist nach Konstruktion $T_i^{(j)}$ die \glqq kleinste\grqq{} Nullstelle von $\widetilde{\Phi}(\xi_i^{(j)})$. Wir definieren f"ur $k, n \in \N$ die messbaren Abbildungen
\[ \tau_k^{(n)}: C(\R^d) \to \R^d, \ \tau_k^n(\xi) = \begin{cases} t_k, & \text{falls } \xi(t_k) > -\frac 1 n\\ (\infty, \ldots, \infty), & \text{sonst} \end{cases}\]
und schlie"slich $\displaystyle \tau: C(\R^d) \to \R^d, \ \xi \mapsto \lim_{n \to \infty} \inf_{k \in \N} \tau_k^{(n)}(\xi)$.\\
$\tau$ ordnet nun jeder (im Unendlichen negativen) Funktion ihre \glqq kleinste\grqq{} Nullstelle zu und ist weiterhin messbar, da punktweise Limites und Infima abz"ahlbar vieler messbarer Funktionen wieder messbar sind. \\
Da nun $\Phi(\xi) = (\tau\circ\widetilde{\Phi}(\xi),\ \xi)$ gilt, ist $\Phi$ messbar.
\medskip

Dadurch erhalten wir einen Poisson"=Punktprozess \[\sum_{i \in \N} \delta_{(T_i^{(j)}, X_i^{(j)} + \xi_i^{(j)}(\cdot))}\] mit transformiertem Intensit"atsma"s 
\[ \Psi(A) = \int_\R \frac 1{m^d} e^{-x} \Prob( x + \xi \in \Phi^{-1}(A))dx, \quad A \in \mathcal{B}^d \times \mathcal{C}, \]
wobei $\Prob$ das zu $W(\cdot)$ geh"orende Wahrscheinlichkeitsma"s sei.
Seien nun
\begin{align*}
&U_t: C(\R^d) \to C(\R^d), \ f(\cdot) \mapsto f(\cdot - t) \\
\text{und } &V_t: \R^d \times C(\R^d) \to \R^d \times C(\R^d), \ (x, f(\cdot)) \mapsto (x + t,f(\cdot - t))
\end{align*}
Verschiebungen um ein $t \in \R^d$.\\
Dann gilt
\begin{align*}
(\Phi \circ U_t) (X_i^{(j)} + \xi_i^{(j)}(\cdot)) ={}& \Phi(X_i^{(j)} + \xi_i^{(j)}(\cdot - t))\\
 ={}& (T_i^{(j)} + t, X_i^{(j)}+\xi_i^{(j)}(\cdot - t))\\
 ={}& V_t(T_i^{(j)}, X_i^{(j)}+\xi_i^{(j)}(\cdot))\\
 ={}& (V_t \circ \Phi) (X_i^{(j)} + \xi_i^{(j)}(\cdot)).
\end{align*}
Da $\xi(\cdot)$ nun nach Bemerkung \ref{z_br_stat} Brown-Resnick-station"ar ist, ist das Intensit"atsma"s von $\sum \delta_{X_i + \xi_i(\cdot)}$ nach Satz \ref{br_stat} translationsinvariant (bzgl. $U_t$) und damit auch das Intensit"atsma"s von $\sum \delta_{X_i^{(j)} + \xi_i^{(j)}(\cdot)}$, denn beide Ma"se unterscheiden sich nur um den Faktor $m^d$. Aufgrund der Tatsache, dass $\Phi$ mit den Verschiebungsoperatoren kommutiert, ist dann auch das transformierte Intensit"atsma"s $\Psi$ translationsinvariant (bzgl. $V_t$).
\medskip

Daraus folgt f"ur jedes $j \in \Z^d$:
\begin{align}
& \max_{\substack{i \in \N \\ T_i^{(j)} \in \left(-\frac{m p}{2}, \frac{m p}{2} \right]}} \left(X_i^{(j)} + \xi_i^{(j)}(\cdot - p \cdot j)\right) \nonumber \\
={} & \max_{\substack{i \in \N \\ T_i^{(j)} + p\cdot j \in \left(p\cdot j-\frac{mp}{2}, p \cdot j+\frac{mp}{2} \right]}} \left(X_i^{(j)} + \xi_i^{(j)}(\cdot - p\cdot j)\right) \nonumber \\
 \stackrel{d}{=}{} & \max_{\substack{ i \in \N \\ T_i^{(j)} \in \left(p\cdot j-\frac{mp}{2}, p \cdot j+\frac{mp}{2} \right] }} \left(X_i^{(j)} + \xi_i^{(j)}(\cdot)\right) \label{eq:transinv}
\end{align}

Betrachtet man nun jede Seite der Gleichung \eqref{eq:transinv} f"ur verschiedene $j \in \Z^d$, so erh"alt man jeweils stochastisch unabh"angige Prozesse. Daraus folgt:
\begin{align*}
Z_3(\cdot) ={}&\max_{j \in \Z^d} \max_{\substack{ i \in \N \\ T_i^{(j)} \in \left(-\frac{m p}{2}, \frac{m p}{2} \right]}} \left(X_i^{(j)} + \xi_i^{(j)}(\cdot - p\cdot j)\right) \\
&\stackrel{d}{=} \max_{j \in \Z^d} \max_{\substack{ i \in \N \\ T_i^{(j)} \in \left(p\cdot j-\frac{mp}{2}, p \cdot j+\frac{mp}{2} \right] }} \left(X_i^{(j)} + \xi_i^{(j)}(\cdot)\right)
\end{align*}
Weiterhin ist
\begin{align}
& \max_{\substack{ i \in \N \\ T_i^{(j)} \in \left(p\cdot j-\frac{mp}{2}, p \cdot j+\frac{mp}{2} \right]}} \left(X_i^{(j)} + \xi_i^{(j)}(\cdot)\right) \nonumber \\
 \stackrel{d}{=} &
\max_{\substack{ i \in \N \\ T_i^{(j \bmod m)} \in \left(p\cdot j-\frac{mp}{2}, p \cdot j+\frac{mp}{2} \right]}} \left(X_i^{(j \bmod m)} + \xi_i^{(j \bmod m)}(\cdot)\right) \label{eq:jstat},
\end{align}
wobei die Modulo-Operation hier komponentenweise zu verstehen ist.
Da die entsprechenden Mengen $\{ i: T_i^{(j \bmod m)} \in \left(p\cdot j-\frac{mp}{2}, p \cdot j+\frac{mp}{2} \right] \}$ disjunkt sind, sind die Prozesse auf jeder der beiden Seiten von \eqref{eq:jstat} f"ur $j_l \in \Z^d, \ l \in \N,$ mit $j_l \equiv j \bmod m$ f"ur alle $l \in \N$ jeweils stochastisch unabh"angig.
Somit gilt auch
\begin{align*}
Z_3(\cdot) &\stackrel{d}{=} \max_{j \in \Z^d} \max_{\substack{i \in \N \\ T_i^{(j \bmod m)} \in \left(p\cdot j-\frac{mp}{2}, p \cdot j+\frac{mp}{2} \right]}} \left(X_i^{(j \bmod m)} + \xi_i^{(j \bmod m)}(\cdot)\right)\\
={}&\max_{k \in \{0, \ldots, m-1\}^d}  \max_{\substack{ j \in \Z^d\\ j \bmod m \equiv k }} \left( \max_{\substack{ i \in \N \\ T_i^{(k)} \in \left(p\cdot j-\frac{mp}{2}, p \cdot j+\frac{mp}{2} \right]}} \left(X_i^{(k)} + \xi_i^{(k)}(\cdot)\right) \right)
\end{align*}
\begin{align*}
={}&\max_{k \in \{0,\ldots\, m-1 \}^d} \max_{\substack{ i \in \N \\ T_i^{(k)} \in \R^d}} \left(X_i^{(k)} + \xi_i^{(k)}(\cdot)\right),\\
  & \qquad \qquad
     \text{da } \bigcup_{\substack{j \in \Z^d\\ j \bmod m \equiv k}} \left(p\cdot j-\frac{mp}{2}, p \cdot j+\frac{mp}{2} \right] = \R^d \\
 ={}& \max_{k \in \{0,\ldots\, m-1 \}^d} \max_{i \in \N} \left(X_i^{(k)} + \xi_i^{(k)}(\cdot)\right) \stackrel{d}{=} Z(\cdot).
\end{align*}
Dabei benutzen wir im letzten Schritt, dass $\sum_{k \in \{0, \ldots, m-1\}^d} \sum_{i \in \N} \delta_{X_i^{(k)}}$ als Summe von stochastisch unabh"angigen Poisson"=Punktprozessen auch ein Poisson"=Punktprozess ist mit Intensit"atsma"s $\sum_{k \in \{0, \ldots, m-1\}^d} \frac 1 {m^d} e^{-x} dx$ $ = e^{-x} dx$.
\end{bew}
\medskip

Um andere Darstellungen des Brown-Resnick-Prozesses zu erhalten, verwenden wir das nachfolgende Lemma, das eine Folgerung aus Satz 13 aus \cite{kabluchko-2008} ist. Der Beweis verl"auft daher analog.

\begin{lemma} \label{mmm}
Seien $W(\cdot)$ ein pfadstetiges, intrinsisch station"ares Gau"ssches Zufallsfeld auf \ $\R^d$ mit Variogramm $\gamma$, Erwartungswert 0 und Varianz $\sigma^2(\cdot)$,\\
$W_i(\cdot) \sim_{u.i.v.} W(\cdot), \ i\in \N,$ sowie $T_i := \inf \left( \emph{\argsup}_{t \in \R^d} \left(W_i(t) - \frac{\sigma^2(t)}{2}\right)\right)$,\\
$M_i := \sup_{t \in \R^d} \left(W_i(t) - \frac{\sigma^2(t)}{2}\right)$ und schlie"slich $F_i(\cdot) :=  W_i(\cdot+T_i) - \frac {\sigma^2(\cdot+T_i)} 2 - M_i$.
Es gelte \eqref{eq:lim}.
\medskip

Weiterhin sei durch $\sum_{i \in \N} \delta_{X_i}$ ein von $\{W_i(\cdot), \ i \in \N\}$ unabh"angiger Poisson"=Punktprozess mit Intensit"atsma"s $e^{-x}dx$ gegeben.

Dann ist $\sum_{i \in \N} \delta_{(T_i, X_i + M_i, F_i)}$ ein Poisson"=Punktprozess auf \ $\R^d \times \R \times C(\R^d)$,
dessen Intensit"atsma"s in ein Produkt zerf"allt.
\end{lemma}

\begin{bew}
Wir definieren zun"achst 
\begin{align*}
&\xi(\cdot) := W(\cdot) - \frac {\sigma^2(\cdot)} 2,\\
&M := \sup_{t \in \R^d} \xi(t),\\
&T := \inf \left( \argsup \left(\xi(t)\right)\right)\\
\text{und } & F(\cdot) := \xi(\cdot+T)-M.
\end{align*}

Um zu zeigen, dass $\sum_{i \in \N} \delta_{(T_i, X_i + M_i, F_i)}$  ein Poisson"=Punktprozess ist, beweisen wir, dass die Abbildung $\Phi: \R \times C(\R^d) \to \R^d \times \R \times C(\R^d), \ (u,\xi) \mapsto (T,u+M,F)$ messbar ist.\\
Wir wissen bereits, dass die Abbildung $\xi \mapsto T$ messbar ist (vgl. Beweis Satz \ref{z3}).
Weiterhin betrachten wir eine Abz"ahlung $t_1, t_2, \ldots$ von $\Q^d$. Nun sind die Abbildungen $f_i: C(\R^d) \to \R, \ \xi(\cdot) \ \mapsto \xi(t_i)$, $i \in \N$, messbar und somit auch die Abbildung auf das punktweise Supremum $\xi \mapsto M = \sup_{i \in \N} f_i(\xi)$. 
Ebenso ist auch die Abbildung $\widetilde{\Phi}: \R^d \times C(\R^d) \to C(\R^d), \ (t, \xi(\cdot)) \mapsto \xi(\cdot + t)$ messbar, denn f"ur alle $t_1, \ldots, t_m, \ m \in \N$ und jede offene Teilmenge $B \subset \R^m$ gilt aufgrund der Stetigkeit, dass
$$ \widetilde{\Phi}^{-1} (C_{t_1, \ldots, t_m}(B)) = \R^d \times \bigcap_{t \in \Q^d} C_{t_1-t, \ldots, t_m-t}(B) \in \mathcal{B}^d \times \mathcal{C}.$$
Somit ist die Abbildung $\xi(\cdot) \mapsto (T, \xi(\cdot)-M) \stackrel{\widetilde{\Phi}}{\mapsto} F$ messbar und daher auch $\Phi$.

Damit ist $\sum_{i \in \N} \delta_{(T_i, X_i + M_i, F_i)}$ ein Poisson"=Punktprozess auf $\R^d \times \R \times C(\R^d)$ mit Intensit"atsma"s
   $$\Psi(A) = \int_{\R} e^{-x} \Prob(A-(0,x,0)) dx, \quad \ A \in \mathcal{B}^d \times \mathcal{B} \times \mathcal{C}$$
wobei $\Prob$ das zu ($T$, $M$, $F$) geh"orende Wahrscheinlichkeitsma"s sei.

Nun gilt f"ur alle $y,z \in \R$, $t \in \R^d$, $ f \in C(\R^d)$ und f"ur jede Borelmenge $A \subset \R^d \times \R \times C(\R^d)$ die Gleichung
\begin{align*}
\int_{A+(0,z,0)} e^{y} d\Psi(t,y,f) ={} & \int_{\R} e^{-x} \int_{A+(0,z,0)} e^{y} d\Prob(t,y-x,f) \ dx \\
 ={} & \int_{\R} e^{-x} \int_{A} e^{y+z}  d\Prob(t,y+z-x,f) \ dx\\
 ={}& \int_{\R} e^{-x-z} \int_{A} e^{y+z} d\Prob(t,y-x,f) \ dx\\
 ={} & \int_{\R} e^{-x} \int_{A} e^{y} d\Prob(t,y-x,f) \ dx = \int_{A} e^{y} d\Psi(t,y,f),
\end{align*}
also ist $e^y d\Psi(t,y,f)$ ein translationsinvariantes Ma"s in der zweiten Komponente.

Seien weiterhin
\begin{align*}
&U_t: C(\R^d) \to C(\R^d), \ f(\cdot) \mapsto f(\cdot - t)\\
\text{und} \quad &V_t: \R^d \times \R \times C(\R^d) \to \R^d \times \R \times C(\R^d), \ (x, y, f(\cdot)) \mapsto (x + t, y, f(\cdot))
\end{align*}
Verschiebungen um ein $t \in \R^d$.\\
Dann gilt:
\[\begin{array}{ccl}
(\Phi \circ U_t) (X_i + \xi_i(\cdot)) &={}& \Phi(X_i + \xi_i(\cdot - t))\\
 &={}& (T_i + t,\ X_i + M_i,\ \xi_i(\cdot - t + T_i + t) - M_i)\\
 &={}& V_t(T_i,\ X_i +M_i,\ \xi_i(\cdot - t + T_i + t) - M_i)\\
 &={}& (V_t \circ \Phi) (X_i + \xi_i(\cdot)).
\end{array}\]
Da der Poisson"=Punktprozess $\sum_{i \in \N} \delta_{X_i + W_i(\cdot) - \frac {\sigma^2(\cdot)} 2}$ nach Bemerkung \ref{z_br_stat} und Satz \ref{ppp_br_stat} ein translationsinvariantes Intensit"atsma"s hat (bzgl. $U_t$) und die Translationsoperatoren mit $\Phi$ kommutieren, ist auch das Bildma"s $\Psi$ bzgl. $V_t$, also in der ersten Komponente, translationsinvariant, d.\,h.
$$\Psi(A+(t,0,0))= \Psi(A) \quad \text{f"ur alle } A \in \mathcal{B}^d \times \mathcal{B} \times \mathcal{C}, \ t \in \R^d.$$

Nun ist
\begin{align}
& \Psi([0,1]^d \times [0,1] \times C(\R^d)) \nonumber\\
={} & \int_\R e^{-x} \Prob\left(\sup_{t \in \R^d} \left(W(t) - \frac {\sigma^2(t)} 2\right) \in [-x, -x +1], \ T \in [0,1]^d\right) dx \nonumber \\
\leq{} & \int_\R e^{-x} \Prob\left(\sup_{t \in [0,1]^d} W(t) \geq -x\right) dx 
={}  \int_\R e^{x} \Prob\left(\sup_{t \in [0,1]^d} W(t) \geq x\right) dx. \label{eq:abschaetzungwuerfel}
\end{align}

Da sich au"serdem $\widetilde{M} := \sup_{t \in [0,1]^d} W(t)$ als Supremum der Gau"sschen Zufallsvariablen $W(t_i), \ i \in \N, \ t_i \in (\Q \cap [0,1])^d,$ schreiben l"asst und aufgrund der Pfad"-stetigkeit von $W(\cdot)$ die Beziehung $\Prob(|\widetilde{M}| < \infty) = 1 > 0$ gilt, folgt nach Satz 5 aus \cite{landau-1970}, dass ein $\varepsilon > 0$ existiert mit
\begin{equation} \E(\exp(\varepsilon |\widetilde{M}|^2)) < \infty. \label{eq:endlerw} \end{equation}

Andererseits ist
\begin{align*}
\E(\exp(\varepsilon |\widetilde{M}|^2)) ={} & \int_0^\infty \Prob\left(\exp(\varepsilon |\widetilde{M}|^2) > z\right) dz \\
 ={} & 1 + \int_1^\infty \Prob\left(\exp(\varepsilon |\widetilde{M}|^2) > z\right) dz \\
={} & 1 + \int_0^\infty 2 \varepsilon z \exp(\varepsilon z^2) \Prob\left(|\widetilde{M}| > z \right) dz.
\end{align*}

Es folgt mit \eqref{eq:abschaetzungwuerfel} die Ungleichungskette
\begin{align*}
\Psi([0,1]^d \times [0,1] \times C(\R^d)) \leq{}& \int_\R e^{x} \Prob\left(\sup_{t \in [0,1]^d} W(t) \geq x\right) dx\\
\leq{} & \int_{-\infty}^{\frac 1 \varepsilon} e^{x} dx + \int_{\frac 1 \varepsilon}^\infty e^{x} \Prob(|\widetilde{M}| \geq x) dx\\
\leq{} & \exp\left(\frac 1 \varepsilon \right) + \int_0^\infty 2 \varepsilon x \exp(\varepsilon x^2) \Prob\left(|\widetilde{M}| > x \right) dx \\
<{} & \exp\left(\frac 1 \varepsilon \right) + \E(\exp(\varepsilon |\widetilde{M}|^2)) < \infty \quad \text{nach \eqref{eq:endlerw}}.
\end{align*}

Somit ist das Ma"s $e^{y} d\Psi(t,y,f)$ in den ersten beiden Komponenten translationsinvariant und jede Menge,
 deren Projektion auf diese Komponenten im Einheitsw"urfel enthalten ist, hat unter diesem Ma"s eine endliche Masse.
Also ist f"ur $A \in \mathcal{C}$, $B \in \mathcal{B}^d \times \mathcal{B}$ durch
$$\Psi_A(B) := \int_{B \times A} e^y d\Psi(t,y,f)$$ ein translationsinvariantes Ma"s auf
$(\R^d \times \R, \mathcal{B}^d \times \mathcal{B})$ definiert, das auf dem Einheitsw"urfel endlich ist.
Daher ist dieses Ma"s ein Vielfaches des Lebesgue-Ma"ses, wobei der (endliche) Vorfaktor von $A \in \mathcal{C}$ abh"angt,
d.\,h. $\Psi_A(B) = \lambda(A) \int_B e^y dy \ dt$ mit $0 < \lambda(A) < \lambda(C(\R^d)) < \infty$.\\
Weiterhin ist durch $Q(\cdot) := \frac {\lambda(\cdot)} {\lambda(C(\R^d))}$ ein Wahrscheinlichkeitsma"s auf
$(C(\R^d), \mathcal{C})$ gegeben.
Mit $0 < \lambda^*:=\lambda(C(\R^d)) < \infty$ gilt also $$ d\Psi(t,y,f) =  \lambda^* e^{-y} dy \ dt \times dQ$$.
\end{bew}
\medskip

\begin{satz} \label{z4}
Seien $\lambda^*$ und $Q$ wie im Beweis von Lemma \ref{mmm}. Weiterhin sei durch 
$\Pi := \sum_{i \in \N} \delta_{(S_i, U_i)}$ ein Poisson"=Punktprozess mit Intensit"atsma"s $\lambda^* e^{-u} du \ ds$ gegeben 
und $\tilde F_i \sim_{u.i.v.} Q$ unabh\"angig von $\Pi$.
\medskip

Dann ist
\begin{align*}
Z_{4,\lambda^*}(t) := & \max_{i \in \N} \Bigg( U_i + \tilde F_i(t-S_i) \Bigg), \ t \in \R^d,
\end{align*}
ein Brown-Resnick-Prozess zum Variogramm $\gamma$ ist, d.\,h. $Z_{4,\lambda^*}(\cdot) \stackrel{d}{=} Z(\cdot)$.
\end{satz}

\begin{bew}
Seien $X_i$, $W_i$, $T_i$, $M_i$, $F_i$ wie in Lemma \ref{mmm}.
Dann sind nach Lemma \ref{mmm} und Lemma \ref{2dim-ppp} sowohl $\sum_{i \in \N} \delta_{(T_i,X_i+M_i,\tilde F_i)}$ als auch
$\Pi$ Poisson-Punktprozesse auf $\R^d \times \R \times C(\R^d)$ mit Intensit\"atsma"s $dt \times e^{-u} du \times Q(dF)$.
Weiterhin lassen sich $Z(\cdot)$ und $Z_{4,\lambda^*}(\cdot)$ auf dieselbe Art und Weise aus diesen Punktprozessen gewinnen, n\"amlich
$Z(\cdot) = \max_{i \in \N} \Gamma(T_i, X_i +M_i, F_i)$ bzw. $Z_{4,\lambda^*} = \max_{i \in \N} \Gamma(S_i, U_i, \tilde F_i)$
mit $\Gamma: \R^d \times \R \times C(\R^d) \to C(\R^d), \ (T,U,F) \mapsto U + F(\cdot-T)$. Damit sind $Z(\cdot)$
 und $Z_{4,\lambda^*}(\cdot)$ in Verteilung gleich.
\end{bew}
\medskip

\begin{bem} \label{rem:lattice}
 Ein \"ahnliches Resultat erhalten wir, wenn wir die Prozesse aus Satz \ref{z4} auf $p\Z^d$, $p>0$, einschr\"anken.

 Dann ist f\"ur
 \begin{align*}
 T_i^{(p)} ={} & \inf \left( \textrm{\argsup}_{t \in p\Z^d} \left(W_i(t) - \frac{\sigma^2(t)}{2}\right)\right),\\
 M_i^{(p)} ={} & \sup_{t \in p\Z^d} \left(W_i(t) - \frac{\sigma^2(t)}{2}\right)\\
 \textrm{and } \quad F_i^{(p)}(\cdot) ={} &  W_i(\cdot+T_i^{(p)}) - \frac {\sigma^2(\cdot+T_i^{(p)})} 2 - M_i^{(p)}, \quad t \in p\Z^d,
 \end{align*}
 das zuf\"allige Ma"s $\sum_{i \in \N} \delta_{(T_i^{(p)}, X_i + M_i^{(p)}, F_i^{(p)})}$ ein Poisson-Punktprozess
 auf $p\Z^d \times \R \times \R^{p\Z^d}$ mit Intensit\"atsma"s $\lambda^{(p)} p^d \delta_t \times e^{-u} d u \times 
 Q^{(p)}(d F)$ f\"ur ein $\lambda^{(p)} > 0$ und ein Wahrscheinlichkeitsma"s $\tilde Q^{(p)}$. Eine \"aquivalente Darstellung
 $Z_{4,\lambda^{(p)}}$ von $Z$ ergibt sich analog zu Satz \ref{z4}.
 Ebenso k\"onnen alle bisherigen Resultate auch auf Prozesse mit eingeschr\"anktem Definitionsbereich \"ubertragen werden.
\end{bem}

Um $Z_4$  zu approximieren, m\"ussen Prozesse mit der Verteilung $Q$ bzw. $Q^{(p)}$ simuliert werden. Dazu ist es notwendig, diese
Ma"se explizit zu kennen. Man beachte dabei, dass $Q$ im Allgemeinen nicht die Verteilung von $W(\cdot + T) - \sigma^2(\cdot+T) - M$ ist
(und $Q^{(p)}$ nicht die Verteilung von $W(\cdot + T^{(p)}) - \sigma^2(\cdot+T^{(p)}) - M^{(p)}$).
Dies erkennt man am folgenden Resultat, nach dem $Q^{(p)}$ die Verteilung von $W(\cdot) - \frac{\sigma^2(\cdot)} 2 \ | \ T^{(p)} = 0$
ist, falls $W(0)=0$. Im Fall einer Brownschen Bewegung $W$ ist aber eine verschobene Brownsche Bewegung anders verteilt als eine Brownsche Bewegung
bedingt darauf, dass sich das Maximum im Ursprung befindet.

\begin{satz} \label{thm:conditional1}
 Sei $W(\cdot)$ wie in Lemma \ref{mmm} und
 $$T^{(p)} = \inf \left( \emph{\argsup}_{t \in p\Z^d} \left(W(t) - \frac{\sigma^2(t)}{2}\right)\right).$$
 Wir nehmen weiterhin $W(0) = 0$ an.
 Dann ist  $Q^{(p)}$ die Verteilung
 $$W(\cdot) - \frac{\sigma^2(\cdot)} 2 \ | \ T^{(p)} = 0.$$
\end{satz}

\begin{bew}
  Seien $A \in \mathcal{B}(\R^{\lambda \Z^d})$ und $V \in \mathcal{B}(\R)$ so, dass
 $$0 < \int_V e^{-x} d x < \infty.$$
 Weiterhin sei $\Pi = \sum_{i \in \N} \delta_{(T_i^{(p)}, X_i + M_i^{(p)}, F_i^{(p)})}$ der Poisson-Punktprozess
 auf $p\Z^d \times \R \times \R^{p\Z^d}$ wie in Bemerkung \ref{rem:lattice}.
 Dann ist
 \begin{equation}
  {Q}^{(p)}(A) = \Prob(\Pi(\{0\} \times V \times A) = 1 \ | \ \Pi(\{0\} \times V \times \R^{p\Z^d}) =1), \label{eq:condprob1}
 \end{equation}
 da das Intensit\"atsma"s von $\Pi$ ein Produktma"s ist.

Wenn wir annehmen, dass die Punkte $(T_i^{(p)}, X_i + M_i^{(p)}, F_i^{(p)})$ so indiziert sind, dass die Folge $(X_i)_{i \in  \N}$
realisierungsweise monoton fallend ist (siehe Kapitel 3), erhalten wir
 \begin{align}
  & \Prob(\Pi(\{0\} \times V \times A) = 1 \ | \ \Pi(\{0\} \times V \times \R^{p\Z^d}) =1) \nonumber \\
 ={} & \sum_{i \in \N} \Prob(T_i^{(p)} = 0, \ X_i+M_i^{(p)} \in V \ | \ \#\{i: \ (T_i^{(p)}, X_i + M_i^{(p)}) \in \{0\} \times V \} =1) \nonumber\\
  & \qquad \cdot \Prob(F_i^{(p)} \in A \ | \ T_i^{(p)} = 0, \ X_i+M_i^{(p)} \in V, \nonumber \\
  & \hspace{2.5cm} \#\{i: \ (T_i^{(p)}, X_i + M_i^{(p)}) \in \{0\} \times V \} =1) \label{eq:condprob2}
 \end{align}
 mit
 \begin{align*}
   & \Prob(F_i^{(p)} \in A \ | \ T_i^{(p)} = 0, \ X_i+M_i^{(p)} \in V, \\
     & \hspace{2cm} \#\{i: \ (T_i^{(p)}, X_i + M_i^{(p)}) \in \{0\} \times V\} =1)\\
  ={}& \Prob(F_i^{(p)} \in A \ |  \ T_i^{(p)} = 0, \ X_i+M_i^{(p)} \in V, \\
     & \hspace{5cm} \ (T_j^{(p)},X_j+M_j^{(p)}) \notin \{0\} \times V \ \forall j\neq i)  \displaybreak[0]\\
  ={} & \Prob(F_i^{(p)} \in A \ | \ T_i^{(p)} = 0, \ X_i \in V, \ (T_j^{(p)},X_j+M_j^{(p)}) \notin \{0\} \times V \ \forall j\neq i) \\
  ={} & \Prob(F_i^{(p)} \in A \ | \ T_i^{(p)} = 0),
 \end{align*}
 wobei wir ausnutzen, dass $W_i$ stochastisch unabh\"angig von $X_i$, $X_j$ und $W_j$ f\"ur alle $j \neq i$ ist.

 Verwenden wir \eqref{eq:condprob1}, \eqref{eq:condprob2} und
 $$ \sum_{i \in \N} \Prob\left(T_i^{(p)} = 0, \  X_i+M_i^{(p)} \in V \ | \ \Pi(\{0\} \times V \times \R^{p\Z^d}) =1\right) = 1,$$
 so erhalten wir
 \begin{align*}
  \tilde{Q}^{(p)}(A) ={} \Prob(F_i^{(p)} \in A \ | \ T_i = 0) {} = {} \Prob\left( W(\cdot) - \frac{\sigma^2(\cdot)}2 \in A \ | \ T^{(p)} =0\right)
 \end{align*}
 f\"ur alle $A \in \mathcal{B}(\R^{\lambda \Z^d})$.
\end{bew}
\medskip

Beschr\"anken wir uns auf die Annahmen von Lemma \ref{mmm}, so k\"onnen wir $Q$ als die Verteilung von $F_i$
bedingt auf $X_i + M_i$ und $T_i$ beschreiben.
Seien $\Pi = \sum_{i \in \N} \delta_{(T_i,X_i+M_i,F_)i}$ und $E \in \mathcal{B}(\R^d \times \R)$ so, dass
$ \int_{E} e^{-x} (d t \times d x) \in (0,\infty)$.
Weiterhin seien $N = \Pi(E \times \mathcal{C}(\R^d))$ und $i_1, \ldots, i_N$ so, dass
$(T_{i_k}, X_{i_k}+M_{i_k}) \in E$ f\"ur $k=1,\,\ldots,N$.
Mit $G_1, \ldots, G_N$ bezeichnen wir eine zuf\"allige Permutation von $F_{i_1}, \ldots, F_{i_N}$, wobei
jede Permutation die gleiche Wahrscheinlichkeit habe.

\begin{satz} \label{thm:conditional2}
 Bedingt auf $N=n$ sind $G_1, \ldots, G_n$ unabh\"angig identisch gem\"a"s $Q$ verteilt.
\end{satz}
\begin{bew}
 Wir wollen zeigen, dass alle endlich-dimensionalen R\"ander von $G_1, \ldots, G_n$ Produkte der eindimensionalen R\"ander mit Verteilung
 $Q$ sind. Durch Zerlegen der Mengen in $\mathcal{C}$ und Umindizieren gen\"ugt es, $\Prob(G_1 \in A_1, \ldots, G_{n_1} \in A_1, G_{n_1+1} \in A_2, \ldots, G_{n_1+n_2+\ldots+n_l} \in A_l \ | \ N=n)$
 equals $\prod_{i=1}^l \tilde Q(A_i)^{n_i}$
 f\"ur paarweise disjunkte Mengen $A_1, \ldots, A_l \in \mathcal{C}$, $n_1,\ldots, n_l \in \N$ mit $n_1 + \ldots + n_l \leq n$
 zu zeigen.
 Seien $m = n_1 + \ldots + n_l$ und $A = \bigcup_{i=1}^l A_i$. Dann gilt
 \begin{align*}
  & \Prob(G_1 \in A_1, \ldots, G_{n_1} \in A_1, G_{n_1+1} \in A_2, \ldots, G_{m} \in A_l \ | \ N=n)\\
  ={} & \sum_{\substack{k_1\geq n_1, \ldots, k_l\geq n_l \\ k_1+\ldots+k_l\leq n}} \Prob(G_1 \in A_1, \ldots, G_{m} \in A_l
         \ | \ \bigcap_{j=1}^l \Pi(E \times A_j) = k_j, \ N=n)\\
  & \hspace{2.2cm} \cdot \Prob(\Pi(E \times A_j) = k_j, \ j=1,\ldots,l \ | \ N=n) \displaybreak[0]\\
  ={} & \sum_{\substack{k_1\geq n_1, \ldots, k_l\geq n_l \\ k_1+\ldots+k_l\leq n}} \frac{k_1}{n} \cdots
      \frac{k_1-n_1+1}{n-n-1+1} \frac{k_2}{n-n_1} \ldots \frac{k_l-n_l}{n-m+1}\\
  & \hspace{2.2cm}\cdot \binom{n}{k_1, \ldots, \ k_l, \ n-k_1-\ldots-k_l} \cdot \tilde Q(A_1)^{k_1} \cdots Q(A_l)^{k_l}\\
  & \hspace{7.5cm} \cdot (1 - Q(A))^{n-k_1-\ldots-k_l} \displaybreak[0]\\
  ={} & \left( \prod_{i=1}^l Q(A_i)^{n_i}\right) \cdot \Bigg( \sum_{\substack{k_1\geq n_1, \ldots, k_l\geq n_l
                             \\ k_1+\ldots+k_l\leq n}}\binom{n-m}{k_1-n_1, \ldots, \ k_l-n_l, \ n-k_1-\ldots-k_l} \\
   & \hspace{2.8cm} \cdot Q(A_1)^{k_1-n_1} \cdots Q(A_l)^{k_l-n_l} \cdot (1 - Q(A) )^{n-k_1-\ldots-k_l} \Bigg) \displaybreak[0]\\
  ={} & \left( \prod_{i=1}^l Q(A_i)^{n_i}\right) \cdot \Bigg( \sum_{\substack{k_1\geq 0, \ldots, k_l\geq 0
                             \\ k_1+\ldots+k_l\leq n-m}}\binom{n-m}{k_1, \ldots, \ k_l, \ n-m-k_1-\ldots-k_l} \\
   & \hspace{2.8cm} \cdot Q(A_1)^{k_1} \cdots Q(A_l)^{k_l} \cdot (1 - Q(A))^{n-m-k_1-\ldots-k_l} \Bigg) \displaybreak[0]\\
={} &  \prod_{i=1}^l Q(A_i)^{n_i}.
 \end{align*}
\end{bew}
\bigskip

%% file: Fehlerabschaetzungen.tex
\chapter{Simulationsmethoden und Fehlerabsch"atzungen}

Aus der Verteilungsgleichheit der Prozesse $Z(\cdot), Z_1(\cdot), Z_2(\cdot), Z_3(\cdot)$ und $Z_{4,\lambda*}(\cdot)$ erhalten wir f"unf verschiedene Darstellungen eines Brown-Resnick-Prozesses, aus denen sich jeweils unterschiedliche Simulationsmethoden ergeben. Dabei tritt stets ein Poisson-Punktprozess auf $\R$ mit Intensit"atsma"s $\lambda e^{-x} dx, \ \lambda > 0$ auf.
Daher betrachten wir zun"achst, wie man diesen simulieren kann.
\bigskip

\section{Simulation von Poisson-Punktprozessen} 

Sei $\{N(x), x \geq 0 \}$ ein (homogener) Poisson-Prozess auf $\R$ mit Intensit"at $\lambda$, $\lambda > 0$, d.\,h. $N(0)=0$, $\{N(x), x \geq 0 \}$ hat unabh"angige und station"are Zuw"achse und f"ur alle $x, t \geq 0$ ist $N(x+t) -N(x)$ Poisson-verteilt mit Parameter $\lambda t$ (vgl. Definition 2.1.1 in \cite{ross-1996}).
Seien $Y_n, \ n \in \N$ die zugeh"origen \glqq Wartezeiten\grqq{}, d.\,h. $Y_n := \inf \{ x \geq 0: N(x) \geq n \} - \inf \{ x \geq 0: N(x) \geq n-1 \}$. Dann sind $Y_n, \ n \in \N$ nach Proposition 2.2.1 aus \cite{ross-1996} unabh"angig identisch exponentialverteilt mit Parameter $\lambda$.
Seien nun $S_n := \sum_{i=1}^n Y_i, \ n \in \N$, d.\,h. $S_n := \inf \{ x \geq 0: N(x) \geq n \}$. Dann gilt andererseits

\begin{lemma}
Es ist $\Pi = \sum_{j \in \N} \delta_{S_j}$ ein Poisson-Punktprozess auf\ $\R_{\geq 0}$ mit Intensit"at $\lambda$.
\end{lemma}
\begin{bew}
Nach Konstruktion gilt f"ur $0 \leq a \leq b$, dass 
\[ \Pi((a,b]) = \{ \# j: S_j \in (a,b] \} = N(b) - N(a). \]
Die Poisson-Punktprozesseigenschaften folgen nun aus der Definition eines Poissonprozesses:
\begin{enumerate}
\item $\Prob(\Pi((a,b])=n) = \Prob(N(b) - N(a) = n) = e^{-\lambda(b-a)} \frac {(\lambda(b-a))^n}{n!}$.
\item F"ur $A_1, A_2, \ldots \subset \R, \ A_i \cap A_j = \emptyset, \ i \neq j$ sind $\Pi(A_1), \Pi(A_2), \ldots$ stochastisch unabh"angig, da $\{N(x),x \geq 0\}$ stochastisch unabh"angige Zuw"achse hat.
\end{enumerate}
\end{bew}

Die Punkte eines Poisson-Punktprozesses auf $\R$ mit Intensit"at $\lambda$ lassen sich also als Summen von unabh"angig identisch exponentialverteilten Zufallsvariablen mit Parameter $\lambda$ simulieren.
\medskip

Aus diesem Lemma erh"alt man nun folgendes

\begin{coro} \label{PPP-e-x}
Es ist $\widetilde{\Pi} = \sum_{j \in \N} \delta_{\log(\frac 1 {S_j})}$ ein Poisson-Punktprozess auf\ $\R$ mit Intensit"atsma"s $\lambda e^{-x} dx$.
\end{coro}
\begin{bew}
Wir stellen zun"achst fest, dass die Abbildung $$H: \R_{> 0} \rightarrow \R, \ x \mapsto \log\left(\frac 1 x\right)$$ ein Hom"oomorphismus ist. Daher ist f"ur jede Borel-Menge $A \subset \R$ auch $H^{-1}(A)$ eine Borel-Menge in $\R_{\geq 0}$ und es gilt $\widetilde{\Pi}(A) = \Pi(H^{-1}(A))$. Auch sind die Urbilder disjunkter Mengen unter der Abbildung $H$ disjunkt. Somit bleiben die Eigenschaften eines Poisson-Punktprozesses erhalten.
\medskip

Zu berechnen bleibt das Intensit"atsma"s von $\widetilde{\Pi}$.
Seien $a < b \in \R$. Dann gilt
\begin{align*}
&\Prob( \widetilde{\Pi}([a,b]) = 0) = \Prob( \Pi([e^{-b}, e^{-a}]) = 0)\\
 ={}& \exp(- \lambda (e^{-a} - e^{-b})) = \exp\left( - \int_a^b \lambda e^{-x}dx\right).
\end{align*}
Also ist das Intensit"atsma"s $\lambda e^{-x} dx$ -- wie behauptet.
\end{bew}

Einen Poisson-Punktprozess wie in Satz \ref{z2} k"onnen wir mit Hilfe des folgenden Korollars simulieren, welches ein Analogon zu Lemma 3 aus \cite{schlather-2002} ist.
\medskip

\begin{coro} \label{PPP-dt-e-x}
Sei $I \subset \R^d$ ein endliches Intervall, $U_i \sim_{u.i.v.} U(I), \ i \in \N$.
Dann ist $\widehat{\Pi} = \sum_{i \in \N} \delta_{(\log(\frac 1 {S_i}), U_i)}$ ein Poisson-Punktprozess auf \ $\R \times I$ mit Intensit"atsma"s $\frac \lambda {|I|} e^{-x} dx \ dt$.
\end{coro}
\begin{bew}
Folgt direkt aus Korollar \ref{PPP-e-x} und Lemma \ref{2dim-ppp}.
\end{bew}

\section{Fehlerabsch"atzungen}

Im Folgenden werden wir den Brown-Resnick-Prozess gem"a"s $Z(\cdot),$ $Z_1(\cdot),$ $Z_2(\cdot),$ $Z_3(\cdot)$ bzw. $Z_{4,\lambda^*}(\cdot)$ approximieren. F"ur das Berechnen der Approximationsfehler (im Sinne der Wahrscheinlichkeit, dass die Approximation nicht der korrekten Realisierung des Prozesses entspricht) werden sich die folgenden Lemmata als n"utzlich erweisen.

\begin{lemma} \label{fehlerlemma}
Sei $\{W(t), t \in \R\}$ die eindimensionale Standard"=Brownsche"=Bewegung, $W_i(\cdot) \sim_{u.i.v.} W, \ i \in \N$ und $\sum_{i \in \N} \delta_{X_i}$ ein davon unabh"angiger Poisson"=Punktprozess auf \ $\R$ mit Intensit"atsm"a"s $\lambda e^{-x} dx, \ \lambda > 0$. Weiterhin seien $t_u \leq 0 \leq t_o$, $C,x \in \R$, $C>x$. Dann gilt

\begin{align*}
&\Prob\left(\exists i \in \N: X_i < x, \sup_{t \in [t_u,t_o]} X_i + W_i(t) > C\right) \\
\leq {} &  \lambda e^{-C} \frac {2|t_u|} {C-x} \exp\left(\frac {|t_u|} 2\right) \left(1- \Phi\left( \frac{C-x-|t^*|}{\sqrt{|t_u|}}\right) \right)\\
&+ \lambda e^{-C} \frac {2|t_o|} {C-x} \exp\left(\frac {|t_o|} 2\right) \left(1- \Phi\left( \frac{C-x-|t^*|}{\sqrt{|t_o|}}\right) \right).
\end{align*}
\end{lemma}

\begin{bew}
Wir betrachten den Poisson-Punktprozess $\Pi := \sum_{i \in \N} \delta_{(X_i, X_i + W_i(\cdot))}$ auf $\R \times C(\R^d)$ mit Intensit"atsma"s $$\Lambda(A_1 \times A_2) = \int_{A_1} e^{-y} \Prob(y+W(\cdot) \in A_2) dy$$ f"ur $A_1 \times A_2 \in \mathcal{B} \times \mathcal{C}$ (vgl. Kapitel 2), wobei $\Prob$ das zu $W(\cdot)$ geh"orende Wahrscheinlichkeitsma"s sei.
\medskip

F"ur $t^* \in \R$ definieren wir nun die Menge 
\[A_{t^*} := \left\{ (X,W(\cdot)) \in \R \times C(\R^d): \ X < x, \sup_{t \in [0,|t^*|]} X + W(t) > C\right\}.\]
Dann gilt:
\begin{align*}
& \Lambda(A_{t^*}) = \int_{-\infty}^x \lambda e^{-y} \frac 2 {\sqrt{\pi}} \int_{\frac{C-y}{\sqrt{2|t^*|}}}^\infty e^{-z^2} dz \ dy  \hspace{0.6cm} \text{nach Formel 1.1.1.4 aus \cite{borodin-2002}}\\
\leq{} & \int_{-\infty}^x \lambda e^{-y} \frac 2 {\sqrt{\pi}} \frac {\sqrt{2|t^*|}}{2(C-y)} \int_{\frac{C-y}{\sqrt{2|t^*|}}}^\infty 2ze^{-z^2} dz \ dy\\
={} & \int_{-\infty}^x \lambda e^{-y} \frac 1 {\sqrt{\pi}} \frac {\sqrt{2|t^*|}}{C-y} \exp\left(-\frac{(C-y)^2}{2|t^*|}\right) dy \\
={} & \int_{C-x}^{\infty} \lambda e^{-C} e^{y} \frac 1 {\sqrt{\pi}} \frac {\sqrt{2|t^*|}}{y} \exp\left(-\frac{y^2}{2|t^*|}\right) dy\\
\leq{} & \frac {\lambda e^{-C}} {\sqrt{\pi}} \frac {\sqrt{2|t^*|}}{C-x} \int_{C-x}^\infty  \exp\left(-\frac{y^2}{2|t^*|}+\frac{2|t^*|y}{2|t^*|}\right) dy\\
={} & \lambda e^{-C} \frac {2|t^*|} {C-x} \exp\left(\frac {|t^*|} 2\right) \int_{C-x}^\infty \frac 1 {\sqrt{2|t^*|\pi}} \exp\left(-\frac{y^2}{2|t^*|}+\frac{2|t^*|y}{2|t^*|}-\frac{|t^*|^2}{2|t^*|}\right) dy \\
={} & \lambda e^{-C} \frac {2|t^*|} {C-x} \exp\left(\frac {|t^*|} 2\right) \int_{C-x}^\infty \frac 1 {\sqrt{2|t^*|\pi}} \exp\left(-\frac{(y-|t^*|)^2}{2|t^*|}\right) dy\\
={} & \lambda e^{-C} \frac {2|t^*|} {C-x} \exp\left(\frac {|t^*|} 2\right) \left(1- \Phi\left( \frac{C-x-|t^*|}{\sqrt{|t^*|}}\right) \right).
\end{align*}

Damit erh"alt man die Ungleichung
\begin{align*}
& \Prob\left(\exists i \in \N: x_i < x, \sup_{t \in [0,|t^*|]} X_i + W_i(t) > C\right)\\
={}& \Prob\left(\Pi(A_{|t^*|}) >0\right) = 1 - \exp\left(-\Lambda(A_{|t^*|})\right)\\
\leq{} & \Lambda(A_{|t^*|}), \hspace{3cm} \text{da } e^{-x} \geq 1 - x  \quad \forall x > 0\\
\leq{} & \lambda e^{-C} \frac {2|t^*|} {C-x} \exp\left(\frac {|t^*|} 2\right) \left(1- \Phi\left( \frac{C-x-|t^*|}{\sqrt{|t^*|}}\right) \right).\\
\end{align*}
Weiterhin ist
\begin{align*}
&\Prob\left(\exists i \in \N: X_i < x, \sup_{t \in [t_u,t_o]} X_i + W_i(t) > C\right) \\
\leq{} & \Prob\left(\exists i \in \N: X_i < x, \sup_{t \in [0,|t_u|]} X_i + W_i(t) > C\right) \\
& + \Prob\left(\exists i \in \N: X_i < x, \sup_{t \in [0,|t_o|]} X_i + W_i(t) > C\right)
\end{align*}
und somit folgt die behauptete Ungleichung.
\end{bew}
\medskip

Dieses Lemma erm"oglicht uns eine Aussage "uber das Verhalten von \[ \Prob\left(\exists i \in \N: X_i < x, \sup_{t \in [0,|t^*|]} X_i + W_i(t) > C\right)\] f"ur gro"se $t^*$.
\medskip

\begin{coro}
Unter den Voraussetzungen von Lemma \ref{fehlerlemma} gibt es $L_1, L_2, t_0 >0$ so, dass f"ur alle $|t^*| > t_0$ gilt:
\begin{align*}
&L_1 \exp\left(\frac {|t^*|} 2\right) \\
\leq{} & - \log\left(1- \Prob\left(\exists i \in \N: X_i < x, \sup_{t \in [0,|t^*|]} X_i + W_i(t) > C\right)\right)\\
\leq{} & L_2 |t^*| \exp\left(\frac {|t^*|} 2\right).
\end{align*}
\end{coro}
\begin{bew}
Seien $\Lambda, A_{t^*}$ wie im Beweis von Lemma \ref{fehlerlemma}. Dann gilt nach eben diesem:
\[\Lambda(A_{t^*}) \leq \lambda e^{-C} \frac {\sqrt{2|t^*|^3}} {C-x} \exp\left(\frac {|t^*|} 2\right) \Bigg(1 + \frac{1}{C-x} \exp\left(-\frac{(C-x)^2}{2|t^*|}\right)\Bigg).\]

Andererseits ist
\[\Lambda(A_{t^*}) = \int_{-\infty}^x \lambda e^{-y} \frac 2 {\sqrt{\pi}} \int_{\frac{C-y}{\sqrt{2|t^*|}}}^\infty e^{-z^2} dz \ dy  \hspace{1cm} \text{nach Formel 1.1.1.4 aus \cite{borodin-2002}}. \]

Definieren wir $\Psi(t) := \int_t^\infty \varphi(s) ds$ mit $\varphi(t) := \frac 1 {\sqrt{2\pi}} \exp\left(-\frac {s^2} 2\right), \ t \geq 0$, so gilt
\begin{align*}
\Psi(t) ={} & \int_t^\infty \frac{1}{\sqrt{2\pi}} \frac s s e^{-\frac{s^2}2} ds\\
={} & \frac{1}{\sqrt{2\pi}} \frac 1 t e^{- \frac {t^2} 2} - \int_t^\infty \frac{1}{\sqrt{2\pi}} \frac s {s^3} e^{-\frac{s^2}2} ds\\
={} & \frac{1}{\sqrt{2\pi}} \frac 1 t e^{- \frac {t^2} 2} - \left[ -\frac{1}{\sqrt{2\pi}} \frac 1 {s^3} e^{- \frac {s^2} 2} \right]_{s=t}^\infty + \int_t^\infty \frac{1}{\sqrt{2\pi}} \frac 3 {s^4} e^{-\frac{s^2}2} ds\\
\geq{} & \left( \frac 1 t - \frac 1 {t^3} \right) \varphi(t) \hspace{4cm} \forall t > 0.
\end{align*}

(Eine pr"azisere Absch"atzung, die man durch weitere Iterationsschritte der partiellen Integration erh"alt, findet sich im Beweis von Satz 4 aus \cite{ballani-2009}.)

Wegen $ t^{-3} = o(t^{-1})$ gibt es also ein $K_1 > 0$ so, dass $\Psi(t) \geq \frac {K_1} t \exp\left(-\frac {t^2} 2\right)$ f"ur alle $t \geq 1$. Da andererseits $\frac {\Psi(t)} {\sqrt{2\pi} \varphi(t)}$ eine stetige Funktion auf $\R_{\geq 0}$ und $[0,1]$ kompakt ist, existiert $K_2 := \min_{t \in [0,1]} \frac {\Psi(t)} {\sqrt{2\pi} \varphi(t)} > 0$ und somit folgt f"ur alle $t > 0$ die Ungleichung
\[ \Psi(t) \geq \min \left(K_2, \frac {K_1} t \right) \exp\left(-\frac {t^2} 2\right). \]

Damit erhalten wir f"ur $|t^*| > 1$ die Absch"atzung:
\begin{align*}
\Lambda(A_{t^*}) ={} & \int_{-\infty}^x \lambda e^{-y} 2  \int_{\frac{C-y}{\sqrt{|t^*|}}}^\infty \frac 1 {\sqrt{2\pi}} e^{-\frac {z^2} 2} dz \ dy  \\
\geq{} & \int_{-\infty}^x \lambda e^{-y} 2  \min\left(K_2, \frac{K_1 \sqrt{|t^*|}}{C-y}\right) \exp\left(-\frac {(C-y)^2} {2|t^*|}\right) \ dy\\
={} & 2\lambda e^{-C} \int_{C-x}^\infty e^{y} \min\left(K_2, \frac{K_1 \sqrt{|t^*|}}{y}\right) \exp\left(-\frac {y^2} {2|t^*|}\right) \ dy\\
={} & 2\lambda e^{-C} e^{\frac {|t^*|} 2} \int_{C-x}^\infty \min\left(K_2, \frac{K_1 \sqrt{|t^*|}}{y}\right) \exp\left(-\frac {(y-|t^*|)^2} {2|t^*|}\right) \ dy\\
={} & 2\lambda e^{-C} e^{\frac {|t^*|} 2} \int_{\frac{C-x} {\sqrt{|t^*|}} - \sqrt{|t^*|}}^\infty  \min\left(K_2\sqrt{|t^*|}, \frac{K_1 |t^*|}{y\sqrt{|t^*|} + |t^*|}\right) \exp\left(-\frac {y^2} 2\right) \ dy
\end{align*}
\begin{align*}
\geq{} & 2\lambda e^{-C} e^{\frac {|t^*|} 2} \int_{\frac{C-x} {\sqrt{|t^*|}} - 1}^{\frac{C-x} {\sqrt{|t^*|}}} \sqrt{|t^*|} \min\left(K_2, \frac{K_1 \sqrt{|t^*|}}{y\sqrt{|t^*|} + |t^*|}\right) \exp\left(-\frac {y^2} 2\right) \ dy\\
\geq{} & 2\lambda e^{-C} \exp\left(\frac {|t^*|} 2\right) \min\left(K_2\sqrt{|t^*|}, \frac{K_1 |t^*|}{C - x + |t^*|}\right)\\
& \hspace{4cm} \cdot \exp\left(-\frac {\max\{(C-x)^2,(C-x-\sqrt{|t^*|})^2\}} {2|t^*|} \right)\\
\geq{} &  \frac{2\lambda e^{-C} K_1 |t^*|}{C - x + |t^*|} \exp\left(\frac {|t^*|} 2\right) \exp\left(-\frac {\max\{(C-x)^2,(C-x-\sqrt{|t^*|})^2\}} {2|t^*|} \right)
\end{align*}
f"ur $|t^*|$ gro"s genug.

Mit $\exp\left(-\frac {(C-x)^2} {2|t^*|}\right) \stackrel{|t^*| \to \infty}{\longrightarrow} 1$ bzw. $\exp\left(-\frac {((C-x)-\sqrt{|t^*|})^2} {2|t^*|}\right) \stackrel{|t^*| \to \infty}{\longrightarrow} \exp(-\frac 1 2)$ und
$$\Prob\left(\exists i \in \N: X_i < x, \sup_{t \in [0,|t^*|]} X_i + W_i(t) > C\right) = 1 - \exp(-\Lambda(A_{t^*}))$$ folgt die Behauptung.
\end{bew}
\bigskip

\begin{lemma} \label{BBsup}
Seien $W_i^{(j)} \sim W, \ i \in \N, \ j \in \N,$ unabh"angig identisch verteilte Standard"=Brownsche"=Bewegungen auf \ $\R$. Seien weiterhin $\sum_{i \in \N} \delta_{X_i^{(j)}}$ f"ur $j \in \N$ stochastisch unabh"angige Poisson"=Punktprozesse auf \ $\R$ mit Intensit"atma"s $e^{-x} dx$ sowie $0 < a_1 < a_2 < a_3 < \ldots$ "aquidistante Punkte mit Abstand $h > 0$ und $a_1 > 4+h, \ a_1 \in h\Z$, sowie $L > 0$, $L \in h\Z$, $C \in \R$. Dann gilt:
\begin{align*}
& \Prob\left( \exists i,j \in \N: C < X_i^{(j)} + \sup_{t \in [a_j, a_j + L]\cap h\Z} \left( W_i^{(j)}(t) - \frac t 2 \right) \leq X_i^{(j)} + W_i^{(j)}(0) \right) \\
\leq{} & \frac {8e^{-C}} h \Bigg( \exp\left(-\frac {\sqrt{a_1-h}} 2 \right)\\
 & + \frac {25L} {9} \sqrt{\frac{2L}\pi} \exp\left(-\frac{(\frac 3 {10} \sqrt{a_1-h} + \frac 1 5 \sqrt{(a_1-h)\vee25} -\log(\frac{\sqrt{(a_1-h)\vee 25}} 2))^2} {2L}\right) \\
& \hspace{1cm} + \frac 1 {\sqrt{2\pi}} \left(1 + \frac 1 {\sqrt{a_1-h}-2}\right) \exp\left(- \frac{(\sqrt{a_1-h}-2)^2} 8\right) \Bigg).
\end{align*} 
\end{lemma}
\begin{bew}
Zun"achst definieren wir "ur $j \in \N$ die Poisson"=Punktprozesse \[\Pi_j := \sum_{i \in \N} \delta_{X_i^{(j)} + W_i^{(j)}}\] auf $C(\R^d)$ sowie die Mengen \[A_j := \left\{ x+w(\cdot) \in C(\R^d): C < x + \sup_{t \in [a_j, a_j + L]\cap h\Z} \left( w(t) - \frac t 2 \right) \leq x + w(0) \right\}.\]
Sei weiterhin $\Prob$ das zu $W(\cdot)$ geh"orende Wahrscheinlichkeitsma"s.

Dann hat $\Pi_j, \ j \in \N,$ das Intensit"atsma"s
$\Lambda_j(A) := \int_{\R} e^{-x} \Prob(W \in A-x) dx$
 f"ur $A \in {\cal C}$ (vgl. Kapitel 2).

Wir wollen nun $\Lambda_j(A_j)$ nach oben hin absch"atzen.
Es gilt:
\begin{align*}
&\Lambda_j(A_j) = \int_{\R} e^{-x} \Prob\left(C - x < \sup_{t \in [a_j, a_j+L]\cap h\Z} \left(W(t) - \frac t 2\right) \leq W(0)\right) dx \\
={} & \int_{\R} e^{-x} \Prob\left(C - x < \sup_{t \in [a_j, a_j+L] \cap h\Z} \left(W(t) - \frac t 2\right) \leq  0\right) dx,\\
& \hspace{8cm} \text{ da } W(0) = 0 \text{ $\Prob$-f.s.} \\
\leq{} & \int_{C + \frac {\sqrt{a_j}} 2 +  \log(\frac {\sqrt{a_j}} 2)}^\infty e^{-x} dx \\
& + \int_C^{C + \frac {\sqrt{a_j}} 2 +  \log(\frac {\sqrt{a_j}} 2)} e^{-x} \Prob\left(C - x < \sup_{t \in [a_j, a_j+L] \cap h\Z} \left(W(t) - \frac t 2\right) \leq  0\right) dx\displaybreak[0]\\
={} & e^{-C} \frac 2 {\sqrt{a_j}} \exp\left(-\frac {\sqrt{a_j}} 2 \right) \\
  & + e^{-C} \int_{- \frac {\sqrt{a_j}} 2 -  \log(\frac {\sqrt{a_j}} 2)}^0 e^{x} \Prob\left(x < \sup_{t \in [a_j, a_j+L] \cap h\Z} \left(W(t) - \frac t 2\right) \leq  0\right) dx\displaybreak[0]\\
\leq{} & e^{-C} \frac 2 {\sqrt{a_j}} \exp\left(-\frac {\sqrt{a_j}} 2 \right) \\
  & + e^{-C} \int_{- \frac {\sqrt{a_j}} 2 -  \log(\frac {\sqrt{a_j}} 2)}^0 e^{x} \Prob\Bigg(0 \geq W(a_j) - \frac {a_j} 2  > - {\sqrt{a_j}}, \\
  & \hspace{6.7cm} \sup_{t \in [a_j, a_j+L] \cap h\Z} \left(W(t) - \frac t 2\right) > x\Bigg) dx\\
  & + e^{-C} \int_{- \frac {\sqrt{a_j}} 2 -  \log(\frac {\sqrt{a_j}} 2)}^0 e^{x} \Prob\Bigg(W(a_j) - \frac {a_j} 2 \leq -{\sqrt{a_j}}, \\
  & \hspace{6.7cm} \sup_{t \in [a_j, a_j+L] \cap h\Z} \left(W(t) - \frac t 2\right) > x\Bigg) dx\displaybreak[0]\\
\leq{} & \frac {2 e^{-C}} {\sqrt{a_j}} \exp\left(-\frac {\sqrt{a_j}} 2 \right)\\
  & + e^{-C} \int_{- \frac {\sqrt{a_j}} 2 -  \log(\frac {\sqrt{a_j}} 2)}^0 e^{x} \Prob\left(0 \geq W(a_j) - \frac {a_j} 2 > -{\sqrt{a_j}}\right) dx\\
  & + e^{-C} \int_{- \frac {\sqrt{a_j}} 2 -  \log(\frac {\sqrt{a_j}} 2)}^0 e^{x} \Prob\Bigg(W(a_j) - \frac {a_j} 2 \leq -{\sqrt{a_j}},\\
  & \hspace{5.5cm} \sup_{t \in [0,L]} \left(\tilde{W}(t) - \frac t 2\right) > x-W(a_j) + \frac {a_j} 2\Bigg) dx,\displaybreak[0]\\
  & \hspace{2.5cm} \text{wobei $\tilde{W}$ eine von $W$ unabh"angige Brownsche Bewegung sei.}
\end{align*}

Nun verwenden wir, dass $W(a_j) - \frac {a_j} 2 \sim {\cal N}\left(-\frac {a_j} 2, a_j\right)$ ist und erhalten damit

\begin{align}
&\Lambda_j(A_j) \leq{} e^{-C} \frac 2 {\sqrt{a_j}} \exp\left(-\frac {\sqrt{a_j}} 2 \right) \nonumber \\
  & + e^{-C} \int_{- \frac {\sqrt{a_j}} 2 -  \log(\frac {\sqrt{a_j}} 2)}^0 e^{x} \int_{-{\sqrt{a_j}}}^0 \frac 1 {\sqrt{2 \pi a_j}} \exp\left(-\frac {(y + \frac {a_j} 2)^2} {2a_j} \right) dy \ dx \nonumber \\
  & + e^{-C} \int_{- \frac {\sqrt{a_j}} 2 -  \log(\frac {\sqrt{a_j}} 2)}^0 e^{x} \int_{-\infty}^{-{\sqrt{a_j}}} \frac 1 {\sqrt{2 \pi a_j}} \exp\left(-\frac {(y + \frac {a_j} 2)^2} {2a_j} \right) \nonumber \\
  & \hspace{6.2cm} \Prob\left(\sup_{t \in [0,L]} \left(\tilde{W}(t) - \frac t 2\right) > x-y\right) dy \ dx \nonumber\\
\leq{} & \frac {2e^{-C}} {\sqrt{a_j}} \exp\left(-\frac {\sqrt{a_j}} 2 \right) + e^{-C} \int_{- \frac {\sqrt{a_j}} 2 -  \log(\frac {\sqrt{a_j}} 2)}^0 e^{x} dx \int_{-1 + \frac {\sqrt{a_j}} 2}^{\frac {\sqrt{a_j}} 2} \frac 1 {\sqrt{2 \pi}} \exp\left(-\frac {y^2} {2} \right) dy  \nonumber \\
  & + e^{-C} \int_{- \frac {\sqrt{a_j}} 2 -  \log(\frac {\sqrt{a_j}} 2)}^0 e^{x} \int_{-\infty}^{-{\sqrt{a_j}}} \frac 1 {\sqrt{2 \pi a_j}} \exp\left(-\frac {(y + \frac {a_j} 2)^2} {2a_j} \right) \nonumber \\
  & \hspace{6.5cm} \Prob\left(\sup_{t \in [0,L]} W(t) > x-y\right) dy \ dx \nonumber \\
\leq{} & e^{-C} \frac 2 {\sqrt{a_j}} \exp\left(-\frac {\sqrt{a_j}} 2 \right) \nonumber\\
   & + e^{-C} \left(1 - \exp\left(-\frac{\sqrt{a_j}} 2 - \log\left(\frac{\sqrt{a_j}}2\right)\right)\right) \cdot \int_{-1 + \frac {\sqrt{a_j}} 2}^{\frac {\sqrt{a_j}} 2} \frac 1 {\sqrt{2 \pi}} \exp\left(-\frac {y^2} {2} \right) dy \nonumber \\
  & + e^{-C} \int_{- \frac {\sqrt{a_j}} 2 -  \log(\frac {\sqrt{a_j}} 2)}^0 e^{x} \int_{-\infty}^{-{\sqrt{a}}} \frac 1 {\sqrt{2 \pi a_j}} \exp\left(-\frac {(y + \frac {a_j} 2)^2} {2a_j} \right) \nonumber \\
  & \hspace{2.5cm} \cdot \frac 2 {\sqrt{\pi}} \int_{\frac{x-y}{\sqrt{2L}}}^\infty \exp(-z^2) dz \ dy \ dx  \text{\ nach Formel 1.1.1.4 aus \cite{borodin-2002}} \nonumber \\
\leq{} & e^{-C} \frac 2 {\sqrt{a_j}} \exp\left(-\frac {\sqrt{a_j}} 2 \right) + e^{-C} \int_{-1 + \frac {\sqrt{a_j}} 2}^{\frac {\sqrt{a_j}} 2} \frac 1 {\sqrt{2 \pi}} \exp\left(-\frac {y^2} {2} \right) dy \nonumber \\
  & + e^{-C} \int_{- \frac {\sqrt{a_j}} 2 -  \log(\frac {\sqrt{a_j}} 2)}^0 e^{x} \int_{-\infty}^{-{\sqrt{a_j}}} \frac 1 {\sqrt{2 \pi a_j}} \exp\left(-\frac {(y + \frac {a_j} 2)^2} {2a_j} \right) \nonumber \\
  & \hspace{3cm} \cdot \frac {\sqrt{2L}} {\sqrt{\pi}(x-y)} \exp\left(-\frac {(x-y)^2}{2L}\right)  dy \ dx \label{eq:drittersum}
\end{align}

F"ur den dritten Summanden aus \eqref{eq:drittersum} l"asst sich weiterhin absch"atzen:
\begin{align*}
& e^{-C} \int_{- \frac {\sqrt{a_j}} 2 -  \log(\frac {\sqrt{a_j}} 2)}^0 e^{x} \int_{-\infty}^{-{\sqrt{a_j}}} \frac 1 {\sqrt{2 \pi a_j}} \exp\left(-\frac {(y + \frac {a_j} 2)^2} {2a_j} \right) \\
  & \hspace{3cm} \cdot \frac {\sqrt{2L}} {\sqrt{\pi}(x-y)} \exp\left(-\frac {(x-y)^2}{2L}\right)  dy \ dx \displaybreak[0]\\
\leq{} &  e^{-C} \int_{- \frac {\sqrt{a_j}} 2 -  \log(\frac {\sqrt{a_j}} 2)}^0 e^{x} \int_{-\infty}^{-{\sqrt{a_j}}} \frac 1 {\sqrt{a_j}} \exp\left(-\frac {(y + \frac {a_j} 2)^2} {2a_j} \right) \\
  & \hspace{3cm} \cdot \frac 1 {\pi} \frac {\sqrt{L}} {\frac {\sqrt{a_j}} 2 -  \log(\frac {\sqrt{a_j}} 2)} \exp\left(-\frac {(x-y)^2}{2L}\right)  dy \ dx, \\
  & \hspace{2cm} \text{da $x-y \geq - \frac {\sqrt{a_j}} 2 -  \log\left(\frac {\sqrt{a_j}} 2\right) + \sqrt{a_j} = \frac {\sqrt{a_j}} 2 -  \log\left(\frac {\sqrt{a_j}} 2\right)$ ist}\displaybreak[0]\\
\leq{} &  e^{-C} \int_{- \frac {\sqrt{a_j}} 2 -  \log(\frac {\sqrt{a_j}} 2)}^0 e^{x} \int_{-\infty}^{-{\sqrt{a_j}}} \frac{\sqrt{L}}{\pi\sqrt{a_j}} \exp\left(-\frac {(y + \frac {a_j} 2)^2} {2a_j} \right)\\
& \hspace{6.5cm} \cdot \exp\left(-\frac {(x-y)^2}{2L}\right)  dy \ dx,\displaybreak[0]\\
& \hspace{1cm} \text{da f"ur alle \ $x > 1$ \ die Ungleichung \ $x - \log(x) > 1 - \log(1) = 1$ \ gilt}\\
& \hspace{1cm} \text{und $\frac{\sqrt{a_j}} 2 > \frac {\sqrt{4}} 2 = 1$ ist} \displaybreak[0]\\
\leq{} &  e^{-C} \int_{- \frac {\sqrt{a_j}} 2 -  \log(\frac {\sqrt{a_j}} 2)}^0 e^{x} dx \int_{-\infty}^{-{\sqrt{a_j}}} \frac 1 {\sqrt{a_j}} \exp\left(-\frac {(y + \frac {a_j} 2)^2} {2a_j} \right) dy \\
& \hspace{6cm} \cdot \frac{\sqrt{L}}{\pi} \exp\left(-\frac {(\frac {\sqrt{a_j}} 2  -  \log(\frac {\sqrt{a_j}} 2))^2}{2L}\right).
\end{align*}

Nun ist die Funktion $\frac 2 5 x -  \log(x)$ auf $\R$ minimal an der Stelle $x=\frac 5 2$. Damit ist $\frac 2 5 x - \log(x) > 0$ f"ur alle $x \in \R$, denn es ist $\frac 2 5 \frac 5 2 = \log(e) > \log\left(\frac 5 2\right)$, und diese Funktion ist isoton auf $\R_{\geq \frac 5 2}$. Daraus folgt
 
\begin{align*}
& \left(\frac {\sqrt{a_j}} 2 - \log\left(\frac{\sqrt{a_j}} 2\right)\right)^2 = \left(\frac {3\sqrt{a_j}} {10} + \left(\frac {2\sqrt{a_j}} {10} - \log\left(\frac {\sqrt{a_j}} 2\right)\right)\right)^2\displaybreak[0]\\
={} & \frac 9 {100} a_j +  \frac 3 5 \sqrt{a_j} \left(\frac {2\sqrt{a_j}} {10} - \log\left(\frac{\sqrt{a_j}} 2\right)\right) + \left(\frac {\sqrt{a_j}} 2 - \log\left(\frac{\sqrt{a}} 2\right)\right)^2\\
\geq{} & \frac 9 {100} a_j  +  \frac 3 5 \sqrt{a_1-h} \left(\frac 2 5 \left(\frac{\sqrt{a_1-h}} {2} \vee \frac 5 2\right)  - \log\left(\frac{\sqrt{a_1-h}} {2} \vee \frac 5 2\right)\right) \\
& + \left(\frac 2 5 \left(\frac{\sqrt{a_1-h}} {2} \vee \frac 5 2\right)  - \log\left(\frac{\sqrt{a_1-h}} {2} \vee \frac 5 2\right)\right) ^2 \displaybreak[0]\\
={} & \frac 9 {100} a_j +  \underbrace{\frac 3 5 \sqrt{a_1-h} \left(\frac {\sqrt{(a_1-h)\vee 25}} 5  - \log\left(\frac{\sqrt{(a_1-h)\vee 25}} 2\right)\right)}_{K_{a_1}^{(1)}} \\
& + \underbrace{\left(\frac {\sqrt{(a_1-h)\vee 25}} 5 - \log\left(\frac{\sqrt{(a_1-h)\vee 25}} 2\right)\right)^2}_{=:K_{a_1}^{(2)}}
\end{align*}
f"ur alle $j \in \N$ und man erh"alt mit $K_{a_1} := K_{a_1}^{(1)} + K_{a_1}^{(2)}$ die Ungleichung
\begin{align*}
&  e^{-C} \int_{- \frac {\sqrt{a_j}} 2 -  \log(\frac {\sqrt{a_j}} 2)}^0 e^{x} dx \int_{-\infty}^{-{\sqrt{a}}} \frac 1 {\sqrt{a_j}} \exp\left(-\frac {(y + \frac {a_j} 2)^2} {2a_j} \right) dy \\
& \hspace{6cm} \cdot \frac{\sqrt{L}} {\pi} \exp\left(-\frac {(\frac{\sqrt{a_j}} 2 -  \log(\frac{\sqrt{a_j}}2))^2}{2L}\right)\displaybreak[0]\\
\leq{} &  e^{-C} \left(1 - \exp\left(- \frac {\sqrt{a_j}} 2 -  \log\left(\frac {\sqrt{a_j}} 2\right)\right)\right) \int_{-\infty}^{\infty} \frac {\sqrt{2\pi}} {\sqrt{2\pi a_j}} \exp\left(-\frac {(y + \frac {a_j} 2)^2} {2a_j} \right) dy \\
& \hspace{1cm}\cdot \frac{\sqrt{L}}{\pi} \exp\left(-\frac {\frac 9 {100} a_j} {2L}\right) \cdot \exp\left(-\frac{K_{a_1}} {2L}\right)\displaybreak[0]\\
\leq{} &  e^{-C} \sqrt{\frac{2L}{\pi}} \exp\left(-\frac {\frac 9 {100} a_j} {2L}\right) \cdot \exp\left(-\frac{K_{a_1}} {2L}\right).
\end{align*}

Einsetzen in \ref{eq:drittersum} liefert schlie"slich
\begin{align}
\Lambda_j(A_j) \leq{} & e^{-C} \frac 2 {\sqrt{a_j}} \exp\left(-\frac {\sqrt{a_j}} 2 \right) + e^{-C} \int_{-1 + \frac {\sqrt{a_j}} 2}^{\frac {\sqrt{a_j}} 2} \frac 1 {\sqrt{2 \pi}} \exp\left(-\frac {y^2} {2} \right) dy \nonumber \\
  & + e^{-C} \sqrt{\frac{2L}{\pi}} \exp\left(-\frac {\frac 9 {100} a_j} {2L}\right) \exp\left(-\frac{K_{a_1}} {2L}\right) \label{eq:mass}.
\end{align}

Nun ist
\begin{align*}
& \Prob\left( \exists i,j \in \N: C < X_i^{(j)} + \sup_{t \in [a_j, a_j + L]} \left( W_i^{(j)}(t) - \frac t 2 \right) \leq X_i^{(j)} + W_i^{(j)}(0) \right) \displaybreak[0]\\
={} & \Prob\left( \exists j \in \N: \Pi_j(A_j) > 0 \right) ={}  1 - \exp\left(- \sum_{j \in \N} \Lambda_j(A_j)\right) \displaybreak[0]\\
\leq{} & \sum_{j \in \N} \Lambda_j(A_j), \qquad \text{da } 1 - e^{-x} \leq x \text{ f"ur alle } x > 0 \text{ gilt.} \displaybreak[0]\\
\leq{} & \sum_{j \in \N} \Bigg(e^{-C} \frac 2 {\sqrt{a_j}} \exp\left(-\frac {\sqrt{a_j}} 2 \right) + e^{-C} \int_{-1 + \frac {\sqrt{a_j}} 2}^{\frac {\sqrt{a_j}} 2} \frac 1 {\sqrt{2 \pi}} \exp\left(-\frac {y^2} {2} \right) dy \nonumber \\
  & \hspace{1cm}+\sqrt{\frac{2L}{\pi}} \exp\left(-\frac {\frac 9 {100} a_j} {2L}\right) \exp\left(-\frac{K_{a_1}}{2L}\right) \Bigg) \qquad \text{nach \eqref{eq:mass}}\displaybreak[0]\\
\leq{} & \int_{a_1-h}^\infty \frac 1 h \Bigg(e^{-C} \frac 2 {\sqrt{a}} \exp\left(-\frac {\sqrt{a}} 2 \right) + e^{-C} \int_{-1 + \frac {\sqrt{a}} 2}^{\frac {\sqrt{a}} 2} \frac 1 {\sqrt{2 \pi}} \exp\left(-\frac {y^2} {2} \right) dy \nonumber \\
  & \hspace{2cm}+\exp\left(-\frac{K_{a_1}} {2L}\right) \sqrt{\frac{2L}{\pi}} \exp\left(-\frac {\frac 9 {100} a} {2L}\right) \Bigg) da, \\
  & \hspace{4cm} \text{da der Integrand streng monoton fallend ist}
\end{align*}
\begin{align*}
\leq{} & e^{-C} \frac 8 h \exp\left(-\frac {\sqrt{a_1-h}} 2 \right)\\
  & + e^{-C} \frac 1 h \int_{-1 + \frac {\sqrt{a_1-h}} 2}^\infty \int_{a_1-h}^\infty \mathbf{1}_{y \in [-1 + \frac {\sqrt{a}} 2, \frac {\sqrt{a}} 2]} \frac 1 {\sqrt{2 \pi}} \exp\left(-\frac {y^2} {2} \right) da \ dy \\
  & +  e^{-C} \exp\left(-\frac{K_{a_1}} {2L}\right) \frac {200L} {9h} \sqrt{\frac{2L}{\pi}} \exp\left(-\frac {\frac 9 {100} (a_1 - h)} {2L}\right)\displaybreak[0]\\
\leq{} &  \frac {e^{-C}} h \Bigg( 8 \exp\left(-\frac {\sqrt{a_1-h}} 2 \right) +  \int_{-1 + \frac {\sqrt{a_1-h}} 2}^\infty (8y+4) \frac 1 {\sqrt{2 \pi}} \exp\left(-\frac {y^2} {2} \right) dy \\
  & + \frac {200L} {9} \sqrt{\frac{2L}{\pi}} \\
  & \hspace{0.8cm} \cdot \exp\left(-\frac{(\frac 3 {10} \sqrt{a_1-h} + \frac 1 5 \sqrt{(a_1-h)\vee25} -\log(\frac{\sqrt{(a_1-h)\vee 25}} 2))^2} {2L}\right),\displaybreak[0]\\
 &\hspace{4cm} \text{da } y \in \left[- 1 + \frac {\sqrt{a}} 2, \frac {\sqrt{a}} 2\right] \iff a \in [4y^2,4(y+1)^2]\displaybreak[0]\\
\leq{} & \frac {8e^{-C}} h \Bigg( \exp\left(-\frac {\sqrt{a_1-h}} 2 \right)\\
  & + \frac {25L} {9} \sqrt{\frac{2L}{\pi}} \exp\left(-\frac{(\frac 3 {10} \sqrt{a_1-h} + \frac 1 5 \sqrt{(a_1-h)\vee25} -\log(\frac{\sqrt{(a_1-h)\vee 25}} 2))^2} {2L}\right)\\
  & \hspace{1cm} + \frac 1 {\sqrt{2\pi}} \left(1 + \frac 1 {\sqrt{a_1-h}-2}\right) \exp\left(- \frac{(\sqrt{a_1-h}-2)^2} 8\right)\Bigg)
\end{align*}
\end{bew}
\bigskip

Versch"arfen wir die Voraussetzungen ein wenig, so erhalten wir

\begin{lemma} \label{BBsup2}
Seien $W_i^{(j)} \sim W, \ i \in \N, \ j \in \N,$ unabh"angig identisch verteilte Standard"=Brownsche"=Bewegungen auf \ $\R$. Seien weiterhin $\sum_{i \in \N} \delta_{X_i^{(j)}}$ f"ur $j \in \N$ stochastisch unabh"angige Poisson"=Punktprozesse auf \ $\R$ mit Intensit"atma"s $e^{-x} dx$ sowie $0 < a_1 < a_2 < a_3 < \ldots$ "aquidistante Punkte mit Abstand $h > 0$ und $a_1 > 16+h, \ a_1 \in h\Z$, sowie $L > 0$, $L \in h\Z$, $C \in \R$. Dann gilt:
\begin{align*}
& \Prob\left( \exists i,j \in \N: C < X_i^{(j)} + \sup_{t \in [a_j, a_j + L]\cap h\Z} \left( W_i^{(j)}(t) - \frac t 2 \right) \leq X_i^{(j)} + W_i^{(j)}(0) \right) \\
\leq{} &  \frac {2e^{-C}} h \Bigg( \exp\left(- {\sqrt{a_1-h}} \right) + \frac {25L} {9} \sqrt{\frac{2L}{\pi}} \exp\left(-\frac{(\sqrt{a_1-h} - \log(\sqrt{a_1-h}))^2} {2L}\right)\\
  &  + \frac 8 {\sqrt{2\pi}} \left(1 + \frac 2 {\sqrt{a_1-h}-4}\right) \exp\left(- \frac{(\sqrt{a_1-h}-4)^2} 8\right) \Bigg).
\end{align*} 
\end{lemma}

\begin{bew}
Der Beweis verl"auft analog zum Beweis von \ref{BBsup}. Da sich jedoch einige Konstanten "andern und Integrationsgrenzen anders gew"ahlt werden, wird er hier dennoch in G"anze ausgef"uhrt. Seien also $\Pi_j, A_j, \Lambda_j$ f"ur $j \in \N$ genau wie dort definiert und $\Prob$ das zu $W(\cdot)$ geh"orende Wahrscheinlichkeitsma"s.

Wiederum wollen wir $\Lambda_j(A_j)$ nach oben hin absch"atzen. Es gilt:
\begin{align*}
&\Lambda_j(A_j) = \int_{\R} e^{-x} \Prob\left(C - x < \sup_{t \in [a_j, a_j+L]\cap h\Z} \left(W(t) - \frac t 2\right) \leq W(0)\right) dx \\
={} & \int_{\R} e^{-x} \Prob\left(C - x < \sup_{t \in [a_j, a_j+L] \cap h\Z} \left(W(t) - \frac t 2\right) \leq  0\right) dx,\\
& \hspace{7.5cm} \text{ da } W(0) = 0 \text{ $\Prob$-f.s.}\displaybreak[0]\\
\leq{} & \int_{C + {\sqrt{a_j}}  +  \log({\sqrt{a_j}})}^\infty e^{-x} dx \\
& + \int_C^{C + {\sqrt{a_j}} +  \log({\sqrt{a_j}})} e^{-x} \Prob\left(C - x < \sup_{t \in [a_j, a_j+L] \cap h\Z} \left(W(t) - \frac t 2\right) \leq  0\right) dx\displaybreak[0]\\
={} & e^{-C} \frac 1 {\sqrt{a_j}} \exp\left(-{\sqrt{a_j}} \right) \\
  & + e^{-C} \int_{- {\sqrt{a_j}} -  \log({\sqrt{a_j}} )}^0 e^{x} \Prob\left(x < \sup_{t \in [a_j, a_j+L] \cap h\Z} \left(W(t) - \frac t 2\right) \leq  0\right) dx\displaybreak[0]\\
\leq{} & e^{-C} \frac 1 {\sqrt{a_j}} \exp\left(-{\sqrt{a_j}}  \right) 
  + e^{-C} \int_{- {\sqrt{a_j}} -  \log({\sqrt{a_j}})}^0 e^{x} \Prob\Bigg(0 \geq W(a_j) - \frac {a_j} 2  > - 2{\sqrt{a_j}}, \\
  & \hspace{6.5cm} \sup_{t \in [a_j, a_j+L] \cap h\Z} \left(W(t) - \frac t 2\right) > x\Bigg) dx\\
  & + e^{-C} \int_{- {\sqrt{a_j}} -  \log({\sqrt{a_j}})}^0 e^{x} \Prob\Bigg(W(a_j) - \frac {a_j} 2 \leq -2{\sqrt{a_j}}, \\
  & \hspace{6.5cm} \sup_{t \in [a_j, a_j+L] \cap h\Z} \left(W(t) - \frac t 2\right) > x\Bigg) dx\displaybreak[0]\\
\leq{} & \frac 1 {\sqrt{a_j}} \exp\left(-{\sqrt{a_j}} \right) + e^{-C} \int_{- {\sqrt{a_j}} -  \log({\sqrt{a_j}} )}^0 e^{x} \Prob\left(0 \geq W(a_j) - \frac {a_j} 2 > -2{\sqrt{a_j}}\right) dx\\
  & + e^{-C} \int_{- {\sqrt{a_j}} -  \log({\sqrt{a_j}} )}^0 e^{x} \Prob\Bigg(W(a_j) - \frac {a_j} 2 \leq -2{\sqrt{a_j}},\\
  & \hspace{5.5cm} \sup_{t \in [0,L]} \left(\tilde{W}(t) - \frac t 2\right) > x-W(a_j) + \frac {a_j} 2\Bigg) dx,\\
  & \hspace{2.5cm} \text{wobei $\tilde{W}$ eine von $W$ unabh"angige Brownsche Bewegung sei.}
\end{align*}

Mit $W(a_j) - \frac {a_j} 2 \sim {\cal N}\left(-\frac {a_j} 2, a_j\right)$ erhalten wir

\begin{align}
&\Lambda_j(A_j) \leq{} e^{-C} \frac 1 {\sqrt{a_j}} \exp\left(-{\sqrt{a_j}} \right) \nonumber \\
  & + e^{-C} \int_{- {\sqrt{a_j}} -  \log( {\sqrt{a_j}} )}^0 e^{x} \int_{-2{\sqrt{a_j}}}^0 \frac 1 {\sqrt{2 \pi a_j}} \exp\left(-\frac {(y + \frac {a_j} 2)^2} {2a_j} \right) dy \ dx \nonumber \\
  & + e^{-C} \int_{- {\sqrt{a_j}} -  \log( {\sqrt{a_j}} )}^0 e^{x} \int_{-\infty}^{-2{\sqrt{a_j}}} \frac 1 {\sqrt{2 \pi a_j}} \exp\left(-\frac {(y + \frac {a_j} 2)^2} {2a_j} \right) \nonumber \\
  & \hspace{6cm} \Prob\left(\sup_{t \in [0,L]} \left(\tilde{W}(t) - \frac t 2\right) > x-y\right) dy \ dx \nonumber \displaybreak[0]\\
\leq{} & \frac {e^{-C}} {\sqrt{a_j}} \exp\left(-{\sqrt{a_j}} \right) + e^{-C} \int_{- {\sqrt{a_j}} -  \log( {\sqrt{a_j}} )}^0 e^{x} dx \cdot \int_{-2 + \frac {\sqrt{a_j}} 2}^{\frac {\sqrt{a_j}} 2} \frac 1 {\sqrt{2 \pi}} \exp\left(-\frac {y^2} {2} \right) dy  \nonumber \\
  & + e^{-C} \int_{- {\sqrt{a_j}} -  \log({\sqrt{a_j}})}^0 e^{x} \int_{-\infty}^{-2{\sqrt{a_j}}} \frac 1 {\sqrt{2 \pi a_j}} \exp\left(-\frac {(y + \frac {a_j} 2)^2} {2a_j} \right) \nonumber \\
  & \hspace{6cm} \Prob\left(\sup_{t \in [0,L]} W(t) > x-y\right) dy \ dx \nonumber \displaybreak[0]\\
\leq{} & e^{-C} \frac 1 {\sqrt{a_j}} \exp\left({\sqrt{a_j}} 1 \right) + e^{-C} \int_{-2 + \frac {\sqrt{a_j}} 2}^{\frac {\sqrt{a_j}} 2} \frac 1 {\sqrt{2 \pi}} \exp\left(-\frac {y^2} {2} \right) dy \nonumber \\
  & + e^{-C} \int_{-{\sqrt{a_j}} -  \log({\sqrt{a_j}})}^0 e^{x} \int_{-\infty}^{-2{\sqrt{a_j}}} \frac 1 {\sqrt{2 \pi a_j}} \exp\left(-\frac {(y + \frac {a_j} 2)^2} {2a_j} \right) \nonumber \\
  & \hspace{2.2cm} \cdot \frac 2 {\sqrt{\pi}} \int_{\frac{x-y}{\sqrt{2L}}}^\infty \exp(-z^2) dz \ dy \ dx \quad \text{nach Formel 1.1.1.4 aus \cite{borodin-2002}} \nonumber \displaybreak[0]\\
\leq{} & e^{-C} \frac 1 {\sqrt{a_j}} \exp\left(-{\sqrt{a_j}} \right) + e^{-C} \int_{-2 + \frac {\sqrt{a_j}} 2}^{\frac {\sqrt{a_j}} 2} \frac 1 {\sqrt{2 \pi}} \exp\left(-\frac {y^2} {2} \right) dy \nonumber \\
  & + e^{-C} \int_{- {\sqrt{a_j}} -  \log( {\sqrt{a_j}} )}^0 e^{x} \int_{-\infty}^{-2{\sqrt{a_j}}} \frac 1 {\sqrt{2 \pi a_j}} \exp\left(-\frac {(y + \frac {a_j} 2)^2} {2a_j} \right) \nonumber \\
  & \hspace{2cm} \cdot \frac {\sqrt{2L}} {\sqrt{\pi}(x-y)} \exp\left(-\frac {(x-y)^2}{2L}\right)  dy \ dx \label{eq:drittersum2}
\end{align}

F"ur den dritten Summanden aus \eqref{eq:drittersum2} l"asst sich weiterhin absch"atzen
\begin{align*}
& e^{-C} \int_{-{\sqrt{a_j}} -  \log({\sqrt{a_j}})}^0 e^{x} \int_{-\infty}^{-2{\sqrt{a_j}}} \frac 1 {\sqrt{2 \pi a_j}} \exp\left(-\frac {(y + \frac {a_j} 2)^2} {2a_j} \right) \\
  & \hspace{3cm} \cdot \frac {\sqrt{2L}} {\sqrt{\pi}(x-y)} \exp\left(-\frac {(x-y)^2}{2L}\right)  dy \ dx \\
\leq{} &  e^{-C} \int_{- {\sqrt{a_j}} -  \log({\sqrt{a_j}})}^0 e^{x} \int_{-\infty}^{-2\sqrt{a_j}} \frac 1 {\sqrt{a_j}} \exp\left(-\frac {(y + \frac {a_j} 2)^2} {2a_j} \right) \\
  & \hspace{3cm} \cdot \frac 1 \pi{}\frac {\sqrt{L}} {{\sqrt{a_j}} -  \log( {\sqrt{a_j}} )} \exp\left(-\frac {(x-y)^2}{2L}\right)  dy \ dx
\end{align*}
\begin{align*}
\leq{} &  e^{-C} \int_{- {\sqrt{a_j}} -  \log({\sqrt{a_j}} )}^0 \int_{-\infty}^{-2{\sqrt{a_j}}} \frac{e^{x}}{\pi}\sqrt{\frac L {a_j}} \exp\left(-\frac {(y + \frac {a_j} 2)^2} {2a_j} \right) \\
& \hspace{7cm} \cdot \exp\left(-\frac {(x-y)^2}{2L}\right)  dy \ dx,\\
& \hspace{1cm} \text{da f"ur alle \ $x > 1$ \ die Ungleichung \ $x - \log(x) > 1 - \log(1) = 1$ \ gilt} \\
\leq{} &  e^{-C} \int_{- {\sqrt{a_j}} -  \log({\sqrt{a_j}} )}^0 e^{x} dx \int_{-\infty}^{-2{\sqrt{a_j}}} \frac 1 {\sqrt{a_j}} \exp\left(-\frac {(y + \frac {a_j} 2)^2} {2a_j} \right) dy \\
& \hspace{6cm} \cdot \frac{\sqrt{L}}{\pi} \exp\left(-\frac {(\sqrt{a_j}  -  \log(\sqrt{a_j}))^2}{2L}\right).
\end{align*}

Mit der gleichen Argumentation wie im Beweis von Lemma \ref{BBsup} ist die Funktion $x \mapsto \frac 2 5 x - \log(x)$ positiv und isoton auf $\R_{\geq \frac 5 2}$ und somit insbesondere auf $$[\sqrt{a_1-h}, \infty) \subset [4, \infty)$$ isoton.
\medskip

Daraus ergibt sich
\begin{align*}
& \left( {\sqrt{a_j}} - \log\left({\sqrt{a_j}} \right)\right)^2 = \left(\frac {3\sqrt{a_j}} {5} + \left(\frac {2\sqrt{a_j}} {5} - \log\left({\sqrt{a_j}} \right)\right)\right)^2\\
={} & \frac 9 {25} a_j +  \frac 6 5 \sqrt{a_j} \left(\frac {2\sqrt{a_j}} 5 - \log({\sqrt{a_j}})\right) + \left(\frac {2\sqrt{a_j}} 5 - \log\left({\sqrt{a}} \right)\right)^2\\
\geq{} & \frac 9 {25} a_j  +  \underbrace{\frac 6 5 \sqrt{a_1-h} \left(\frac {2\sqrt{(a_1-h)}} 5  - \log\left({\sqrt{(a_1-h)}} \right)\right)}_{K_{a_1}^{(1)}} \\
& + \underbrace{\left(\frac {2\sqrt{(a_1-h)}} 5 - \log\left({\sqrt{(a_1-h)}}\right)\right)^2}_{=:K_{a_1}^{(2)}}
\end{align*}
f"ur alle $j \in \N $ und mit $K_{a_1} := K_{a_1}^{(1)} +K_{a_1}^{(2)}$ ergibt sich die Ungleichung
\begin{align*}
&  e^{-C} \int_{- \frac {\sqrt{a_j}} 2 -  \log(\frac {\sqrt{a_j}} 2)}^0 e^{x} dx \int_{-\infty}^{-{\sqrt{a}}} \frac 1 {\sqrt{a_j}} \exp\left(-\frac {(y + \frac {a_j} 2)^2} {2a_j} \right) dy \\
& \hspace{6cm} \cdot \frac{\sqrt{L}}{\pi} \exp\left(-\frac {({\sqrt{a_j}} -  \log({\sqrt{a_j}}))^2}{2L}\right)\\
\leq{} &  e^{-C} e^0 \cdot  \sqrt{2\pi} \cdot \frac{\sqrt{L}}{\pi} \exp\left(-\frac {\frac 9 {25} a_j} {2L}\right) \exp\left(\frac {K_{a_1}} {2L}\right).
\end{align*}

Einsetzen in \ref{eq:drittersum2} f"uhrt schlie"slich zu

\begin{align}
\Lambda_j(A_j) \leq{} & e^{-C} \frac 1 {\sqrt{a_j}} \exp\left(-{\sqrt{a_j}} \right) + e^{-C} \int_{-2 + \frac {\sqrt{a_j}} 2}^{\frac {\sqrt{a_j}} 2} \frac 1 {\sqrt{2 \pi}} \exp\left(-\frac {y^2} {2} \right) dy \nonumber \\
  & + e^{-C} \sqrt{\frac{2L}{\pi}} \exp\left(-\frac {\frac 9 {25} a_j} {2L}\right) \exp\left(\frac {K_{a_1}} {2L}\right)
   \label{eq:mass2}.
\end{align}

Nun ist
\begin{align*}
& \Prob\left( \exists i,j \in \N: C < X_i^{(j)} + \sup_{t \in [a_j, a_j + L]} \left( W_i^{(j)}(t) - \frac t 2 \right) \leq X_i^{(j)} + W_i^{(j)}(0) \right) \displaybreak[0]\\
={} & \Prob\left( \exists j \in \N: \Pi_j(A_j) > 0 \right) ={}  1 - \exp\left(- \sum_{j \in \N} \Lambda_j(A_j)\right) \displaybreak[0]\\
\leq{} & \sum_{j \in \N} \Lambda_j(A_j), \quad \text{da } 1 - e^{-x} > x \text{ f"ur alle } x > 0 \text{ gilt} \displaybreak[0]\\
\leq{} & \sum_{j \in \N} \Bigg(e^{-C} \frac 1 {\sqrt{a_j}} \exp\left(-{\sqrt{a_j}}\right) + e^{-C} \int_{-2 + \frac {\sqrt{a_j}} 2}^{\frac {\sqrt{a_j}} 2} \frac 1 {\sqrt{2 \pi}} \exp\left(-\frac {y^2} {2} \right) dy \nonumber \\
  & \hspace{1cm}+\sqrt{\frac {2L}{\pi}} \exp\left(-\frac {\frac 9 {25} a_j} {2L}\right)\exp\left(\frac{K_{a_1}}{2L}\right) \Bigg) \qquad \text{nach \eqref{eq:mass2}}\displaybreak[0]\\
\leq{} & \int_{a_1-h}^\infty \frac 1 h \Bigg(e^{-C} \frac 1 {\sqrt{a}} \exp\left(-{\sqrt{a}}\right) + e^{-C} \int_{-2 + \frac {\sqrt{a}} 2}^{\frac {\sqrt{a}} 2} \frac 1 {\sqrt{2 \pi}} \exp\left(-\frac {y^2} {2} \right) dy \nonumber \\
  & \hspace{2cm}+ \sqrt{\frac{2L}{\pi}} \exp\left(-\frac {\frac 9 {25} a} {2L}\right) \exp\left(\frac{K_{a_1}}{2L}\right) \Bigg) da \displaybreak[0]\\
\leq{} & e^{-C} \frac 2 h \exp\left(-{\sqrt{a_1-h}} \right)\\
  & + e^{-C} \frac 1 h \int_{-2 + \frac {\sqrt{a_1-h}} 2}^\infty \int_{a_1-h}^\infty \mathbf{1}_{y \in [-2 + \frac {\sqrt{a}} 2, \frac {\sqrt{a}} 2]} \frac 1 {\sqrt{2 \pi}} \exp\left(-\frac {y^2} {2} \right) da \ dy \\
  & +  e^{-C} \frac {50L} {9h} \sqrt{\frac{2L}{\pi}} \exp\left(-\frac {\frac 9 {25} (a_1 - h)} {2L}\right) \exp\left(\frac{K_{a_1}}{2L}\right)\displaybreak[0]\\
\leq{} &  \frac {e^{-C}} h \Bigg(2 \exp\left(-{\sqrt{a_1-h}} \right) +  \int_{- 2 + \frac {\sqrt{a_1-h}} 2}^\infty (16y+16) \frac 1 {\sqrt{2 \pi}} \exp\left(-\frac {y^2} {2} \right) dy \\
  & + \frac {50L} {9} \sqrt{\frac{2L}{\pi}} \exp\left(-\frac{(\sqrt{a_1-h} - \log(\sqrt{a_1-h}))^2} {2L}\right),\displaybreak[0]\\
 &\hspace{4cm} \text{da } y \in \left[- 2 + \frac {\sqrt{a}} 2, \frac {\sqrt{a}} 2\right] \iff a \in [4y^2,4(y+2)^2] \displaybreak[0]\displaybreak[0]\\
\leq{} & \frac {2e^{-C}} h \Bigg( \exp\left(- {\sqrt{a_1-h}} \right) + \frac {25L} {9} \sqrt{\frac{2L}{\pi}} \exp\left(-\frac{(\sqrt{a_1-h} - \log(\sqrt{a_1-h}))^2} {2L}\right)\displaybreak[0]\\
  &  + \frac 8 {\sqrt{2\pi}} \left(1 + \frac 2 {\sqrt{a_1-h}-4}\right) \exp\left(-\frac{(\sqrt{a_1-h}-4)^2} 8\right) \Bigg)
\end{align*}
\end{bew}
\medskip

\subsection{Approximation von $Z(\cdot)$}

Zun"achst wollen wir einen Brown-Resnick-Prozess gem"a"s der Definition \ref{brown-resnick} approximieren. Dazu erzeugen wir zun"achst
 unabh"angig identisch standardexponentialverteilte Zufallsvariablen $Y_i,\ i \in \N,$ und setzen 
$X_k := \log\left(\frac{1}{\sum_{i=1}^k Y_i}\right)$. Dann ist $\sum_{i \in \N} \delta_{X_i}$ ein Poisson"=Punktprozess
 mit Intensit"at $e^{-x} dx$ nach Korollar \ref{PPP-e-x}.
\medskip

Weiterhin sei zu jedem $X_i$ ein (stochastisch unabh"angiger) Pfad $\{W_i(t) - \frac {\sigma^2(t)} 2: t \in \R^d \}$ gegeben. Definieren wir
$Z^{(k)}(t) := \max_{i=1,\ldots,k} X_i + W_i(t) - \frac {\sigma^2(t)} 2, \ k \in \N, \ t \in \R^d$, so gilt nach Konstruktion $Z^{(k)}(t) \stackrel{k \rightarrow \infty}{\longrightarrow} Z(t)$ f"ur alle $t \in \R^d$.
\medskip

Wir wollen nun den Fehler der Approximation nach $k$ Schritten berechnen. Dazu betrachten wir den Spezialfall, dass $W(\cdot)$ die eindimensionale Standard"=Brownsche"=Bewegung ist. Zudem schr"anken wir die Approximation
 auf ein Intervall $[-b,b] \subset \R, \ b> 0$ ein.
Sei nun $C:= \inf_{t \in [a,b]} ( Z^{(k)}(t) + \frac {\sigma^2(t)} 2 )$.
\medskip

\begin{prop} \label{fehler_0}
Seien $x(k) < c(k) < 0$. Dann gilt
\begin{align} 
& \Prob(Z^{(k)}(t) \neq Z(t) \text{ f\"ur ein } t \in [-b,b] \ |\ X_k \leq x(k), C_k > c(k)) \displaybreak[0] \nonumber\\
\leq{} & 4 e^{-c(k)} \frac b {c(k)-x(k)} \exp\left(\frac b 2\right)
\left(1- \Phi\left( \frac{C-x-b}{\sqrt{b}}\right) \right). \label{eq:absch01}
\end{align}
Weiterhin
\begin{align}
 \Prob(C_k \leq c(k)) \leq 1 - \left(1 - \exp\left(-\exp\left(-\frac{c(k)} 2\right)\right)\right) \cdot \left(2\Phi\left(-\frac{c(k)} {2\sqrt{b}}\right) - 1\right)^2 \label{eq:absch02}
\end{align}
und
\begin{align}
 \Prob(X_k > x(k)) \leq \frac{\exp(-(k-1) x(k))}{(k-1)!} (1 - \exp(-\exp(-x(k)))). \label{eq:absch03}
\end{align}
Die unbedingte Fehlerwahrscheinlichkeit $\Prob(Z^{(k)}(t) \neq Z(t) \text{ f\"ur ein } t \in [-b,b])$ kann durch diese
Summe der rechten Seiten von \eqref{eq:absch01}, \eqref{eq:absch02} und \eqref{eq:absch03} beschr\"ankt werden.
\end{prop}

\begin{bew}
Es gilt:
\begin{align*}
& Z^{(k)}(t) \neq Z(t) \text{ f"ur ein } t \in [-b,b] \ |\ X_k \leq x(k), C_k > c(k) \\
\iff {}& Z^{(k)}(t) < Z(t) \text{ f"ur ein } t \in [-b,b] \ |\ X_k \leq x(k), C_k > c(k) \displaybreak[0]\\  
\Longrightarrow {}& \exists m>k: \ X_m + W_m(t) > Z^{(k)}(t) + \frac {\sigma^2(t)} 2 \geq c(k) \text{ f"ur ein } t \in [-b,b] \\
  & \hspace{8cm} |\ X_k \leq x(k), C_k > c(k)\\
\Longrightarrow {}& \exists m>k: \ \sup_{t \in [-b,b]} (X_m + W_m(t)) > c(k) \ | \ X_k \leq x(k).
\end{align*}

Weiterhin ist die Folge $(X_m)_{m \in \N}$ nach Konstruktion monoton fallend, d.\,h. aus $m > k$ folgt $X_m \leq X_k \leq x(k)$.

Daraus erhalten wir folgende Fehlerabsch"atzung:

\begin{align*}
& \Prob(Z^{(k)}(t) \neq Z(t) \text{ f"ur ein } t \in [-b,b]\ |\ X_k \leq x(k), C_k > c(k))\\
\leq{} & \Prob(\exists i\in \N: X_i \leq x(k), \ \sup_{t \in [-b,-b]} (X_i + W_i(t)) > c(k))\\
\leq{} & 4 e^{-c(k)} \frac b {c(k)-x(k)} \exp\left(\frac {b} 2\right)
           \left(1- \Phi\left( \frac{C-x-b}{\sqrt{b}}\right) \right),
\end{align*}
wobei der letzte Schritt aus Lemma \ref{fehlerlemma} folgt.

Weiterhin gilt
 \begin{align*}
  \Prob(C_k > c(k)) \geq{}& \Prob(X_1 > c(k)/2) \Prob\left(\inf_{t \in [-b,b]} W(t) > c(k)/2\right)\\
  ={} & (1 - \exp(-\exp(-c(k)/2))) \cdot \Prob\left(\inf_{t \in [0,b]} W(t) > c(k)/2\right)^2 \\
  ={} & (1 - \exp(-\exp(-c(k)/2))) \cdot \left(2\Phi\left(-\frac{c(k)} {2\sqrt{b}}\right) - 1\right)^2
 \end{align*}
  nach Formel 1.1.1.4 in \cite{borodin-2002},
und
\begin{align*}
  \Prob(X_k > x(k)) ={} & \Prob(e^{-X_k} < e^{-x(k)})\\
                    ={} & \int_0^{e^{-x(k)}} t^{k-1} \frac{\exp(-t)}{(k-1)!} d t\\
                    \leq{} & \frac{\exp(-(k-1) x(k))}{(k-1)!} \int_0^{e^{-x(k)}} \exp(-t) d t\\
                    ={} & \frac{\exp(-(k-1) x(k))}{(k-1)!} (1 - \exp(-\exp(-x(k)))).
\end{align*}
Dann k\"onnen wir die Ungleichung 
\begin{align*}
&\Prob(Z^{(k)}(t) \neq Z(t) \text{ f\"ur ein } t \in [-b,b])\\
={} & \Prob(Z^{(k)}(t) \neq Z(t) \text{ f\"ur ein } t \in [-b,b] \ |\ X_k \leq x(k), C_k > c(k))\\
    & \quad \cdot \Prob(X_k \leq x(k), C_k > c(k)) \\
& +{} \Prob(Z^{(k)}(t) \neq Z(t) \text{ f\"ur ein } t \in [-b,b] \ |\ X\{X_k > x(k)\} \cup \{C_k \leq c(k)\})\\
  & \quad \cdot \Prob(\{X_k > x(k)\} \cup \{C_k \leq c(k)\}) \displaybreak[0] \\
\leq{} & \Prob(Z^{(k)}(t) \neq Z(t) \text{ f\"ur ein } t \in [-b,b] \ |\ X_k \leq x(k), C_k > c(k))\\
& {}+{} \Prob(\{X_k > x(k)\} \cup \{C_k \leq c(k)\})\\
\leq{} & \Prob(Z^{(k)}(t) \neq Z(t) \text{ f\"ur ein } t \in [-b,b] \ |\ X_k \leq x(k), C_k > c(k))\\
& {}+{} \Prob(C_k \leq c(k)) {}+{} \Prob(X_k > x(k)).
\end{align*}
\end{bew}
\medskip

Um die Konvergenz $Z^{(k)}$ gegen $Z$ im Sinne von $$ \lim_{k \to \infty} \Prob(Z^{(k)}(t) \neq Z(t) \text{ f\"ur ein } t \in [-b,b]) = 0$$
einzusehen, m\"ussen $c(k)$ und $x(k)$ so gew\"ahlt werden, dass jede der rechten seiten von \eqref{eq:absch01}, \eqref{eq:absch02}
und \eqref{eq:absch03} f\"ur $k \to \infty$ gegen Null geht. Dies ist z.\,B. f\"ur $c(k) = - \log\log(\log(k)/2)$
und $x(k)  = - \log(k)/2$ erf\"ullt.
\bigskip

\subsection{Approximation von $Z_1(\cdot)$}
Nun soll ein Brown-Resnick-Prozess gem"a"s Korollar \ref{z1} approximiert werden. Wir erzeugen dazu u.\,i.\,v. Exp($\frac 1 n$)-Zufallsvariablen $Y_i^{(j)}, \ i \in \N,\ j \in \{1,\ldots,n\},$ und definieren \[X_k^{(j)}:= \log\left(\frac{1}{\sum_{i=1}^k Y_i^{(j)}}\right).\]
Dann ist f"ur jedes $j \in \{1,\ldots,n\}$ nach Korollar \ref{PPP-e-x} der Punktprozess $\sum_{i \in \N} \delta_{X_i^{(j)}}$ ein Poisson"=Punktprozess auf $\R$ mit Intensit"atsma"s $\frac 1 n e^{-x}$.\\
F"ur jedes Paar $(i,j) \in \N \times \{1,\ldots,n\}$ sei $W_i^{(j)}(\cdot)$ unabh"angig identisch verteilt wie $W(\cdot)$. 

Seien weiterhin $ k \in \N$ sowie -- wie in Korollar \ref{z1} -- $h_1, \ldots, h_n \in \R^d$ gegeben.
Dann definieren wir
\[Z_1^{(k)}(t) := \max_{j=1,\ldots,n} \max_{i=1,\ldots,k} X_i^{(j)} + W_i^{(j)}(t-h_j) - \frac {\sigma^2(t-h_j)} 2. \]
Wiederum gilt $\lim_{k \to \infty} Z_1^{(k_1,\ldots,k_n)}(t) = Z_1(t)$ f"ur jedes $t \in \R^d$.

F"ur die Fehlerabsch"atzung beschr"anken wir uns erneut auf den Fall, dass $W(\cdot)$ eine eindimensionale Standard-Brownsche-Bewegung
 ist und betrachten die Simulation auf einem Intervall $[-b,b] \subset \R$, $b>0$. Weiterhin sei
 $h_1 < h_2 < \ldots < h_n$  und die Menge $\{h_j, \ 1 \leq j \leq n\}$ sei symmetrisch
Seien $X_k:= \max_{j=1,\ldots,n} X_{k}^{(j)}$, $C_k:= \inf_{t \in [a, b]} \left(Z_1^{(k)}(t) + \min_{j=1,\ldots,n} \frac {\sigma^2(t-h_j)} 2 \right)$.

\begin{prop} \label{fehler_1}
Sei $ x(k) < c(k) < 0$.
Dann gilt
\begin{align*}
& \Prob(Z_1^{(k)}(t) \neq Z_1(t) \text{ f\"ur ein } t \in [-b,b] \ |\ X_k \leq x(k), C_k > c(k)) \displaybreak[0] \nonumber\\
\leq{} & 4 e^{-c(k)} \frac {b-h_1} {c(k)-x(k)} e^{(b-h_1)/2}
\left(1- \Phi\left( \frac{C-x-(b-h_1)}{\sqrt{b-h_1}}\right) \right).
\end{align*}
Weiterhin ist
 \begin{align*}
 & \Prob(C_k \leq c(k))\\
 \leq{} & 1 - \left(1-\exp(-\exp(-c(k)/2))\right) \\
 & \hspace{0.2cm} \cdot\left(
 \Phi\left(-\frac{c(k)+ (b-h_1)}{2\sqrt{b-h_1}}\right) + e^{-\frac{c(k)}2}\left(1-\Phi\left(-\frac{c(k)-(b-h_1)}{2\sqrt{b-h_1}}\right)\right)\right)^2
 \end{align*}
 und
 $$ \Prob(X_k > x(k)) \leq{} n^{k+1} \frac{\exp(-(k-1) x(k))}{(k-1)!} (1 - \exp(-\exp(-nx(k)))).$$
\end{prop}
\begin{bew}
Es gilt
\begin{align*}
& Z_1^{(k)}(t) \neq Z_1(t) \text{ f"ur ein } t \in [-b,b] \ |\ X_k \leq x(k), C_k > c(k) \\
\iff& Z_1^{(k)}(t) < Z_1(t) \text{ f"ur ein } t \in [-b,b] \ |\ X_k \leq x(k), C_k > c(k) \\
\Longrightarrow{}& \exists j \in \{1,\ldots,n\},\ m>k: \ X_m^{(j)} + W_m^{(j)}(t-h_j) > Z_1^{(k)}(t) + \frac {\sigma^2(t-h_j)} 2 \geq c(k) \\
                 & \hspace{4.5cm} \text{f"ur ein } t \in [-b,b] \ |\ X_k \leq x(k), C_k > c(k)  \\
\Longrightarrow{}& \exists j \in \{1,\ldots,n\},\ m>k: \ \sup_{t \in [-b-h_n,b-h_1]} (X_m^{(j)} + W_m^{(j)}(t)) > c(k)
\ |\ X_k \leq x(k) .
\end{align*}

Nun ist jede Folge $(X_n^{(j)})_{n \in \N}$ monoton fallend und daher folgt aus $m > k$, dass $X_m^{(j)} \leq X_{k}^{(j)} \leq x(k)$.
Wir erhalten somit unter Verwendung von $-b-h_n = - (b-h_1)$, dass
\begin{align*}
& \Prob(Z_1^{(k_1,\ldots,k_n)}(t) \neq Z_1(t) \text{ f"ur ein } t \in [-b,b] \ |\ X_k \leq x(k), C_k > c(k) )\\
\leq{} &\Prob(\exists j \in \{1,\ldots,n\},\ m>k: \ \sup_{t \in [-b-h_n,b- h_1]} (X_m^{(j)} + W_m^{(j)}(t)) > c(k) \ |\ X_k \leq x(k))\\
\leq {} & \sum_{j=1}^n \Prob(\exists i\in \N: \ X_i^{(j)} \leq x, \ \sup_{t \in [-b-h_n,b- h_1]} (X_i^{(j)} + W_i^{(j)}(t)) > c(k))
\end{align*}
\begin{align*}
\leq {} & 2 e^{-c(k)} \sqrt{ \frac 2 \pi} \frac {\sqrt{(b-h_1)^3}} {c(k)-x(k)} e^{(b-h_1)/2}
\left(1 + \frac{1}{c(k)-x(k)} \exp\left(-\frac{(c(k)-x(k))^2}{2(b-h_1)}\right)\right),
\end{align*}
wobei im letzten Schritt Lemma \ref{fehlerlemma} verwendet wurde.

F\"ur $c(k) < 0$ gilt zudem mit $X_1 = \max_{j=1,\ldots,n} X_i^{(j)}$, dass
\begin{align*}
 \Prob(C_k > c(k)) \geq{}& \Prob(X_1 > c(k)/2) \Prob\left(\inf_{t \in [-(b-h_1),b-h_1]} (W(t)-|t|/2) > c(k)/2\right)\\
 ={}& \Prob(X_1 > c(k)/2) \Prob\left(\inf_{t \in [0,b-h_1]} (W(t)-|t|/2) > c(k)/2\right)^2\\
 ={}&  \left(1-\exp(-\exp(-c(k)/2))\right) \cdot \Bigg(\Phi\left(-\frac{c(k)}{2\sqrt{b-h_1}} - \frac{\sqrt{b-h_1}} 2\right)\\
& \hspace{0.7cm} + \exp(-c(k)/2) \left(1-\Phi\left(-\frac{c(k)}{2\sqrt{b-h_1}} + \frac{\sqrt{b-h_1}} 2\right)\right) \Bigg)^2
\end{align*}
nach Formel 2.1.1.4 in \cite{borodin-2002}.
Weiterhin ist
\begin{align*}
  \Prob(X_k^{(j)} > x(k)) ={} & \Prob(e^{-X_k^{(j)}} < e^{-x(k)})\\
                    ={} & \int_0^{e^{-x(k)}} n^{k} t^{k-1} \frac{\exp(-nt)}{(k-1)!} d t\displaybreak[0]\\
                    \leq{} & n^{k} \frac{\exp(-(k-1) x(k))}{(k-1)!}
                                    \int_0^{e^{-x(k)}} \exp(-nt) d t\\
                    ={} & n^{k}\frac{\exp(-(k-1) x(k))}{(k-1)!} (1 - \exp(-\exp(-nx(k))))
\end{align*}
f\"ur $j=1,\ldots,n$.
Die Behauptung folgt mit $\Prob(X_k > x(k)) \leq \sum_{j=1}^n \Prob(X_k^{(j)} > x(k))$.
\end{bew}
\medskip

Konvergenz erh\"alt man beispielsweise wiederum f\"ur $c(k) = - \log\log(\log(k)/2)$
und $x(k)  = - \log(k)/2$.
\bigskip

\subsection{Approximation von $Z_2(\cdot)$}
Es soll ein Brown-Resnick-Prozess gem"a"s Korollar \ref{z2} approximiert werden.
 Wir erzeugen dazu unabh"angig identisch standardexponentialverteilte Zufallsvariablen $Y_i,\ i \in \N,$
 und definieren $X_k:= \log\left(\frac 1 {\sum_{i=1}^k Y_i}\right)$ sowie unabh"angig identisch auf $I$
 gleichverteilte Zufallsvariablen $H_i,\ i \in \N$. Dann ist nach Korollar \ref{PPP-dt-e-x}
 durch $\sum_{i \in \N} \delta_{(X_i, H_i)}$ ein Poisson"=Punktprozess auf $\R\times I$
 mit Intensit"atsma"s $\frac 1 {|I|} e^{-x} dx \ dh$ gegeben. Weiterhin seien $W_1(\cdot), W_2(\cdot), \ldots$
 unabh"angig identisch verteilt wie $W(\cdot)$.

F"ur $k \in \N$ definieren wir
\[Z_2^{(k)}(t) := \max_{i=1,\ldots,k} \left( X_i + W_i(t-H_i) - \frac {\sigma^2(t-H_i)} 2 \right). \]
Dann gilt $\lim_{k \to \infty} Z_2^{(k)}(t) = Z_2(t)$ f"ur jedes $t \in \R^d$.

F"ur die Fehlerabsch"atzung beschr"anken wir uns erneut auf den Fall, dass $W(\cdot)$ eine eindimensionale
 Standard-Brownsche-Bewegung ist und betrachten die Simulation nur auf einem Intervall $[-b,b] \subset \R$, $b>0$.

Sei $C_k= \inf_{t \in [a, b]} (Z_2^{(k)}(t) + \min_{h \in I} \frac {\sigma^2(t-h)} 2 )$. Ferner nehmen wir an, dass
$I$ symmetrisch sei. Schlie"slich definieren wir $t_o := b - \min_{h \in I} h$.

\begin{prop} \label{fehler_2}
Unter den genannten Voraussetzungen und Notationen gilt
\begin{align*}
& \Prob(Z_2^{(k)}(t) \neq Z(t) \text{ f"ur ein } t \in [-b,b] \ | \ X_k \leq x(k), C_k > c(k)) \\
\leq {} & 4 e^{-c(k)}\frac {t_o} {c(k)-x(k)} \exp\left(\frac {t_o} 2\right) \left(1- \Phi\left( \frac{C-x-t_o}{\sqrt{t_o}}\right) \right).
\end{align*}
Weiterhin ist
 \begin{align*}
 & \Prob(C_k \leq c(k))\\
 \leq{} & 1 - \left(1-\exp(-\exp(-c(k)/2))\right) \\
 & \hspace{0.2cm} \cdot\left(
 \Phi\left(-\frac{c(k)+ t_o}{2\sqrt{b-h_1}}\right) + e^{-\frac{c(k)}2}\left(1-\Phi\left(-\frac{c(k)-t_o}{2\sqrt{b-t_o}}\right)\right)\right)^2
 \end{align*}
 und
 $$ \Prob(X_k > x(k)) \leq{} \frac{\exp(-(k-1) x(k))}{(k-1)!} (1 - \exp(-\exp(-x(k)))).$$
\end{prop}
\begin{bew}
Es ist
\begin{align*}
& Z_2^{(k)}(t) \neq Z_2(t) \text{ f"ur ein } t \in [a,b]  \ | \ X_k \leq x(k), C_k > c(k)\\
\iff {}& Z_2^{(k)}(t) < Z_2(t) \text{ f"ur ein } t \in [a,b]  \ | \ X_k \leq x(k), C_k > c(k) \\  
\Longrightarrow {}& \exists m>k: \ X_m + W_m(t-H_m) > Z_2^{(k)}(t) + \frac {\sigma^2(t-H_m)} 2 \geq c(k)\\
& \hspace{5cm} \text{ f"ur ein } t \in [a,b]  \ | \ X_k \leq x(k), C_k > c(k)\\
\Longrightarrow {}& \exists m>k: \ \sup_{t \in [-t_o,t_o]} (X_m + W_m(t)) > c(k)  \ | \ X_k \leq x(k) \\
\Longrightarrow {}& \exists i \in \N: \ X_i \leq x(k), \ \sup_{t \in [-t_o,t_o]} (X_i + W_i(t)) > c(k).
\end{align*}
Die Aussage der Proposition folgt nun direkt aus Lemma \ref{fehlerlemma}.

Die Absch\"atzung f\"ur $\Prob(C_k \leq c(k))$ l\"asst sich wie in Proposition \ref{fehler_1} beweisen,
die Absch\"atzung f\"ur $\Prob(X_k > x(k))$ wie in Propostion \ref{fehler_0}.
\end{bew}
\bigskip

\subsection{Approximation von $Z_3(\cdot)$}

Wir wollen einen Brown-Resnick-Prozess nach Satz \ref{z3} approximieren. Wir erzeugen zun"achst u.\,i.\,v.
 Exp($\frac 1 {m^d}$)-Zufallsvariablen $Y_i^{(j)}, \ i \in \N, \ j \in \Z^d,$ und definieren 
\[X_k^{(j)}:= \log\left(\frac{1}{\sum_{i=1}^k Y_i^{(j)}}\right).\]
Dann ist f"ur jedes $j \in \Z^d$ nach Korollar \ref{PPP-e-x} der Punktprozess $\sum_{i \in \N} \delta_{X_i^{(j)}}$
 ein Poisson"=Punktprozess auf $\R$ mit Intensit"atsma"s $\frac 1 {m^d} e^{-x}$.
F"ur jedes Paar $(i,j) \in \N \times \Z^d$ sei $W_i^{(j)}(\cdot)$ unabh"angig identisch verteilt wie $W(\cdot)$.

Seien weiterhin $j_{min},\ j_{max} \in \Z^d$ mit $j_{min} < j_{max}$ und $k \in \N$ gegeben.
Dann definieren wir
\[Z_3^{(k)}(t) := 
\max_{j=j_{min},\ldots,j_{max}} \max_{\substack{i=1,\ldots,k\\ T_i^{(j)} \in \left(-\frac{m p}{2}, \frac{m p}{2} \right] }}
 X_i^{(j)} + W_i^{(j)}(t-pj) - \frac {\sigma^2(t-pj)} 2. \]

Nach Konstruktion gilt \[\lim_{\substack{j_{min} \to -\infty, j_{max} \to \infty,\\ k \to \infty}} Z_3^{(k)}(t) = Z_3(t)\] f"ur jedes $t \in \R^d$.

Wir wollen nun auch hier eine Fehlerabsch"atzung f"ur die Approximation durch"-f"uhren.
Dabei beschr"anken wir uns wiederum auf den Fall, dass $W(\cdot)$ eine eindimensionale Standard-Brownsche-Bewegung ist
 und betrachten die Simulation auf einem Intervall $[-b,b] \subset \R, \ b \in p\Z_{>0}$.
Weiterhin sei $m=1$; $j_{min} = -j_{max} \in \Z_{<0}$ seien so gew"ahlt, dass $b+4 < p\cdot j_{max}$.

Zur Vereinfachung betrachten wir den Prozess lediglich auf dem Gitter $p\Z$.
Wir definieren $$T_i^{(j,p)} := \inf\left(\argsup\left( \left( W_i^{(j)}(t) - \frac {|t|} 2 \right)
\Big|_{t \in p\Z} \right) \right)$$
und approximieren $Z_3(\cdot)$ durch 
$$ Z_3^{(k; p)} (t) := \max_{j=j_{min},\ldots,j_{max}} \max_{\substack{i=1,\ldots,k\\ T_i^{(j,p)} \in \left(-\frac{m p}{2},
 \frac{m p}{2} \right] }} X_i^{(j)} + W_i^{(j)}(t-pj) - \frac {\sigma^2(t-pj)} 2$$
f"ur $t \in p\Z$.
\medskip

Seien nun $X_k:= \max_{j=j_{min},\ldots,j_{max}} X_{k}^{(j)}$, $C_k:= \inf_{t \in [-b, b]\cap p\Z}
 Z_3^{(k; p)}(t)$.

Unter der Annahme $X_k < C_k$ gilt f"ur jedes $j \in \{ j_{min}, \ldots, j_{max} \}$, $i > k$
mit $T_i^{(j, p)} \in \left( - \frac p 2, \frac p 2 \right]$, dass $$\max_{t\in [a,b] \cap p\Z} X_i^{(j)} + W(t) - \frac t 2
 = X_i^{(j)} < X_k < C_k.$$

Somit k"onnen alle Pfade $X_i^{(j)} + W_i^{(j)}(\cdot - pj) - \frac {\sigma^2(\cdot-pj)} 2$ f"ur $i > k$
 ohne Verschlechterung der Approximation abgebrochen werden, sobald $X_{k}^{(j)} < C_k$ ist.
\medskip

\begin{prop} \label{fehler_3}
Sei $c(k)<0$. Dann gilt:
\begin{align*}
& \Prob(Z_3^{(k; p)}(t) \neq \tilde{Z}_3^{(p)}(t) \text{ f"ur ein } t \in [-b,b]\cap p\Z \ | X_k \leq c(k), C_k > c(k)) \\
\leq {} &  \frac {16e^{-c(k)}} p \Bigg( \exp\left(-\frac {\sqrt{pj_{max}-b}} 2 \right) + \frac {100} {9} \sqrt{\frac{b^3}{\pi}}\\
& \hspace{0.5cm} \cdot \exp\left(-\frac{(\frac 3 {10} \sqrt{pj_{max}-b} + \frac 1 5 \sqrt{(pj_{max}-b)\vee25} -\log(\frac{\sqrt{(pj_{max}-b)\vee 25}} 2))^2} {4b}\right)\\
 & \hspace{1cm} + \frac 1 {\sqrt{2\pi}} \left(1 + \frac 1 {\sqrt{pj_{max}-o}-2}\right) \exp\left(- \frac{(\sqrt{pj_{max}-o}-2)^2} 8\right)\Bigg)
\end{align*}
Weiterhin gelten
\begin{align*}
 & \Prob(C_k \leq c(k)) \\
  \leq{} & 1 - \left(1-\exp\left(-\frac{2b/p+1}{4} \cdot \exp\left(-\frac{c(k)}2\right) \left(1-\exp\left(-\frac p 2\right)\right)^2\right)\right)\\
 & \cdot \Bigg( 1- \frac{4}{1-\exp(-p/2)}\Bigg( 1 - \Phi\left(-\frac{c(k)}{\sqrt{8b}} - \sqrt{\frac b 2}\right)\\
   & \hspace{2.5cm} - \exp(-c(k)/2)\left(1-\Phi\left(-\frac{c(k)}{\sqrt{8b}} + \sqrt{\frac b 2}\right)\right) \Bigg)^2\Bigg)
 \end{align*}
 und
 \begin{align*}
  \Prob(X_k > c(k)) \leq (2j_{max} +1)\frac{\exp(-(k-1) c(k))}{(k-1)!} (1 - \exp(-\exp(-c(k)))).
 \end{align*}

\end{prop}

\begin{bew}
Gem"a"s obigen Ausf"uhrungen haben wir nur den Fehler zu betrachten, der durch die Nichtber"ucksichtigung der Poisson"=Punktprozesse $\Pi^{(j)}$ mit $j \in \Z \setminus \{j_{min}, \ldots, j_{max} \}$ entsteht.
Wir ben"otigen also
\begin{align*}
&\Prob\Big( \exists (i,j) \in \N \times \Z_{> j_{max}}: \\
& \hspace{1.4cm} c(k) < X_i^{(j)} + \sup_{t \in [-(p j - a), -(p j - b)]\cap p\Z} \left( W_i^{(j)}(t) - \frac {|t|} 2 \right) \leq X_i^{(j)} + W_i^{(j)}(0) \Big)\\
={} &\Prob\Big( \exists i \in \N, \ j \in \Z_{> j_{max}}: \ c(k) < X_i^{(j)} + \sup_{t \in [p j - b, p j +b]\cap p\Z} \left( W_i^{(j)}(t) - \frac t 2 \right) \leq X_i^{(j)} \Big) \displaybreak[0]\\
 \leq{} & \frac {8e^{-c(k)}} p \Bigg( \exp\left(-\frac {\sqrt{pj_{max}-b}} 2 \right) + \frac {100} {9} \sqrt{\frac{b^3}{\pi}} \\
 & \hspace{0.6cm} \cdot \exp\left(-\frac{(\frac 3 {10} \sqrt{pj_{max}-b} + \frac 1 5 \sqrt{(pj_{max}-b)\vee25} -\log(\frac{\sqrt{(pj_{max}-b)\vee 25}} 2))^2} {4b}\right)\\
 & \hspace{1cm} + \frac 1 {\sqrt{2\pi}} \left(1 + \frac 1 {\sqrt{pj_{max}-b}-2}\right) \exp\left(- \frac{(\sqrt{pj_{max}-b}-2)^2} 8\right)\Bigg),
\end{align*}
wobei die letzte Ungleichung aus Lemma \ref{BBsup} folgt. Analog erh"alt man
\begin{align*}
&\Prob\Big( \exists (i,j) \in \N \times \Z_{< j_{min}}: \\
& \hspace{1.6cm} c(k) < X_i^{(j)} + \sup_{t \in [-b-pj, b-pj]\cap p\Z} \left( W_i^{(j)}(t) - \frac t 2 \right) \leq X_i^{j} + W_i^{j}(0) \Big)\displaybreak[0]\\
\leq{} &  \frac {8e^{-c(k)}} p \Bigg( \exp\left(-\frac {\sqrt{-b-pj_{min}}} 2 \right) + \frac{100}{9} \sqrt{\frac{b^3}{\pi}}\\
&\hspace{0.6cm} \cdot \exp\left(-\frac{(\frac 3 {10} \sqrt{-b-pj_{min}} + \frac 1 5 \sqrt{(-b-pj_{min})\vee25} -\log(\frac{\sqrt{(-b-pj_{min})\vee 25}} 2))^2} {4b}\right)\\
 & \hspace{1cm} + \frac 1 {\sqrt{2\pi}} \left(1 + \frac 1 {\sqrt{-b-pj_{min}}-2}\right) \exp\left(- \frac{(\sqrt{-b-pj_{min}}-2)^2} 8\right)\Bigg).
\end{align*}

Die Aussage der Proposition folgt daraus mit
\begin{align*}
& \Prob(Z_3^{(k; p)}(t) \neq \tilde{Z}_3^{(p)}(t) \text{ f"ur ein } t \in [-b,b]\cap p\Z \ | X_k \leq c(k), C_k > c(k)) \\
\leq{}& \Prob\Big( \exists (i,j) \in \N \times \Z \setminus [j_{min}, j_{max}]: \\
& \hspace{1.2cm} c(k) < X_i^{(j)} + \sup_{t \in [-(p j + b), -(p j - b)]\cap p\Z} \left( W_i^{(j)}(t) - \frac {|t|} 2 \right) \leq X_i^{(j)} + W_i^{(j)}(0) \Big).
\end{align*}

Weiterhin definieren wir Zufallsvariablen $Y_i, \ i\in \N$ durch die Poisson-Punktprozesse
$$\sum_{i \in \N} \delta_{Y_i} = \sum_{j=-b/p}^{b/p} \sum_{i \in \N} \delta_{X_i^{(j)}} \mathbf{1}_{T_i^{(j;p)} = 0}.$$
Dieser hat das Intensit\"atsma"s $(2 b/g +1) \Prob(T^{(p)} = 0) e^{-x} d x$.
Da jede Folge $(X_i^{(j)})_{i \in\N}$ monoton fallend ist, k\"onnen wir $Y_1 > Y_2 > \ldots$ annehmen.
Damit erhalten wir
\begin{align*}
& \Prob(C_k > c(k))\\
 \geq{} & \Prob(Y_1 > c(k)/2) \cdot \Prob\left( \inf_{t \in [-2b,2b]\cap p\Z} \left(W(t) - \frac{|t|} 2\right) > \frac{c(k)} 2 \ | \ T^{(p)} =0 \right)\\
 ={} & \left(1-\exp\left(-(2b/p +1) e^{-c(k)/2} \Prob(T^{(p)} = 0)\right)\right) \\
    &  \cdot\left(1 - \Prob\left( \inf_{t \in [0,2b]} (W(t)-|t|/2) < c(k)/2 \right)^2 / \Prob(T^{(p)} = 0) \right) \displaybreak[0]\\
 \geq{} & \left(1-\exp\left(-\frac{2b/p+1}{4} \cdot \exp\left(-\frac{c(k)}2\right) \left(1-\exp\left(-\frac p 2\right)\right)^2\right)\right)\\
& \cdot \Bigg( 1-  \frac{4}{1-e^{-p/2}} \Bigg( 1 - \Phi\left(-\frac{c(k)}{\sqrt{8b}} - \sqrt{\frac b 2}\right)\\
 & \hspace{2.5cm} - \exp(-c(k)/2) \left(1-\Phi\left(-\frac{c(k)}{\sqrt{8b}} + \sqrt{\frac b 2} \right)\right) \Bigg)^2\Bigg)
\end{align*}
nach Formel 2.1.1.4 in \cite{borodin-2002}.

F\"ur die Verteilung von $X_k$ ergibt sich
$$ \Prob(X_k > c(k)) \leq \sum_{j=-j_{max}}^{j_{max}} \Prob(X_k^{(j)} > c(k)) \leq (2 j_{max} +1) \Prob(X_k^{(1)} > c(k)).$$
Daraus l\"asst sich wie in Proposition \ref{fehler_0} die behauptete Absch\"atzung herleiten.
\end{bew}
\bigskip

Setzen wir sogar $p\cdot j_{min} < a-16, \ b+16 < p\cdot j_{max}$ voraus, so erhalten wir analog mit Hilfe von Lemma \ref{BBsup2} die folgende

\begin{prop} \label{fehler_3b}
Unter den genannten Voraussetzungen und Notationen gilt:
\begin{align*}
& \Prob(Z_3^{(k; p)}(t) \neq \tilde{Z}_3^{(p)}(t) \text{ f"ur ein } t \in [-b,b]\cap p\Z \ | X_k \leq c(k), C_k > c(k)) \\
\leq {} & \frac {4e^{-c(k)}} p \Bigg( \exp\left(- {\sqrt{p j_{max} - b}} \right)\\
  &+ \frac {100} {9} \sqrt{\frac{b^3}{\pi}} \exp\left(-\frac{(\sqrt{p j_{max} - b} - \log(\sqrt{p j_{max} - b}))^2} {4b}\right)\\
  &  + \frac 8 {\sqrt{2\pi}} \left(1 + \frac 2 {\sqrt{p j_{max} - b}-4}\right) \exp\left(- \frac{(\sqrt{p j_{max} - b}-4)^2} 8\right) \Bigg)
\end{align*}
\end{prop}

Diese Fehlerabsch"atzungen weisen f"ur $b \to \infty$ ein deutlich langsameres Wachstum auf als die Absch"azungen in den vorhergehenden Abschnitten.
Daher k"onnte es f"ur gro"se Intervalle sinnvoll sein, $Z_3(\cdot)$ anstelle von $Z(\cdot)$ zu approximieren.
Um Konvergenz einzusehen, m\"ussen wir wiederum $c(k)$ und $j_{max} = j_{max}(k)$ so w\"ahlen, dass
jede Fehlerabsch\"atzung aus Proposition \ref{fehler_3} f\"ur $k \to \infty$ gegen Null geht, z.\,B.
 $c(k) = - \log(k)/2$ und $j_{max} = k$.
\bigskip

\subsection{Approximation von $Z_{4,\lambda*}(\cdot)$}

Schlie"slich soll der Brown-Resnick-Prozess auch gem"a"s Satz \ref{z4} approximiert werden.
Sei zun"achst $I \subset \R^d$ ein Intervall. Wir erzeugen nun unabh"angig identisch mit Parameter $\lambda^* |I|$
 exponentialverteilte Zufallsvariablen $V_i,\ i \in \N$ und definieren $U_k:= \log\left(\frac 1 {\sum_{i=1}^n V_i}\right)$
 sowie unabh"angig identisch auf $I$ gleichverteilte Zufallsvariablen $S_i,\ i \in \N$.
 Dann ist nach Korollar \ref{PPP-dt-e-x} durch $\sum_{i \in \N} \delta_{(S_i, U_i)}$
 ein Poisson"=Punktprozess auf $\R \times I$ mit Intensit"atsma"s $\lambda^* e^{-u} du \ ds$ gegeben.
 Au"serdem seien $\tilde F_1(\cdot), \tilde F_2(\cdot), \ldots$ unabh"angig identisch verteilt mit Verteilung $Q$.
\medskip

F"ur $k \in \N$ definieren wir nun
\[ Z_{4,\lambda^*}^{(k,I)}(t) := \max_{i=1,\ldots,k} X_i + \tilde F_i(t-S_i). \]
Nach Konstruktion gilt wiederum $\lim_{\substack{k \to \infty\\ I \to \R^d}} Z_{4,\lambda^*}^{(k,I)} (t) = Z_{4,\lambda^*}(t)$ f"ur alle $t \in \R^d$.
\medskip

F"ur eine Fehlerabsch"atzung beschr"anken wir uns erneut auf den Fall,
 dass $W(\cdot)$ die eindimensionale Standard"=Brownsche"=Bewegung ist und betrachten den Prozess auf einem Gitter $p\Z, \ p>0$.
 Wir wollen $Z_{4,\lambda^{(p)}}$ auf einem Intervall $[-b,b]$ approximieren. Es sei $I=:[\-v,v] \supset [-b,b], \ v \in p\Z$. 

Seien nun $\sum_{i \in \N} \delta_{(R_i,Y_i^{(p)})}$ ein Poisson-Punktprozess auf $([-v,v]\cap p\Z) \times \R$ und
$\sum_{i \in \N} \delta_{(r_i,y_i^{(p)})}$ ein Poisson-Punktprozess auf $(p\Z \setminus [-v,v]) \times \R$
jeweils mit Intensit\"atsma"s $p\lambda^{(p)} \delta_{p\Z}(dr) e^{-u} du$ und fallend in der zweiten Komponente, d.\,h.
$Y_1^{(p)} > Y_2^{(p)} > \ldots$ sowie $y_1^{(p)} > y_2^{(p)} > \ldots$.
 Weiterhin seien $\tilde F_i, \tilde f_i \sim_{u.i.v.} Q$.
Dann gilt
\begin{equation} \label{gitter-repres}
Z_{4,\lambda^{(p)}}(\cdot) = \max\left\{ \max_{i \in \N} (Y_i^{(p)} + \tilde F_i(\cdot-R_i)),
\max_{i \in \N} (y_i^{(p)} + \tilde f_i(\cdot-t_i))\right\} .
\end{equation}

Nun definieren wir die diskrete Approximation
\begin{align*}
Z_{4,\lambda^*}^{(k,v,p)}(t) :={}& \max_{i=1,\ldots,k} \left( Y_i^{(p)} + \tilde F_i(t-R_i) \right)
\end{align*}
f"ur $t \in p\Z$, f"ur welche wir eine Fehlerabsch"atzung vornehmen werden.

Weiterhin sei $C_k:= \min_{t \in [-b,b]\cap p\Z} Z_{4,\lambda^{(p)}}^{(k;v;p)}(t)$.
F\"ur die folgende Proposition verwenden wir, dass $Q^{(p)}$ die Verteilung von $W(t) - \frac{|t|}{2} \ | \ T^{(p)} =0$ ist,
wobei $W$ eine Standard Brownsche Bewegung und $T^{(p)} = {\rm inf}\left( \argsup_{t \in p\Z}\left( W(t) - |t|/2 \right) \right)$
sei.

\begin{prop} \label{fehler_4}
Seien $c(k) <0$ und $v > b+4$. Dann gilt:
\begin{align*}
& \Prob(Z_{4,\lambda^{(p)}}^{(k;v;p)}(t) \neq Z_{4,\lambda^{(p)}}^{(p)}(t) \text{ f"ur ein } t \in [-b,b] \cap p\Z \ | \ Y_k^{(p)} \leq c(k), C_k > c(k)) \\
\leq{}& 64 e^{-c(k)} (1 - e^{- \frac p 2})^{-2} \lambda^* \\
&\cdot \Bigg( \exp\left(-\frac {\sqrt{v-b}} 2 \right) + \frac 1 {\sqrt{2\pi}} \left(1 + \frac 1 {\sqrt{v-b}-2}\right) \exp\left(- \frac{(\sqrt{v-b}-2)^2} 8\right) \\
 & + \frac {100\sqrt{b^3}} {9\sqrt{\pi}} \exp\left(-\frac{(\frac 3 {10} \sqrt{v-b} + \frac 1 5 \sqrt{(v-b)\vee25} -\log(\frac{\sqrt{(v-b)\vee 25}} 2))^2} {4b}\right)\Bigg).
\end{align*}
Weiterhin sind
\begin{align*}
 & \Prob(C_k \leq c(k))\\
 \leq{} & 1 - \left(1-\exp\left(- \lambda^{p} (2b+p) e^{-c(k)/2}\right)\right) \\
       & \cdot \left( 1 - (2(o-b)\lambda^{(p)})^{k}\frac{\exp(- \frac{(k-1)c(k)}2)}{(k-1)!}
          (1 - \exp(-e^{-(v-b)\lambda^{(p)} c(k)})) \right) \\
& \cdot \Bigg( 1-  \frac{4}{1-e^{-p/2}} \Bigg(1 - \Phi\left(-\frac{c(k)}{\sqrt{8b}} - \sqrt{\frac b 2}\right)\\
& \hspace{3.2cm} - \exp(-c(k)/2) \left( 1 - \Phi\left(-\frac{c(k)}{\sqrt{8b}} + \sqrt{\frac b 2} \right)\right) \Bigg)^2 \Bigg)
\end{align*}
und
\begin{align*}
 \Prob(Y_k^{(p)} > c(k)) \leq ((2 v+p)\lambda^{(p)})^{k} \frac{\exp(-(k-1) c(k))}{(k-1)!} (1 - \exp(-e^{-(2v+p)\lambda^{(p)} c(k)})).
\end{align*}
\end{prop}

\begin{bew}
Es gilt wegen \eqref{gitter-repres} und $X_k < c(k) < C_k$, dass
\begin{align*}
& \Prob(Z_{4,\lambda^{(p)}}^{(k;v,p)}(t) \neq Z_{4,\lambda^{(p)}}^{(p)}(t) \text{ f"ur ein } t \in [-b,b] \cap p\Z  \ | \ Y_k^{(p)} \leq c(k), C_k > c(k)) \nonumber \\
\leq{} & \Prob \Bigg( \exists i \in \N: \sup_{t \in [-b, b]\cap p\Z} y_i^{(p)} + W_i(t-r_i) - |t-r_i|/2 > c(k) \ | T^{(p)} = 0 \Bigg)\\
\leq{} & 2 \sum_{r \in v + p\N} \lambda^{(p)} p \cdot \int_{\R} e^{-x} \Prob\Bigg(\sup_{t \in [-b, b] \cap p\Z} x + W(t-s) - \frac {|t-s|} 2 > c(k)  \Bigg| \ T^{(p)} = 0\Bigg) dx\nonumber \displaybreak[0] \\
\leq{} & 2 \frac{p\lambda^{(p)}}{\Prob(T^{(p)}=0)} \Bigg(\sum_{s \in o + p\N} \int_{\R} e^{-x} \Prob\Bigg( c(k) < \sup_{t \in [s-b, s+b] \cap p\Z} x + W(t) - \frac {|t|} 2 \leq x\Bigg) dx\nonumber.
\end{align*}

Analog zu Lemma \ref{BBsup} ergibt sich daraus
\begin{align*}
& \Prob(Z_{4,\lambda^{(p)}}^{(k;v;p)}(t) \neq Z_{4,\lambda^{(p)}}(t) \text{ f"ur ein } t \in [-b,b] \cap p\Z \ | \ X_k \leq c(k), C_k > c(k)) \displaybreak[0]\\
\leq{}& 16e^{-C} \lambda^* \Prob(T^{(p)}=0)^{-1} \\
&\cdot \Bigg( \exp\left(-\frac {\sqrt{v-b}} 2 \right) + \frac 1 {\sqrt{2\pi}} \left(1 + \frac 1 {\sqrt{v-b}-2}\right) \exp\left(- \frac{(\sqrt{v-b}-2)^2} 8\right) \\
 & + \frac {100\sqrt{b^3}} {9\sqrt{\pi}} \exp\left(-\frac{(\frac 3 {10} \sqrt{v-b} + \frac 1 5 \sqrt{(v-b)\vee25} -\log(\frac{\sqrt{(v-b)\vee 25}} 2))^2} {4b}\right)\Bigg)
 \end{align*}

Au"serdem ist
\begin{align*}
  & \Prob(T^{(p)}=0) \\
\geq{} & \Prob\Bigg( \sup_{t \in (-\infty,-p)} \left(W(t) - \frac{|t|} 2 - W(-p) + \frac {|-p|} 2\right) < \frac p 2,
 \ W(-p) - \frac {|-p|} 2 \leq - \frac p 2, \\
 & \hspace{1cm} W(p) - \frac {|p|} 2 \leq - \frac p 2, \ \sup_{t \in (p,\infty)} \left(W(t) - \frac{|t|} 2 - W(p) + \frac {|p|} 2\right) < \frac p 2\Bigg) \\
 ={} & \left[ \Prob(W(p) \leq 0) \cdot \Prob\left( \sup_{t \in (0,\infty)} W(t) - \frac t 2 < \frac p 2\right)\right]^2 \\
 ={} & \left[ \frac 1 2 \cdot (1-e^{- \frac p 2})\right]^2 \quad \text{nach Formel 2.1.1.4 (1) aus \cite{borodin-2002}} \\
 ={} & \frac 1 4 (1 - e^{- \frac p 2})^2.
\end{align*}
Damit folgt die erste Aussage der Proposition.

Weiterhin gilt
\begin{align*}
  \Prob(Y_k^{(p)} > x) ={} & \Prob(e^{-Y_k^{(p)}} < e^{-x})\\
                    ={} & \int_0^{e^{-x}} ((2v/p+1)p\lambda^{(p)})^{k} t^{k-1}  \frac{\exp(-(2v/p+1)p\lambda^{(p)}t)}{(k-1)!} d t\displaybreak[0]\\
                    \leq{} & ((2v+p)\lambda^{(p)})^{k} \frac{\exp(-(k-1) x(k))}{(k-1)!}
                                    \int_0^{e^{-x}} \exp(-(2v+p)\lambda^{(p)}t) d t\\
                    ={} & ((2v+p)\lambda^{(p)})^{k}\frac{\exp(-(k-1) x)}{(k-1)!} (1 - \exp(-\exp(-(2v+p)\lambda^{(p)}x))).
\end{align*}

\"Ahnlich wie im Beweis von Proposition \ref{fehler_3}
definieren wir nun Zufallsvariablen $Y_i', Y_i'',\ i \in \N,$ \"uber die Poisson-Punktprozesse $\sum_{i \in \N} \delta_{Y_i'} = \sum_{i\in\N} \delta_{Y_i^{(p)}} \mathbf{1}_{|s_i| \leq b}$
und $\sum_{i \in \N} \delta_{Y_i''} = \sum_{i\in\N} \delta_{Y_i^{(o;p)}} \mathbf{1}_{|s_i| > b}$ mit $Y_1' \geq Y_2' \geq \ldots$ and $Y_1'' \geq Y_2'' \geq \ldots$.

Dann gilt
\begin{align*}
 & \Prob(C_k > c(k))\\
 \geq{} & \Prob(Y_1' > c(k)/2) \cdot \Prob(Y_k'' < c(k) /2) \\
 &\cdot \Prob\left( \inf_{t \in [-2b,2b]\cap p\Z} \left(W(t) - \frac{|t|} 2\right) > \frac{c(k)} 2 \ | \ T^{(p)} =0 \right) \displaybreak[0] \\
 ={} & \left(1-\exp\left(- \lambda^{(p)} p (2b/p+1) e^{-c(k)/2}\right)\right) \\
       & \cdot \left( 1 - (2(v-b)\lambda^{(p)})^{k}\frac{\exp(-\frac{(k-1) c(k)} 2)}{(k-1)!}
          (1 - \exp(-e^{-(v-b)\lambda^{(p)} c(k)})) \right) \\
    &  \cdot\left(1 - \Prob\left( \inf_{t \in [0,2b]} (W(t)-|t|/2) < c(k)/2 \right)^2 / \Prob(T^{(p)} = 0) \right) \displaybreak[0]\\
 \geq{} & \left(1-\exp\left(-\lambda^{(p)} (2b+p) e^{-c(k)/2}\right)\right) \\
       & \cdot \left( 1 - (2(v-b)\lambda^{(p)})^{k}\frac{\exp(-\frac{(k-1) c(k)}2)}{(k-1)!}
          (1 - \exp(-e^{-(v-b)\lambda^{(p)} c(k)})) \right) \\
& \cdot \Bigg( 1-  \frac{4}{1-e^{-p/2}} \Bigg(1 - \Phi\left(-\frac{c(k)}{\sqrt{8b}} - \sqrt{\frac b 2}\right)\\
& \hspace{3.2cm} - \exp(-c(k)/2) \left( 1 - \Phi\left(-\frac{c(k)}{\sqrt{8b}} + \sqrt{\frac b 2} \right)\right) \Bigg)^2 \Bigg).
\end{align*} 
\end{bew}

Versch"arfen wir die Voraussetzungen mit $ b+16 < o$, so erhalten wir mit einer Rechnung wie im Beweis von Lemma \ref{BBsup2} die folgende

\begin{prop} \label{fehler_4b}
Unter den genannten Voraussetzungen und Notationen gilt:
\begin{align*}
& \Prob(Z_{4,\lambda^{(p)}}^{(k;v;p)}(t) \neq Z_{4,\lambda^{(p)}}^{(p)}(t) \text{ f"ur ein } t \in [-b,b] \cap p\Z \ | \ X_k \leq c(k), C_k > c(k)) \\
\leq{} & 16 e^{-c(k)} (1 - e^{- \frac p 2})^{-2} \lambda^{(p)} \\
&\cdot \Bigg( \exp\left(-\sqrt{v-b}  \right) + \frac 8 {\sqrt{2\pi}} \left(1 + \frac 1 {\sqrt{v-b}-4}\right) \exp\left(- \frac{(\sqrt{v-b}-4)^2} 8\right) \\
 & + \frac {100\sqrt{b^3}} {9\sqrt{\pi}} \exp\left(-\frac{(\sqrt{o-b} -\log(\sqrt{v-b}))^2} {4b}\right) \Bigg).
\end{align*}
\end{prop}

Diese Absch"atzungen haben f"ur $b \to \infty$ das gleiche Verhalten wie im vorherigen Abschnitt.
Daher k"onnte auch die Approximation von $Z_{4, \lambda}$ besonders f"ur gro"se Intervalle sehr interessant sein.
W\"ahlen wir beispielsweise $v = v(k) = k^{1/4}$ und $c(k) = -\log(k)/4$, konvergieren wiederum alle Absch\"atzungen aus Proposition
\ref{fehler_4} gegen $0$ f\"ur $k \to \infty$.

%% file: Simulation_revised.tex
\chapter{Simulationsergebnisse}

F\"ur eine aktualisierte und ausf\"uhrliche Simulationsstudie, die insbesondere die Aussagen der S\"atze \ref{thm:conditional1} und \ref{thm:conditional2} ber\"ucksichtigt
sowie die Methoden anhand mathematisch sch\"arfere Kriterien vergleicht, sei auf

\textsc{Oesting}, Marco; \textsc{Kabluchko}, Zakhar ; \textsc{Schlather}, Martin :
\newblock \emph{Simulation Techniques for Brown Resnick Processes}.
\newblock Submitted to \emph{Extremes} (2010).

verwiesen.
Ein besonderer Dank gilt dabei den Referees dieses Papers, die mit ihren Hinweisen zu einer Verbesserung der Fehlersbsch\"atzungen und der Simulationsstudie beigetragen haben.